\documentclass[12pt,reqno]{amsart}
\usepackage{amsmath,amssymb,url,color,tikz,colortbl}
\usepackage[all]{xy}
\usepackage{young}
\usepackage{enumerate}
\usepackage[margin=1in]{geometry}
\usepackage{hyperref}
\usepackage{amsfonts}
\usepackage{amsmath}
\usepackage{amssymb}
\usepackage{amsthm}
\usepackage{xcolor}
\usepackage{graphicx}
\usepackage{listings}
\usepackage{tram}
\usepackage{manfnt}
\usepackage{url,color}

\def\lcm{\mathop{\rm LCM}}

\usepackage[backend=bibtex]{biblatex}
\usepackage{amsfonts,latexsym,tikz,tikz-cd,hyperref, wrapfig,float,wasysym}

\newtheorem{thm}{Theorem}[section]
\newtheorem{prop}[thm]{Proposition}
\newtheorem{cor}[thm]{Corollary}
\newtheorem{porism}[thm]{Porism}
\newtheorem{lem}[thm]{Lemma}

\theoremstyle{definition}
\newtheorem{remark}[thm]{Remark}
\newtheorem{definition}[thm]{Definition}
\newtheorem{ex}[thm]{Example}

\usepackage{comment}
\newif\ifshow
\showtrue
\ifshow
  
\else
  \excludecomment{longproof}
\fi

\newcommand{\htau}{\tau}

\newcommand{\VV}{{\sf V}}
\newcommand{\NN}{\mathbb{N}}

\newcommand{\hl}{ {\sf h} }
\newcommand{\tl}{ {\sf t} }

\newcommand{\BB}{ {\sf C} }

\newcommand{\jrho}{\rho}

\newcommand{\calF}{\mathcal{F}}

\newcommand{\return}[1]{%
  \hspace{0.2cm}
  \begin{tikzpicture}[baseline={([yshift=-.8ex]current bounding box.center)}]
    \node[font = {\scriptsize}] at (0,0.2) {$#1$};
    \draw (-0.2,0) -- (0.2,0) -- (0.2,-0.3);
    \draw (0.1,-0.2) -- (0.2,-0.3) -- (0.3,-0.2);
  \end{tikzpicture}%
}

\usepackage{mathtools}

\newcommand{\tikzColoredSquare}[4]{%
\path [fill=#4] (#1*0.5-0.25,#2*0.5-0.25) rectangle (#1*0.5+0.25,#2*0.5+0.25);
\node at (#1*0.5,#2*0.5) {$#3$};
}

\usetikzlibrary{patterns}

\newcommand{\tikzColoredSquarenwl}[4]{%
\path [pattern=north west lines, pattern color=#4] (#1*0.5-0.25,#2*0.5-0.25) rectangle (#1*0.5+0.25,#2*0.5+0.25);
\node at (#1*0.5,#2*0.5) {$#3$};
}

\newcommand{\tikzColoredSquarehl}[4]{%
\path [pattern=horizontal lines, pattern color=#4] (#1*0.5-0.25,#2*0.5-0.25) rectangle (#1*0.5+0.25,#2*0.5+0.25);
\node at (#1*0.5,#2*0.5) {$#3$};
}

\newcommand{\tikzColoredSquarech}[4]{%
\path [pattern=crosshatch, pattern color=#4] (#1*0.5-0.25,#2*0.5-0.25) rectangle (#1*0.5+0.25,#2*0.5+0.25);
\node at (#1*0.5,#2*0.5) {$#3$};
}

\usepackage[colorinlistoftodos,color=red!70]{todonotes}

\begin{document}

\title{Rowmotion on the chain of V's poset and whirling dynamics}
\author[Plante]{Matthew Plante}
\address{Department of Mathematics, University of Connecticut, Storrs, CT 06269-1009, USA}
\email{matthew.plante@uconn.edu}

\author[Roby]{Tom Roby}
\address{Department of Mathematics, University of Connecticut, Storrs, CT 06269-1009, USA}
\email{tom.roby@uconn.edu}
\date{\today}

\begin{abstract}
Given a finite poset $P$, we study the \emph{whirling} action on vertex-labelings of $P$
with the elements $\{0,1,2,\dotsc ,k\}$.  When such labelings are (weakly) order-reversing,
we call them \textit{$k$-bounded $P$-partitions}.  
We give a general equivariant bijection between $k$-bounded $P$-partitions and order
ideals of the poset $P\times [k]$ which conveys whirling to the well-studied \textit{rowmotion} operator.  As an
application, we derive periodicity and homomesy results for rowmotion acting on
the \textit{chain of V's} poset $\VV \times [k]$.  We are able to generalize some of these results to
the more complicated dynamics of rowmotion on $\BB_{n}\times [k]$, where $\BB_{n}$ is the
\textbf{claw poset} with $n$ unrelated elements each covering $\widehat{0}$. \vspace{1 ex}

\end{abstract}

\maketitle
\vspace{-0.333 in}
\tableofcontents
\vspace{-0.333 in}
\section{Introduction}\label{sec:intro}

We connect the well-studied operation of \textit{rowmotion} on the order ideals of a finite poset
with the less familiar \textit{whirling} action on $P$-partitions with bounded labels.  One
of our main results is an equivariant bijection that carries one to the other for any finite
poset $P$.  Here we are following in the footsteps of Haddadan~\cite{shahrzad}, who noticed
this connection when $P$ itself is a chain.  We then leverage this to study the rowmotion
action on the ``chain of V's'' poset $\VV_{k}:=\VV \times [k]$ 
(a 3-element V-shaped poset cross a finite chain, see Figure~\ref{fig:vee}), which has surprisingly good dynamical properties.
We also generalize this to the case where we replace $\VV$ with a \textbf{$n$-claw}, $\BB_{n}$, a poset with a single minimal
element covered by exactly $n$ incomparable elements. In both cases we obtain both 
periodicity results and homomesy.  

\subsection{Basic definitions and background in combinatorial dynamics }\label{}

Let $P$ be a finite poset, and $\mathcal{J}\left( P\right)$ be the set of order ideals of
$P$. (For basic poset definitions, we refer the reader to Stanley~\cite[Ch.~3]{Sta11}.)
\emph{Combinatorial rowmotion} is an invertible map
$\rho:\mathcal{J}(P)\rightarrow\mathcal{J}(P)$ which takes each ideal $I\in\mathcal{J}(P)$
to the order ideal generated by the minimal elements of the complement of $I$ in $P$. The
periodicity of this map on products of chains was first studied by Brouwer and
Schrijver~\cite{BS74}, and Cameron and Fon-der-Flaass~\cite{CF95}.  Later Striker and
Williams~\cite{SW12} considered it as one element of the ``toggle group'' of a poset and
related it to a kind of ``promotion'' operator on order ideals.  Around the same time,
Armstrong, Stump, and Thomas~\cite{ast} studied rowmotion on \emph{root posets}, relating it
to ``Kreweras complementation'' on noncrossing partitions, and used this to prove a
conjecture of Panyushev about the equality of the average cardinality of antichains for each
rowmotion orbit.  

Propp and Roby~\cite{PR13} noticed that this conjecture was merely one
instance of a much broader phenomenon which they dubbed \emph{homomesy}. Given
a finite set $S$, a \textquotedblleft statistic\textquotedblright%
\ $f:S\rightarrow\mathbb{C}$, and an invertible map $\varphi$ on $S$, we call
$f$ \emph{homomesic} if the average value of $f$ is the same for every
$\varphi$-orbit $\mathcal{R}$, i.e., $\displaystyle\frac{1}{\#\mathcal{R}}%
\sum_{x\in\mathcal{R}}f(x)=c$, where $c$ is a constant not dependent on the
choice of orbit $\mathcal{R}$. Over the past decade, many authors have proved homomesy
results as one tool to understand various combinatorial dynamical systems.

Along with the work initiated by Reiner, Stanton, and White on the \emph{Cyclic Sieving
Phenomenon}~\cite[]{csp,csp-brief}, the confluence of all this work was the
beginning of \emph{dynamical algebraic combinatorics} as a distinct area
within algebraic combinatorics (with antecedents going back to the
Robinson--Schensted--Knuth correspondence and related operations on Young
tableaux such as promotion, evacuation, and cyclage).  In the past decade, the subfield has
grown in a number of directions, and the study of rowmotion has been of continuing interest.
For more background information, see the survey articles of Hopkins~\cite{Hopkins-order},
Roby~\cite{robydac}, and Striker~\cite{striker-dac}.

Cameron and Fon-der-Flaass~\cite{CF95} were the first to describe rowmotion as a product of
involutions called \textit{toggles}, as detailed in Section~\ref{ss:toggles}.  A natural
generalization of toggling at a poset element $x$ is ``\textit{whirling} at $x$,'' which cycles
the label at $x$ among $j$ possible values.  (Toggles are the case when $j=2$.)  Joseph,
Propp, and Roby defined these and the operation of \textit{whirling} on sets of functions
between finite sets, obtaining various homomesy results for different classes of functions
(injective, surjective, etc.)~\cite{JPR18}. This is described in
Section~\ref{sec:whirling}.  

A bijective function $f:P \to [p]$ (with $\# P = p$) such that $f(x)<f(y)$ 
whenever %%%if and only if (WRONG)
$x<_P y$ is called a \emph{linear extension}.  We denote by $\mathcal{L}(P)$ the set of
\emph{all linear extensions}
of $P$; its cardinality, $e(P)$, is an important numerical invariant of a poset. 
Its refinement, the \emph{order polynomial} $\Omega_{P}(k)$, counts the number of $k$-bounded
$P$-partitions.  For some special posets $P$, mainly ones connected with Lie theory (root
and minuscule posets) and those of partition or shifted shapes, product formulae for
$\Omega_{P}(k)$ are known.  Hopkins surveys these posets, the formulae, and gives the
heuristic: \textit{Posets with order-polynomial product formulae are the same as the posets
with good dynamical behavior.}   The one poset in his list whose rowmotion dynamics were relatively
unexplored is $\VV \times [k]$, a gap this paper fills.  In separate work Hopkins and Rubey
study the dynamics of Sch\"utzenberger promotion on linear extensions of $\VV \times [k]$,
which also exhibits unusually good behavior~\cite{Hopkins-Rubey}.  In his doctoral thesis,
Ben Adenbaum provides a combinatorial approach (along with a density argument) to proving
that this same periodicity extends to \textit{piecewise-linear rowmotion} on the order
polytope $\mathcal{O}(P)$~\cite[Prop.~3.23]{a24}.  Since most dynamical results in the
piecewise-linear realm have been proven first at the (higher) birational level, then
tropicalized to the PL realm, this is a notable methodology.  The authors are unaware of any
successful attempts to lift periodicity of the rowmotion action on $\VV \times [k]$
(even conjecturally) to the birational realm.

\subsection{Organization of this paper}\label{ss:org}
In Section~\ref{sec:intro} after the introduction, we review the
toggling definition of rowmotion. Section~\ref{sec:whirling} describes whirling, and includes the equivariant
bijection which allows us to study rowmotion on $\VV_{k}$ as whirling on $k$-bounded
$P$-partitions.  Section~\ref{sec:vee} contains our main periodicity and homomesy results for rowmotion
on $\VV_{k}$, which use decompositions of the ``orbit board'' of the corresponding whirling
action into ``whorms''.  Finally, Section~\ref{sec:claw} contains the periodicity and homomesy results
which generalize to rowmotion on the ``chain of claws'' graph, $\BB_{n}\times [k]$, and
extends our proofs to this situation. 

\subsection{Rowmotion as a product of toggles}\label{ss:toggles} 
Writing rowmotion as a product of toggles has been quite useful for understanding its
properties and generalizing it to other situations, e.g., lifting to the piece-wise linear
and birational level~\cite[]{EP20}.  Thomas and Williams call this ``Rowmotion in
Slowmotion''~\cite{TW19}.  
\begin{definition}\label{def:rowmotion}
We define the (order-ideal) rowmotion map, $\rho :\mathcal{J}(P)\to\mathcal{J}(P)$ as
follows:  For any $I\in \mathcal{J}(P)$, $\rho (I)$ is
the order ideal generated by the minimal elements of the complement of $I$, as in the example below.
\end{definition}

\begin{ex}
Here is one iteration of $\jrho$ on an order ideal with the action broken down
into its three steps: (1) complement, (2) take minimal elements, (3) saturate down. 
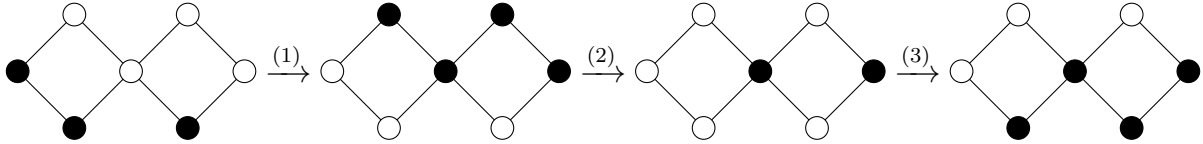
\[
\begin{tikzpicture}[xscale = 0.75, yscale = 0.75, baseline={([yshift=-.8ex]current bounding box.center)}]
    \draw[-] (0,0) -- (1,-1);
    \draw[-] (0,0) -- (1,1);
    \draw[-] (2,0) -- (1,-1);
    \draw[-] (2,0) -- (1,1);
    \draw[-] (2,0) -- (3,-1);
    \draw[-] (2,0) -- (3,1);
    \draw[-] (4,0) -- (3,-1);
    \draw[-] (4,0) -- (3,1);
    
    \draw [fill] (0,0) circle [radius =.2];
    \draw [fill] (1,-1) circle [radius =.2];
    \draw [fill=white] (1,1) circle [radius =.2];
    \draw [fill=white] (2,0) circle [radius =.2];
    \draw [fill=white] (3,1) circle [radius =.2];
    \draw [fill] (3,-1) circle [radius =.2];
    \draw [fill=white] (4,0) circle [radius =.2];
    
\end{tikzpicture}
\xrightarrow[]{(1)}
\begin{tikzpicture}[xscale = 0.75, yscale = 0.75, baseline={([yshift=-.8ex]current bounding box.center)}]
    \draw[-] (0,0) -- (1,-1);
    \draw[-] (0,0) -- (1,1);
    \draw[-] (2,0) -- (1,-1);
    \draw[-] (2,0) -- (1,1);
    \draw[-] (2,0) -- (3,-1);
    \draw[-] (2,0) -- (3,1);
    \draw[-] (4,0) -- (3,-1);
    \draw[-] (4,0) -- (3,1);
    
    \draw [fill=white] (0,0) circle [radius =.2];
    \draw [fill=white] (1,-1) circle [radius =.2];
    \draw [fill] (1,1) circle [radius =.2];
    \draw [fill] (2,0) circle [radius =.2];
    \draw [fill] (3,1) circle [radius =.2];
    \draw [fill=white] (3,-1) circle [radius =.2];
    \draw [fill] (4,0) circle [radius =.2];
    
\end{tikzpicture}
\xrightarrow[]{(2)}
\begin{tikzpicture}[xscale = 0.75, yscale = 0.75, baseline={([yshift=-.8ex]current bounding box.center)}]
    \draw[-] (0,0) -- (1,-1);
    \draw[-] (0,0) -- (1,1);
    \draw[-] (2,0) -- (1,-1);
    \draw[-] (2,0) -- (1,1);
    \draw[-] (2,0) -- (3,-1);
    \draw[-] (2,0) -- (3,1);
    \draw[-] (4,0) -- (3,-1);
    \draw[-] (4,0) -- (3,1);
    
    \draw [fill=white] (0,0) circle [radius =.2];
    \draw [fill=white] (1,-1) circle [radius =.2];
    \draw [fill=white] (1,1) circle [radius =.2];
    \draw [fill] (2,0) circle [radius =.2];
    \draw [fill=white] (3,1) circle [radius =.2];
    \draw [fill=white] (3,-1) circle [radius =.2];
    \draw [fill] (4,0) circle [radius =.2];
    
\end{tikzpicture}
\xrightarrow[]{(3)}
\begin{tikzpicture}[xscale = 0.75, yscale = 0.75, baseline={([yshift=-.8ex]current bounding box.center)}]
    \draw[-] (0,0) -- (1,-1);
    \draw[-] (0,0) -- (1,1);
    \draw[-] (2,0) -- (1,-1);
    \draw[-] (2,0) -- (1,1);
    \draw[-] (2,0) -- (3,-1);
    \draw[-] (2,0) -- (3,1);
    \draw[-] (4,0) -- (3,-1);
    \draw[-] (4,0) -- (3,1);
    
    \draw [fill=white] (0,0) circle [radius =.2];
    \draw [fill] (1,-1) circle [radius =.2];
    \draw [fill=white] (1,1) circle [radius =.2];
    \draw [fill] (2,0) circle [radius =.2];
    \draw [fill=white] (3,1) circle [radius =.2];
    \draw [fill] (3,-1) circle [radius =.2];
    \draw [fill] (4,0) circle [radius =.2];
    
\end{tikzpicture}
\]
\end{ex}

Cameron and Fon-der-Flaass~\cite{CF95} showed that for any finite poset $P$, rowmotion can
be realized as ``toggling once at each element of $P$ along any linear extension (from top
to bottom)''.  Other toggling orders also lead to interesting maps, such as
Striker--Williams ``promotion'' (of order ideals) of a poset, which is toggling from
left-to-right along ``files'' of a poset~\cite{SW12}. To \textit{toggle} an order ideal $I$
at a poset element $x$, means to add $x$ to $I$ if it is not included, or delete $x$ if it is,
but only if the result would also be an order ideal.  Formally: 

\begin{definition}\label{def:toggles}
For each fixed $x\in P$ define the \emph{(order-ideal) toggle}  $\htau_x:\mathcal{J}(P)\to\mathcal{J}(P)$ by
\[
\htau_x(I) = 
\begin{cases}
I\smallsetminus \{x\}&\text{if }x\in I\text{ and }I\smallsetminus \{x\}\in \mathcal{J}(P)\\
I\cup \{x\}&\text{if }x\not\in I\text{ and }I\cup\{x\}\in\mathcal{J}(P)\\
I&\text{otherwise. }
\end{cases}
\]
\end{definition}
\noindent It is an easy exercise to show that order-ideal toggles \cite[\S 2]{CF95} are involutions,
and that toggles at incomparable elements commute (a special case of
Prop~\ref{prop:inccommute}).   

\begin{ex}
We will toggle node-by-node down the following fixed linear extension: at each step we consider whether or not to toggle the red node in or out. 
$
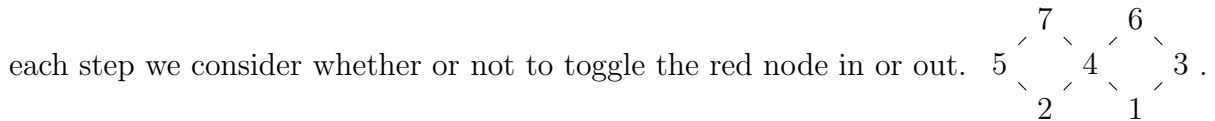
\begin{tikzpicture}[scale =0.6, baseline={([yshift=-.8ex]current bounding box.center)}]

    \node at (-1,1) (2) {$2$};
    \node at (-2,2) (4) {$5$};
    \node at (0,2) (6) {$4$};
    \node at (-1,3) (12) {$7$};
    \node at (1,1) (3) {$1$};
    \node at (2,2) (9) {$3$};
    \node at (1,3) (18) {$6$};
    
    \draw[-] (2) -- (4);
    \draw[-] (2) -- (6);
    \draw[-] (3) -- (6);
    \draw[-] (3) -- (9);
    \draw[-] (4) -- (12);
    \draw[-] (6) -- (12);
    \draw[-] (6) -- (18);
    \draw[-] (9) -- (18);
\end{tikzpicture}.
$

For this linear extension we toggle the elements row-by-row from top-to-bottom,
left-to-right within each row (although the order within rows is irrelevant, since toggles at
incomparable elements commute, hence the name ``rowmotion'').    

\begin{align*}
\begin{tikzpicture}[xscale = 0.75, yscale = 0.75, baseline={([yshift=-.8ex]current bounding box.center)}]
    \draw[-] (0,0) -- (1,-1);
    \draw[-] (0,0) -- (1,1);
    \draw[-] (2,0) -- (1,-1);
    \draw[-] (2,0) -- (1,1);
    \draw[-] (2,0) -- (3,-1);
    \draw[-] (2,0) -- (3,1);
    \draw[-] (4,0) -- (3,-1);
    \draw[-] (4,0) -- (3,1);
    \draw [fill=white] (0,0) circle [radius =.2];
    \draw [fill=] (1,-1) circle [radius =.2];
    \draw [red,fill=white] (1,1) circle [radius =.2];
    \draw [fill=] (2,0) circle [radius =.2];
    \draw [fill] (3,1) circle [radius =.2];
    \draw [fill] (3,-1) circle [radius =.2];
    \draw [fill] (4,0) circle [radius =.2];
\end{tikzpicture}
\xrightarrow{\htau_7}
\begin{tikzpicture}[xscale = 0.75, yscale = 0.75, baseline={([yshift=-.8ex]current bounding box.center)}]
    \draw[-] (0,0) -- (1,-1);
    \draw[-] (0,0) -- (1,1);
    \draw[-] (2,0) -- (1,-1);
    \draw[-] (2,0) -- (1,1);
    \draw[-] (2,0) -- (3,-1);
    \draw[-] (2,0) -- (3,1);
    \draw[-] (4,0) -- (3,-1);
    \draw[-] (4,0) -- (3,1);
    \draw [fill=white] (0,0) circle [radius =.2];
    \draw [fill=] (1,-1) circle [radius =.2];
    \draw [fill=white] (1,1) circle [radius =.2];
    \draw [fill=] (2,0) circle [radius =.2];
    \draw [red,fill] (3,1) circle [radius =.2];
    \draw [fill] (3,-1) circle [radius =.2];
    \draw [fill] (4,0) circle [radius =.2];
\end{tikzpicture}
\xrightarrow{\htau_6}
\begin{tikzpicture}[xscale = 0.75, yscale = 0.75, baseline={([yshift=-.8ex]current bounding box.center)}]
    \draw[-] (0,0) -- (1,-1);
    \draw[-] (0,0) -- (1,1);
    \draw[-] (2,0) -- (1,-1);
    \draw[-] (2,0) -- (1,1);
    \draw[-] (2,0) -- (3,-1);
    \draw[-] (2,0) -- (3,1);
    \draw[-] (4,0) -- (3,-1);
    \draw[-] (4,0) -- (3,1);
    \draw [red,fill=white] (0,0) circle [radius =.2];
    \draw [fill=] (1,-1) circle [radius =.2];
    \draw [fill=white] (1,1) circle [radius =.2];
    \draw [fill] (2,0) circle [radius =.2];
    \draw [fill=white] (3,1) circle [radius =.2];
    \draw [fill] (3,-1) circle [radius =.2];
    \draw [fill] (4,0) circle [radius =.2];
\end{tikzpicture}
\xrightarrow{\htau_5}
\begin{tikzpicture}[xscale = 0.75, yscale = 0.75, baseline={([yshift=-.8ex]current bounding box.center)}]
    \draw[-] (0,0) -- (1,-1);
    \draw[-] (0,0) -- (1,1);
    \draw[-] (2,0) -- (1,-1);
    \draw[-] (2,0) -- (1,1);
    \draw[-] (2,0) -- (3,-1);
    \draw[-] (2,0) -- (3,1);
    \draw[-] (4,0) -- (3,-1);
    \draw[-] (4,0) -- (3,1);
    \draw [fill] (0,0) circle [radius =.2];
    \draw [fill] (1,-1) circle [radius =.2];
    \draw [fill=white] (1,1) circle [radius =.2];
    \draw [red,fill] (2,0) circle [radius =.2];
    \draw [fill=white] (3,1) circle [radius =.2];
    \draw [fill] (3,-1) circle [radius =.2];
    \draw [fill] (4,0) circle [radius =.2];
\end{tikzpicture}\\
\xrightarrow{\htau_4}
\begin{tikzpicture}[xscale = 0.75, yscale = 0.75, baseline={([yshift=-.8ex]current bounding box.center)}]
    \draw[-] (0,0) -- (1,-1);
    \draw[-] (0,0) -- (1,1);
    \draw[-] (2,0) -- (1,-1);
    \draw[-] (2,0) -- (1,1);
    \draw[-] (2,0) -- (3,-1);
    \draw[-] (2,0) -- (3,1);
    \draw[-] (4,0) -- (3,-1);
    \draw[-] (4,0) -- (3,1);
    \draw [fill] (0,0) circle [radius =.2];
    \draw [fill] (1,-1) circle [radius =.2];
    \draw [fill=white] (1,1) circle [radius =.2];
    \draw [fill=white] (2,0) circle [radius =.2];
    \draw [fill=white] (3,1) circle [radius =.2];
    \draw [fill] (3,-1) circle [radius =.2];
    \draw [red,fill] (4,0) circle [radius =.2];
\end{tikzpicture}
\xrightarrow{\htau_3}
\begin{tikzpicture}[xscale = 0.75, yscale = 0.75, baseline={([yshift=-.8ex]current bounding box.center)}]
    \draw[-] (0,0) -- (1,-1);
    \draw[-] (0,0) -- (1,1);
    \draw[-] (2,0) -- (1,-1);
    \draw[-] (2,0) -- (1,1);
    \draw[-] (2,0) -- (3,-1);
    \draw[-] (2,0) -- (3,1);
    \draw[-] (4,0) -- (3,-1);
    \draw[-] (4,0) -- (3,1);
    \draw [fill] (0,0) circle [radius =.2];
    \draw [red,fill] (1,-1) circle [radius =.2];
    \draw [fill=white] (1,1) circle [radius =.2];
    \draw [fill=white] (2,0) circle [radius =.2];
    \draw [fill=white] (3,1) circle [radius =.2];
    \draw [fill] (3,-1) circle [radius =.2];
    \draw [fill=white] (4,0) circle [radius =.2];
\end{tikzpicture}
\xrightarrow{\htau_2}
\begin{tikzpicture}[xscale = 0.75, yscale = 0.75, baseline={([yshift=-.8ex]current bounding box.center)}]
    \draw[-] (0,0) -- (1,-1);
    \draw[-] (0,0) -- (1,1);
    \draw[-] (2,0) -- (1,-1);
    \draw[-] (2,0) -- (1,1);
    \draw[-] (2,0) -- (3,-1);
    \draw[-] (2,0) -- (3,1);
    \draw[-] (4,0) -- (3,-1);
    \draw[-] (4,0) -- (3,1);
    \draw [fill] (0,0) circle [radius =.2];
    \draw [fill] (1,-1) circle [radius =.2];
    \draw [fill=white] (1,1) circle [radius =.2];
    \draw [fill=white] (2,0) circle [radius =.2];
    \draw [fill=white] (3,1) circle [radius =.2];
    \draw [red,fill] (3,-1) circle [radius =.2];
    \draw [fill=white] (4,0) circle [radius =.2];
\end{tikzpicture}
\xrightarrow{\htau_1}
\begin{tikzpicture}[xscale = 0.75, yscale = 0.75, baseline={([yshift=-.8ex]current bounding box.center)}]
    \draw[-] (0,0) -- (1,-1);
    \draw[-] (0,0) -- (1,1);
    \draw[-] (2,0) -- (1,-1);
    \draw[-] (2,0) -- (1,1);
    \draw[-] (2,0) -- (3,-1);
    \draw[-] (2,0) -- (3,1);
    \draw[-] (4,0) -- (3,-1);
    \draw[-] (4,0) -- (3,1);
    \draw [fill] (0,0) circle [radius =.2];
    \draw [fill] (1,-1) circle [radius =.2];
    \draw [fill=white] (1,1) circle [radius =.2];
    \draw [fill=white] (2,0) circle [radius =.2];
    \draw [fill=white] (3,1) circle [radius =.2];
    \draw [fill=white] (3,-1) circle [radius =.2];
    \draw [fill=white] (4,0) circle [radius =.2];
\end{tikzpicture}
\end{align*}

\end{ex}

\begin{prop}[\protect{\cite[Lemma~1]{CF95}}]
\label{prop:togall}
Let $x_1,x_2,\dots,x_p$ be any linear extension (i.e., any
order-preserving listing of the elements) of a finite poset $P$ with $p$
elements. Then the composite map $\htau_{x_1} \htau_{x_2} \cdots
\htau_{x_p}$ coincides with the rowmotion operation $\jrho$.  
\end{prop}

\section{The whirling map}\label{sec:whirling}

An order ideal in any poset can be considered as a binary labeling of the poset with 0
(indicating elements outside the order ideal) and 1 (indicating those inside of it).
Equivalently, this represents an order-reversing map $f:P\rightarrow \{0,1 \}$.  Then
toggling at an element, simply switches the labels at that element (or leaves them alone if
the result would not be order-reversing).  Similarly, we can define whirling at a poset
element to cycle through a larger set of possible labels until it arrives at one that gives
a legitimate order-reversing map.  To whirl the entire labeling, one whirls once at each
element along a linear extension.  

A notion equivalent to whirling in this context was first defined by James Propp (dubbed
\textit{winching} by Peter Winkler) and was used by S.~Haddadan\cite{shahrzad}, as an aid to
proving certain homomesies for rowmotion on order ideals of type $A$ root and minuscule posets (triangles and
rectangles).  Joseph, Propp, and Roby~\cite{JPR18} later defined whirling in the context of
functions between finite sets.  We review this first, before defining whirling of
``$k$-bounded $P$-partitions'' in Section~\ref{ss:whirlP}. We prove a general equivariant
bijection between rowmotion on $P\times [k]$, where $P$ is \textit{any} finite poset, which
is similar to the equivariant bijections Haddadan used in her work with triangular and
rectangular posets.

\subsection{Whirling functions between finite sets}\label{ss:whirl}

Let $\mathcal{F}\subseteq [k]^{[n]}$ be a family of functions $f:[n]
\to [k]$. For the rest of section~\ref{ss:whirl}, we use $\{1,\dots,k\}=[k]$ to represent the congruence
classes of $\mathbb{Z}/k\mathbb{Z}$, as opposed to the usual $\{0,1,\dots,k-1\}$. For fixed
values of $k$ and $n$, we represent such functions in \emph{one-line} notation, e.g.,  $f =
21344$ represents the function $f\in [4]^{[5]}$ with $f(1)=2$, $f(2)=1$, $f(3)=3$, $f(4)=4$,
and $f(5)=4$.  
\begin{definition}[\protect{\cite[Definition~2.3]{JPR18}}]
For $f\in \mathcal{F}$ we define the \emph{whirl} $w_i:\mathcal{F}\to\mathcal{F}$ at index $i$ as follows: repeatedly add 1 (modulo $k$) to the value of $f(i)$ until we get a function in $\mathcal{F}$. 
\end{definition}

\begin{ex}
Let $\mathcal{F} = \{f\in [4]^{[5]}:f(1)\not =f( 2)\}$. 
If we apply $w_2$ to $f = 21344$, adding 1 in the second position gives $22344$, but this is not in $\mathcal{F}$. Adding 1 again in this position gives the result: $w_2(f)=23344$. 
\end{ex}

\noindent We will now highlight some specific results from the paper where whirling was
first introduced. Let $\mathrm{Inj}_m(n,k)$ be the set of \emph{$m$-injective functions},
that is, functions $f:[n]\to[k]$ such that $\#f^{-1}(t)\leq m$ for all $t\in[k]$. Similarly,
let $\mathrm{Sur}_m(n,k)$ be the set of \emph{$m$-surjective functions}, that is,
$f:[n]\to[k]$ such that $\#f^{-1}(t)\geq m$ for all $t\in[k]$. Note that injective functions
are $1$-injections and surjective functions are $1$-surjections. We also define the
statistic $\eta_j(f)=\# f^{-1}(\{j\})$. 

\begin{thm}[\protect{\cite[Theorem~2.11]{JPR18}\label{thm:JPRthm}}]
Fix $\mathcal{F}$ to be either $\mathrm{Inj}_m(n,k)$ or $\mathrm{Sur}_1(n,k)$ for given
$n,k,m\in \mathbb{P}$. Then under the action of $\mathbf{w}=w_n\circ w_{n-1}\circ \cdots
\circ w_1$ on $\mathcal{F}$, $\eta_j$ is $\frac{n}{k}$-mesic for any $j\in[k]$. 
\end{thm}
\noindent This result is conjectured to hold for $\mathrm{Sur}_m(n,k)$, but is still open for $m>1$.
Proof details can be found in Sections 2.2--2.4 of \cite{JPR18}. 

\begin{ex}
Here is the orbit of $\mathbf{w}$ on $\mathrm{Inj}_1(3,6)$ containing $f=415$.

$$415\xrightarrow{\mathbf{w}} 621\xrightarrow{\mathbf{w}} 342\xrightarrow{\mathbf{w}} 563\xrightarrow{\mathbf{w}} 124\xrightarrow{\mathbf{w}} 356\xrightarrow{\mathbf{w}} 412\xrightarrow{\mathbf{w}} 534\xrightarrow{\mathbf{w}} 651\xrightarrow{\mathbf{w}} 263\return{\mathbf{w}}$$

Figure~\ref{fig:chunk} shows the corresponding \textit{orbit board} (a matrix whose rows are
the successive orbit elements) partitioned into ``chunks,'' each of which contains exactly
the set of numbers $\{1,2,\dotsc ,k \}$.  Notice that each value $1,2,\dots,6$
appear exactly 5 times in this orbit of size 10, in accordance with the $1/2$-mesy of
Theorem~\ref{thm:JPRthm}. See \cite[Section 2]{JPR18} for precise definitions.  
\end{ex}

\begin{figure}\centering
\begin{tikzpicture}[baseline={([yshift=-.8ex]current bounding box.center)}]
    \tikzColoredSquare{0}{0}{4}{white!45!yellow}
    \tikzColoredSquare{1}{0}{1}{white!45!blue}
    \tikzColoredSquare{2}{0}{5}{white!45!yellow}
    
    \tikzColoredSquare{0}{-1}{6}{white!45!yellow}
    \tikzColoredSquare{1}{-1}{2}{white!45!blue}
    \tikzColoredSquare{2}{-1}{1}{white!45!red}
    
    \tikzColoredSquare{0}{-2}{3}{white!45!blue}
    \tikzColoredSquare{1}{-2}{4}{white!45!blue}
    \tikzColoredSquare{2}{-2}{2}{white!45!red}
    
    \tikzColoredSquare{0}{-3}{5}{white!45!blue}
    \tikzColoredSquare{1}{-3}{6}{white!45!blue}
    \tikzColoredSquare{2}{-3}{3}{white!45!red}
    
    \tikzColoredSquare{0}{-4}{1}{white!45!green}
    \tikzColoredSquare{1}{-4}{2}{white!45!green}
    \tikzColoredSquare{2}{-4}{4}{white!45!red}
    
    \tikzColoredSquare{0}{-5}{3}{white!45!green}
    \tikzColoredSquare{1}{-5}{5}{white!45!red}
    \tikzColoredSquare{2}{-5}{6}{white!45!red}
    
    \tikzColoredSquare{0}{-6}{4}{white!45!green}
    \tikzColoredSquare{1}{-6}{1}{white!45!orange}
    \tikzColoredSquare{2}{-6}{2}{white!45!orange}
    
    \tikzColoredSquare{0}{-7}{5}{white!45!green}
    \tikzColoredSquare{1}{-7}{3}{white!45!orange}
    \tikzColoredSquare{2}{-7}{4}{white!45!orange}
    
    \tikzColoredSquare{0}{-8}{6}{white!45!green}
    \tikzColoredSquare{1}{-8}{5}{white!45!orange}
    \tikzColoredSquare{2}{-8}{1}{white!45!yellow}
    
    \tikzColoredSquare{0}{-9}{2}{white!45!yellow}
    \tikzColoredSquare{1}{-9}{6}{white!45!orange}
    \tikzColoredSquare{2}{-9}{3}{white!45!yellow}

\end{tikzpicture}
\vspace{-0.5em}
\caption{The orbit board of $\mathbf{w}$ on $\mathrm{Inj}_1(3,6)$ containing $f=415$,
partitioned into ``chunks''.}\label{fig:chunk}
\end{figure}

\subsection{Whirling $k$-bounded $P$-partitions and rowmotion}\label{ss:whirlP}
Now we extend the definition of whirling to $k$-bounded $P$-partitions, and prove our main
general result relating whirling and rowmotion.    
Throughout the rest of the paper, $P$ will denote a finite poset. Define $[0,k] := \{0,1,2,\dots,k\}$.
%with $n=\#P$, and let $L$ be a bijection from $P$ to $[n]$, called a \emph{labeling} of $P$.  
A \emph{$P$-partition} is a map $\sigma$ from $P$ to $\NN$ such that if $x<_P y$, then $\sigma(x)\geq\sigma(y)$~\cite[Ch. 3]{Sta11}.
\begin{definition}
A $k$-\emph{bounded} $P$-\emph{partition} is a function $f:P\to [0,k]$ such that if $x\leq_P y$, then  $f(x)\ge f(y)$. Let $\mathcal{F}_k(P)$ be the set of all such functions. 
\end{definition}

\noindent 
Throughout the rest of the paper we use $\{0,1,\dots,k\}$ to represent the congruence
classes of $\mathbb{Z}/(k+1)\mathbb{Z}$, as usual.
\begin{definition}
For $f\in \mathcal{F}_k(P)$ and $x\in P$, define $w_x:\mathcal{F}_k(P)\to\mathcal{F}_k(P)$, called the \emph{whirl at }$x$, as follows: repeatedly add 1 ($\textrm{mod } k+1$) to the value of $f(x)$ until we get a function in $\mathcal{F}_k(P)$.  This new function is $w_x(f)$.  \end{definition}
\noindent The case $k=1$ of the above definition recovers toggling of order ideals (Def.~\ref{def:toggles}). 
\begin{prop}\label{prop:inccommute}
If $x,y\in P$ are incomparable, then $w_xw_y(f) = w_yw_x(f)$.
\end{prop}

\begin{proof}
Since $x$ and $y$ are incomparable, there are no inequalities constraining the relationship
between $f(x)$ and $f(y)$. So $w_xw_y = w_yw_x$.
\end{proof} 

\begin{definition}\label{def:whirl}
Let $(x_{1},x_{2},\dots, x_{p})$ be a linear extension of $P$.  Define $w:
\mathcal{F}_k(P)\rightarrow \mathcal{F}_k(P)$ by $w:=w_{x_{1}}w_{x_{2}}\dots w_{x_{p}}$.
The above proposition shows that this is well-defined, since one can get from any linear
extension to any other by a sequence of interchanges of incomparable elements.  
\end{definition}

\begin{ex}
Let $P$ be the $\VV$ poset with labels \begin{tikzpicture}[scale = .65, baseline={([yshift=-.8ex]current bounding box.center)}]
    \node at (-1,1) (a) {$\ell$};
    \node at (0,0) (b) {$c$};
    \node at (1,1) (c) {$r$};
    \draw[-] (a) -- (b);
    \draw[-] (c) -- (b);
\end{tikzpicture}, $k=2$, and $w=w_c w_r w_\ell$.
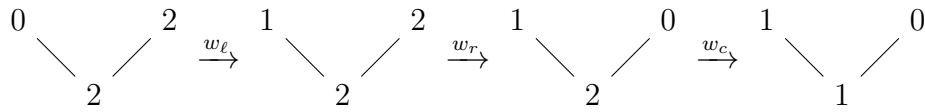
\[
\begin{tikzpicture}[baseline={([yshift=-.8ex]current bounding box.center)}]
    \node at (-1,1) (a) {$0$};
    \node at (0,0) (b) {$2$};
    \node at (1,1) (c) {$2$};
    \draw[-] (a) -- (b);
    \draw[-] (c) -- (b);
\end{tikzpicture}
\xrightarrow{w_\ell}
\begin{tikzpicture}[baseline={([yshift=-.8ex]current bounding box.center)}]
    \node at (-1,1) (a) {$1$};
    \node at (0,0) (b) {$2$};
    \node at (1,1) (c) {$2$};
    \draw[-] (a) -- (b);
    \draw[-] (c) -- (b);
\end{tikzpicture}
\xrightarrow{w_r}
\begin{tikzpicture}[baseline={([yshift=-.8ex]current bounding box.center)}]
    \node at (-1,1) (a) {$1$};
    \node at (0,0) (b) {$2$};
    \node at (1,1) (c) {$0$};
    \draw[-] (a) -- (b);
    \draw[-] (c) -- (b);
\end{tikzpicture}
\xrightarrow{w_c}
\begin{tikzpicture}[baseline={([yshift=-.8ex]current bounding box.center)}]
    \node at (-1,1) (a) {$1$};
    \node at (0,0) (b) {$1$};
    \node at (1,1) (c) {$0$};
    \draw[-] (a) -- (b);
    \draw[-] (c) -- (b);
\end{tikzpicture}
\]
\end{ex}

There is a natural bijection between order ideals of a poset $P$ and $1$-bounded $P$-partitions in $\mathcal{F}_1(P)$. Specifically, a $1$-bounded $P$-partition in $\mathcal{F}_1(P)$ is simply the indicator function of an order ideal $I\in J(P)$. We extend this to an equivariant bijection $\mathcal{F}_k(P)\to\mathcal{J}(\mathcal{P}\times [k])$ which sends $w$ to $\jrho$,  meaning the following diagram commutes.

\[
\begin{tikzpicture}[baseline={([yshift=-.8ex]current bounding box.center)}]
    \node at (0,2) (a) {$\mathcal{F}_k(P)$};
    \node at (0,0) (b) {$\mathcal{J}(P\times[k])$};
    \node at (4,2) (c) {$\mathcal{F}_k(P)$};
    \node at (4,0) (d) {$\mathcal{J}(P\times[k])$};
    \draw[<->] (a) -- (b);
    \draw[<->] (c) -- (d);
    \draw[->] (a) -- (c);
    \draw[->] (b) -- (d);
    
    \node at (2,0.25) (x) {$\jrho$};
    \node at (2,2.25) (x) {$w$};
\end{tikzpicture}
\]

We will call the chains $\{(x,1),(x,2),\dots,(x,k)\}\subseteq P\times [k]$, for $x\in P$, the \emph{fibers} of $P\times [k]$, and construct an equivariant bijection that first sends $w_x$ to order-ideal toggling down the fiber $\{(x,1),(x,2),\dots,(x,k)\}$.

\begin{ex}\label{eg:eqbij} Here is an example of the equivariant bijection between rowmotion
on order ideals of $[3]\times [4]$ and whirling on $\mathcal{F}_4([3])$.  

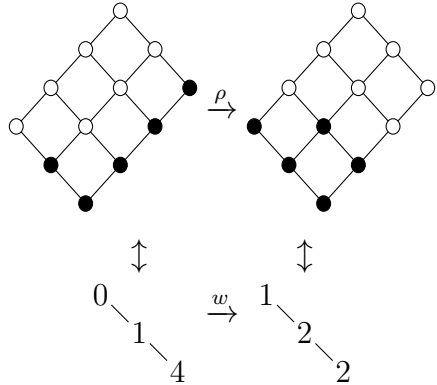
\[
\begin{tikzpicture}[scale = 0.51,xscale = 0.90, baseline={([yshift=-.8ex]current bounding box.center)}]
    \draw  (0,0) -- (-2,2);
    \draw  (1,1) -- (-1,3);
    \draw  (2,2) -- (0,4);
    \draw  (3,3) -- (1,5);
    
    \draw  (0,0) -- (3,3);
    \draw  (-1,1) -- (2,4);
    \draw  (-2,2) -- (1,5);
    
    \draw [fill] (0,0) circle [radius=.2];
    \draw [fill] (-1,1) circle [radius=.2];
    \draw [fill=white] (-2,2) circle [radius=.2];
    
    \draw [fill] (1,1) circle [radius=.2];
    \draw [fill=white] (0,2) circle [radius=.2];
    \draw [fill=white] (-1,3) circle [radius=.2];
    
    \draw [fill] (2,2) circle [radius=.2];
    \draw [fill=white] (1,3) circle [radius=.2];
    \draw [fill=white] (0,4) circle [radius=.2];
    
    \draw [fill] (3,3) circle [radius=.2];
    \draw [fill=white] (2,4) circle [radius=.2];
    \draw [fill=white] (1,5) circle [radius=.2];
\end{tikzpicture}
 \xrightarrow{\rho}
\begin{tikzpicture}[scale = 0.51,xscale = 0.90, baseline={([yshift=-.8ex]current bounding box.center)}]
    \draw  (0,0) -- (-2,2);
    \draw  (1,1) -- (-1,3);
    \draw  (2,2) -- (0,4);
    \draw  (3,3) -- (1,5);
    
    \draw  (0,0) -- (3,3);
    \draw  (-1,1) -- (2,4);
    \draw  (-2,2) -- (1,5);
    
    \draw [fill] (0,0) circle [radius=.2];
    \draw [fill] (-1,1) circle [radius=.2];
    \draw [fill] (-2,2) circle [radius=.2];
    
    \draw [fill] (1,1) circle [radius=.2];
    \draw [fill] (0,2) circle [radius=.2];
    \draw [fill=white] (-1,3) circle [radius=.2];
    
    \draw [fill=white] (2,2) circle [radius=.2];
    \draw [fill=white] (1,3) circle [radius=.2];
    \draw [fill=white] (0,4) circle [radius=.2];
    
    \draw [fill=white] (3,3) circle [radius=.2];
    \draw [fill=white] (2,4) circle [radius=.2];
    \draw [fill=white] (1,5) circle [radius=.2];
\end{tikzpicture}
\]

\[
\begin{tikzpicture}[scale = 0.51, baseline={([yshift=-.8ex]current bounding box.center)}]
    
    \node at (0,2) (w) {$\updownarrow$};
    \node at (1,-1) (b) {$4$};
    \node at (0,0) (c) {$1$};
    \node at (-1,1) (d) {$0$};
    \draw [shorten >= 2mm,shorten <= 2mm] (0,0) -- (-1,1);
    \draw [shorten >= 2mm,shorten <= 2mm] (0,0) -- (1,-1);
\end{tikzpicture}
 \xrightarrow{w}
\begin{tikzpicture}[scale = 0.51, baseline={([yshift=-.8ex]current bounding box.center)}]
    \node at (0,2) (w) {$\updownarrow$};
    \node at (1,-1) (b) {$2$};
    \node at (0,0) (c) {$2$};
    \node at (-1,1) (d) {$1$};
    \draw [shorten >= 2mm,shorten <= 2mm] (0,0) -- (-1,1);
    \draw [shorten >= 2mm,shorten <= 2mm] (0,0) -- (1,-1);
\end{tikzpicture}
\]
\end{ex}

\begin{lem}\label{lem:eqwhirl}
There is an equivariant bijection between $\mathcal{F}_k\left(P\right)$ and $\mathcal{J}(P\times [k])$ which sends $w_x$ to the toggle product $\tau_{(x,1)}\tau_{(x,2)}\dots\tau_{(x,k)}$.
\end{lem}
\begin{proof}
The bijection $\phi$ between $\mathcal{F}_k\left(P\right)$ and $\mathcal{J}(P\times [k])$
will be a natural one, where $\phi(f) = I$ implies $f(x)$ counts the number of elements in the intersection of $I$ and the fiber at $x$ in $P\times [k]$. That is, the intersection of the fiber at $x$ in $P\times [k]$, $\{(x,1),(x,2),\dots,(x,k)\}$, and an order ideal $I$ is $\{(x,1),(x,2),\dots,(x,d)\}$ for some $d\in[k]$, if and only if $f(x) = d$. Otherwise, the intersection is empty if and only if $f(x)=0$.

We will now complete the proof by fixing $f$ and $x$ and showing that 
\[
\phi(w_x(f)) = \tau_{(x,1)}\tau_{(x,2)}\dots\tau_{(x,k)}(\phi(f)).
\]

Let $m=\max(\{f(y):y\lessdot x,\text{ for } y\in P\}\cup\{0\})$ and $M=\min(\{f(y):x\lessdot y,\text{ for } y\in P\}\cup\{k\})$, so $m\leq f(x) \leq M$. 
There are two cases:

\begin{itemize}
    \item \textit{Case 1.} If $f(x) < M$, then $(w_x(f))(x) = f(x)+1$.
    \item \textit{Case 2.} If $f(x) = M$, then $f(x)+1$ is greater than $M$ so we increment through values modulo $k+1$ until we get to $(w_x(f))(x)=m$.
\end{itemize}

On the other hand, we consider the action of $\tau =
\tau_{(x,1)}\tau_{(x,2)}\dots\tau_{(x,k)}$ on $\phi(f)$. Recall that $\tau_{(x,i)}$ does not
change the order ideal whenever $(x,i)$ is not maximal in the order ideal or minimal in the
complement of the order ideal. Because of this, there are only two outcomes of $\tau(S(f))$,
either one element was added to order ideal, or elements were removed until removing anymore
would not result in an order ideal. This reflects the two cases above, which completes the
proof. 
\end{proof} 

The following general theorem is extremely useful for studying the dynamics of rowmotion for certain
kinds of posets.  Such an equivariant bijection was used by
Haddadan~\cite[Lemma~3.1]{shahrzad} in the specific situation she studied of rowmotion and promotion
on triangular and rectangular posets, but has not been stated in anything like this
generality in literature we have seen.

\begin{figure}
\begin{tikzpicture}[baseline={([yshift=-.8ex]current bounding box.center)}]
    \node at (0,-1) (x) {$\VV$};
    
    \node at (-1,1) (a) {$2$};
    \node at (0,0) (b) {$1$};
    \node at (1,1) (c) {$3$};
    \draw[-] (a) -- (b);
    \draw[-] (c) -- (b);
\end{tikzpicture}
\hspace{4cm}
\begin{tikzpicture}[baseline={([yshift=-.8ex]current bounding box.center)}]
    \node at (0,-1) (x) {$\VV\times[4]$};
    
    \node at (-1,1) (a1) {$5$};
    \node at (0,0) (b1) {$1$};
    \node at (1,1) (c1) {$9$};
    
    \draw[-] (a1) -- (b1);
    \draw[-] (c1) -- (b1);

    \node at (-1,2) (a2) {$6$};
    \node at (0,1) (b2) {$2$};
    \node at (1,2) (c2) {$10$};
    
    \draw[-] (a1) -- (a2);
    \draw[-] (b1) -- (b2);
    \draw[-] (c1) -- (c2);
    
    \draw[-] (a2) -- (b2);
    \draw[-] (c2) -- (b2);

    \node at (-1,3) (a3) {$7$};
    \node at (0,2) (b3) {$3$};
    \node at (1,3) (c3) {$11$};
    
    \draw[-] (a3) -- (a2);
    \draw[-] (b3) -- (b2);
    \draw[-] (c3) -- (c2);
    
    \draw[-] (a3) -- (b3);
    \draw[-] (c3) -- (b3);

    \node at (-1,4) (a4) {$8$};
    \node at (0,3) (b4) {$4$};
    \node at (1,4) (c4) {$12$};
    
    \draw[-] (a3) -- (a4);
    \draw[-] (b3) -- (b4);
    \draw[-] (c3) -- (c4);
    
    \draw[-] (a4) -- (b4);
    \draw[-] (c4) -- (b4);
\end{tikzpicture}
\caption{A linear extension of $\VV$ on the left and the linear extension for $\VV\times[4]$
constructed from it on the right, as in the proof of Theorem~\ref{thm:eqwhirl}.}
\label{fig:crossexample}
\end{figure}
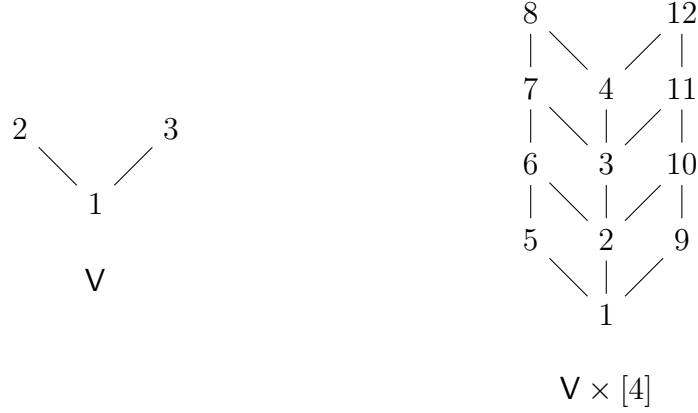

\begin{thm}\label{thm:eqwhirl}
Fix any linear extension $(x_{1},x_{2},\dots ,x_{p})\in\mathcal{L}(P)$. There is an
equivariant bijection between  $\mathcal{F}_k\left(P\right)$ and 
$\mathcal{J}(P\times [k])$
which sends whirling, $w=w_{x_{1}}w_{x_{2}}\cdots w_{x_{p}}$ on $\mathcal{F}_k\left(P\right)$ to rowmotion on
$\mathcal{J}(P\times[k])$.  
\end{thm}

\begin{proof}
The bijection is based on the one from Lemma~\ref{lem:eqwhirl}. We will construct the linear extension that traverses
$P\times[k]$ fiber by fiber, i.e., $(x_i,j) \to x_{k(i-1)+j}$. Then the composite map
$\tau_{x_1}\tau_{x_2}\cdots\tau_{x_{kp}}$ 
on $P\times[k]$ coincides with rowmotion by Proposition~\ref{prop:togall}. By
Lemma~\ref{lem:eqwhirl}, the toggles
$\tau_{x_{k(i-1)+1}}\tau_{x_{k(i-1)+2}}\cdots\tau_{x_{k(i-1)+k}}$ coincide with
$w_{x_i}$. Together these give us that $w=w_{x_{1}}w_{x_{2}}\cdots w_{x_{p}}$, is exactly 
rowmotion. 
\end{proof}
\begin{ex}
Figure~\ref{fig:crossexample} shows a linear extension of the poset, $\VV$ and the
corresponding linear extension for $\VV\times[4]$ on the right, as in the proof of Theorem~\ref{thm:eqwhirl}, of $\VV\times[4]$.
\end{ex}

The following definitions will allow us to partition orbit boards of whirling into subsets
called \emph{whorms}. The partitioning of orbit boards of combinatorial
actions into some special class of disjoint subsets is a useful tool, used by several
authors, including different kinds of ``snakes''~\cite{shahrzad,JR18}, and ``chunks''~\cite[]{JPR18}. 

\begin{definition}
For any $x\in P$ and $f\in \mathcal{F}_k(P)$, define $(x,f)$ to be a \emph{whirl element}. The whirl element $(y,g)$ is \emph{whirl successive} to $(x,f)$ if either:
\begin{enumerate}
    \item $y=x$ and $g(y)=w(f)(x)=f(x)+1$, or 
    \item $x$ covers $y$, $f=g$, and $f(x)=g(y)$.
\end{enumerate}
\end{definition}
\noindent 
We consider whirl-successive elements to be whirl elements which are one step away from
each other, either by moving one covering relation down the poset or by whirling the
function at the element, and ending one label greater.  
\emph{While we must consider the entire $P$-partitions $f$ and $g$ to check whether two whirl
elements are whorm-connected, we think of whirl elements as being simply $(x,f(x))$, the
location and its label, and indicate them in this way in the examples that follow.}  
 
\begin{definition}\label{def:whorms} Two whirl elements $(x,f)$ and $(y,g)$ are
\emph{whorm-connected} if there exists a sequence of whirl-successive elements $\{(x,f) =
(x_0,f_0), (x_1,f_1),\dots , (x_p,f_p) = (y,g)\}$. A \emph{whorm} is a maximal set of
whorm-connected whirl elements, that is, if $(x,f)$ is in a whorm and $(x,f)$ is
whorm-connected to $(y,g)$, then $(y,g)$ is in the whorm.  
\end{definition}

\begin{ex} An orbit of whirling $P$-partitions (for $P=[2]\times [2]$) with its four whorms
indicated by the same color and (redundantly) node-shape.  
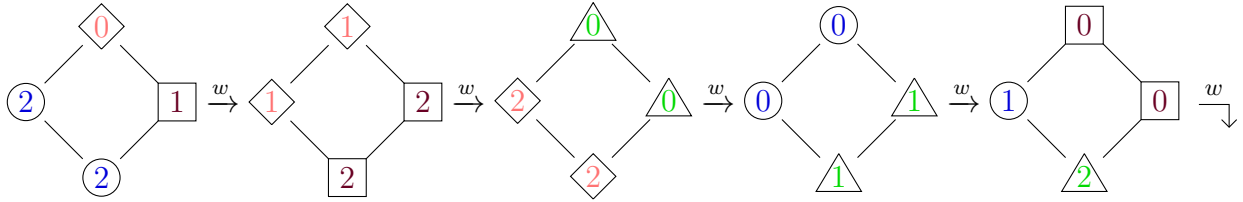
\[
\begin{tikzpicture}[baseline={([yshift=-.8ex]current bounding box.center)}]
    
    \draw (0,0) circle (7pt);
    \node at (0,0) (b) {$\color{black!15!blue}2$};
	\draw (-1,1) circle (7pt);
    \node at (-1,1) (a) {$\color{black!15!blue}2$};
	\draw[-] (1-0.25,1-0.25) -- (1+0.25,1-0.25) -- (1+0.25,1+0.25) -- (1-0.25,1+0.25) -- (1-0.25,1-0.25) ; %square
	\node at (1,1) (c) {$\color{black!45!purple}1$};
    \draw[-] (0,2-0.3) -- (0-0.3,2) -- (0,2+0.3) -- (0+0.3,2+0) -- (0,2-0.3) ; %diamond
    \node at (0,2) (d) {$\color{white!45!red} 0$};
    \draw[-] (a) -- (b);
    \draw[-] (c) -- (b);
    \draw[-] (d) -- (a);
    \draw[-] (d) -- (c);
\end{tikzpicture}
\xrightarrow[]{w}
\begin{tikzpicture}[baseline={([yshift=-.8ex]current bounding box.center)}]
	\draw[-] (0-0.25,0-0.25) -- (0+0.25,0-0.25) -- (0+0.25,0+0.25) -- (0-0.25,0+0.25) -- (0-0.25,0-0.25) ; %square
    \node at (0,0) (b) {$\color{black!45!purple}2$};
    \draw[-] (-1,1-0.3) -- (-1-0.3,1) -- (-1,1+0.3) -- (-1+0.3,1+0) -- (-1,1-0.3) ; %diamond
    \node at (-1,1) (a) {$\color{white!45!red} 1$};
	\draw[-] (1-0.25,1-0.25) -- (1+0.25,1-0.25) -- (1+0.25,1+0.25) -- (1-0.25,1+0.25) -- (1-0.25,1-0.25) ; %square    
	\node at (1,1) (c) {$\color{black!45!purple}2$};
    \draw[-] (0,2-0.3) -- (0-0.3,2) -- (0,2+0.3) -- (0+0.3,2+0) -- (0,2-0.3) ; %diamond
    \node at (0,2) (d) {$\color{white!45!red} 1$};
    \draw[-] (a) -- (b);
    \draw[-] (c) -- (b);
    \draw[-] (d) -- (a);
    \draw[-] (d) -- (c);
\end{tikzpicture}
\xrightarrow[]{w}
\begin{tikzpicture}[baseline={([yshift=-.8ex]current bounding box.center)}]
	
    \draw[-] (0,0-0.3) -- (0-0.3,0) -- (0,0+0.3) -- (0+0.3,0+0) -- (0,0-0.3) ; %diamond
    \node at (0,0) (b) {$\color{white!45!red}2$};
    \draw[-] (-1,1-0.3) -- (-1-0.3,1) -- (-1,1+0.3) -- (-1+0.3,1+0) -- (-1,1-0.3) ; %diamond
    \node at (-1,1) (a) {$\color{white!45!red}2$};
    \draw[-] (1,1+0.3) -- (1-0.30,1-0.20) -- (1+0.30,1-0.20) -- (1,1+0.3) ; %triangle
    \node at (1,1) (c) {$\color{black!15!green}0$};
    \draw[-] (0,2+0.3) -- (0-0.30,2-0.20) -- (0+0.30,2-0.20) -- (0,2+0.3) ; %triangle
    \node at (0,2) (d) {$\color{black!15!green}0$};
    \draw[-] (a) -- (b);
    \draw[-] (c) -- (b);
    \draw[-] (d) -- (a);
    \draw[-] (d) -- (c);
\end{tikzpicture}
\xrightarrow[]{w}
\begin{tikzpicture}[baseline={([yshift=-.8ex]current bounding box.center)}]
	\draw[-] (0,0+0.3) -- (0-0.30,0-0.20) -- (0+0.30,0-0.20) -- (0,0+0.3) ; %triangle
    \node at (0,0) (b) {$\color{black!15!green}1$};
    \draw (-1,1) circle (7pt);
    \node at (-1,1) (a) {$\color{black!15!blue}0$};
	\draw[-] (1,1+0.3) -- (1-0.30,1-0.20) -- (1+0.30,1-0.20) -- (1,1+0.3) ; %triangle
    \node at (1,1) (c) {$\color{black!15!green}1$};
    \draw (0,2) circle (7pt);
    \node at (0,2) (d) {$\color{black!15!blue}0$};
    \draw[-] (a) -- (b);
    \draw[-] (c) -- (b);
    \draw[-] (d) -- (a);
    \draw[-] (d) -- (c);
\end{tikzpicture}
\xrightarrow[]{w}
\begin{tikzpicture}[baseline={([yshift=-.8ex]current bounding box.center)}]
    \draw[-] (0,0+0.3) -- (0-0.30,0-0.20) -- (0+0.30,0-0.20) -- (0,0+0.3) ; %triangle
    \node at (0,0) (b) {$\color{black!15!green}2$};
    \draw (-1,1) circle (7pt);
    \node at (-1,1) (a) {$\color{black!15!blue}1$};
	\draw[-] (1-0.25,1-0.25) -- (1+0.25,1-0.25) -- (1+0.25,1+0.25) -- (1-0.25,1+0.25) -- (1-0.25,1-0.25) ; %square
    \node at (1,1) (c) {$\color{black!45!purple}0$};
    \draw[-] (0-0.25,2-0.25) -- (0+0.25,2-0.25) -- (0+0.25,2+0.25) -- (0-0.25,2+0.25) -- (0-0.25,2-0.25) ; %square
    \node at (0,2) (d) {$\color{black!45!purple}0$};
    \draw[-] (a) -- (b);
    \draw[-] (c) -- (b);
    \draw[-] (d) -- (a);
    \draw[-] (d) -- (c);
\end{tikzpicture}\return{w}
\]

\end{ex}

\section{Dynamics for rowmotion on $\VV\times[k]$}\label{sec:vee}
In this section we consider the dynamics of rowmotion acting on the order ideals of the chain of V's poset
$\VV_{k}$,  establishing its periodicity and finding interesting examples of homomesy. 

\subsection{The chain of V's poset}

While most posets currently known to enjoy nice dynamical properties under rowmotion seem to
be the root and minuscule posets associated with various types of Lie algebras, the current
poset seems to be \textit{sui generis}.  

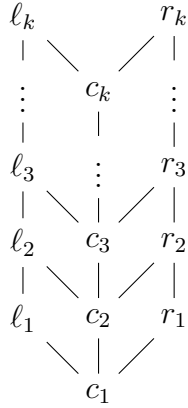
\begin{figure}
    \begin{tikzpicture}[baseline={([yshift=-.8ex]current bounding box.center)}]
    \node at (-1,1) (l2) {$\ell_1$};
    \node at (0,0) (c3) {$c_1$};
    \node at (1,1) (r2) {$r_1$};
    \draw[-] (c3) -- (l2);
    \draw[-] (c3) -- (r2);
    \node at (-1,2) (l3) {$\ell_2$};
    \node at (0,1) (c4) {$c_2$};
    \node at (1,2) (r3) {$r_2$};
    \draw[-] (c4) -- (l3);
    \draw[-] (c4) -- (r3);
    \node at (-1,3) (l4) {$\ell_3$};
    \node at (0,2) (c5) {$c_3$};
    \node at (1,3) (r4) {$r_3$};
    \draw[-] (c5) -- (l4);
    \draw[-] (c5) -- (r4);
    \node at (-1,5) (l5) {$\ell_k$};
    \node at (0,4) (c6) {$c_k$};
    \node at (1,5) (r5) {$r_k$};
    \draw[-] (c6) -- (l5);
    \draw[-] (c6) -- (r5);
    \node at (-1,4) (l) {$\vdots$};
    \node at (0,3) (c) {$\vdots$};
    \node at (1,4) (r) {$\vdots$};
    
    \draw[-] (c3) -- (c4);
    \draw[-] (c4) -- (c5);
    \draw[-] (c5) -- (c);
    \draw[-] (c) -- (c6);
    
    \draw[-] (l2) -- (l3);
    \draw[-] (l3) -- (l4);
    \draw[-] (l4) -- (l);
    \draw[-] (l) -- (l5);
    
    \draw[-] (r2) -- (r3);
    \draw[-] (r3) -- (r4);
    \draw[-] (r4) -- (r);
    \draw[-] (r) -- (r5);
\end{tikzpicture}
\caption{Hasse diagram of $\VV_k$ with our vertex-labeling convention.}\label{fig:vee}
\end{figure}

\begin{definition}
Let $\VV $ be the 3-element poset with Hasse diagram
$
\begin{tikzpicture}[scale = 0.5, baseline={([yshift=-.8ex]current bounding box.center)}]
    \draw[-] (-1,1) -- (0,0);
    \draw[-] (1,1) -- (0,0);
    
    \draw [fill] (0,0) circle [radius =.1];
    \draw [fill] (-1,1) circle [radius =.1];
    \draw [fill] (1,1) circle [radius =.1];
\end{tikzpicture}
$,
and define $\VV_k=\VV\times[k]$, where $[k]$ is the chain poset.  We call $\VV_{k}$ the
\emph{chain of V's poset}.  We will use the labeling convention of
Figure~\ref{fig:vee} throughout this section.  
\end{definition}

Our main goals for this section are the following theorems.  We will leverage the
equivariant bijection (Theorem~\ref{thm:eqwhirl}) and the notion of whorms from Section~\ref{ss:whirlP}.  

\begin{thm}\label{thm:Vorder}
The order of rowmotion on $\mathcal{J}(\VV_k)$ is $2(k+2)$.
\end{thm}

\begin{thm}\label{thm:Vhom}
Let $\chi_s$ be the indicator function for $s\in \VV_k$. We have the following homomesies for the action of $\jrho$ on $\mathcal{J}(\VV_k)$
\begin{enumerate}
    \item The statistic $\chi_{\ell_i}-\chi_{r_i}$ is $0$-mesic for all $i\in [k]$.
    \item The statistic $\chi_{\ell_1} + \chi_{r_1} - \chi_{c_k}$ is $\frac{2(k-1)}{k+2}$-mesic.
\end{enumerate}

\end{thm}

\begin{figure}
\noindent\fbox{\parbox{\textwidth}{
\[
\begin{tikzpicture}[scale = 0.9024, baseline={([yshift=-.8ex]current bounding box.center)}]
    \draw[-] (-1,1) -- (-1,3);
    \draw[-] (1,1) -- (1,3);
    \draw[-] (0,0) -- (0,2);

    \draw[-] (-1,1) -- (0,0);
    \draw[-] (1,1) -- (0,0);
    
    \draw [fill=white] (-1,1) circle [radius =.1];
    \draw [fill=white] (0,0) circle [radius =.1];
    \draw [fill=white] (1,1) circle [radius =.1];
    
    \draw[-] (-1,2) -- (0,1);
    \draw[-] (1,2) -- (0,1);
    
    \draw [fill=white] (-1,2) circle [radius =.1];
    \draw [fill=white] (0,1) circle [radius =.1];
    \draw [fill=white] (1,2) circle [radius =.1];
    
    \draw[-] (-1,3) -- (0,2);
    \draw[-] (1,3) -- (0,2);
    
    \draw [fill=white] (-1,3) circle [radius =.1];
    \draw [fill=white] (0,2) circle [radius =.1];
    \draw [fill=white] (1,3) circle [radius =.1];
\end{tikzpicture}
\xrightarrow{\jrho}
\begin{tikzpicture}[scale = 0.9024, baseline={([yshift=-.8ex]current bounding box.center)}]
    \draw[-] (-1,1) -- (-1,3);
    \draw[-] (1,1) -- (1,3);
    \draw[-] (0,0) -- (0,2);

    \draw[-] (-1,1) -- (0,0);
    \draw[-] (1,1) -- (0,0);
    
    \draw [fill=white] (-1,1) circle [radius =.1];
    \draw [fill] (0,0) circle [radius =.1];
    \draw [fill=white] (1,1) circle [radius =.1];
    
    \draw[-] (-1,2) -- (0,1);
    \draw[-] (1,2) -- (0,1);
    
    \draw [fill=white] (-1,2) circle [radius =.1];
    \draw [fill=white] (0,1) circle [radius =.1];
    \draw [fill=white] (1,2) circle [radius =.1];
    
    \draw[-] (-1,3) -- (0,2);
    \draw[-] (1,3) -- (0,2);
    
    \draw [fill=white] (-1,3) circle [radius =.1];
    \draw [fill=white] (0,2) circle [radius =.1];
    \draw [fill=white] (1,3) circle [radius =.1];
\end{tikzpicture}
\xrightarrow{\jrho}
\begin{tikzpicture}[scale = 0.9024, baseline={([yshift=-.8ex]current bounding box.center)}]
    \draw[-] (-1,1) -- (-1,3);
    \draw[-] (1,1) -- (1,3);
    \draw[-] (0,0) -- (0,2);

    \draw[-] (-1,1) -- (0,0);
    \draw[-] (1,1) -- (0,0);
    
    \draw [fill] (-1,1) circle [radius =.1];
    \draw [fill] (0,0) circle [radius =.1];
    \draw [fill] (1,1) circle [radius =.1];
    
    \draw[-] (-1,2) -- (0,1);
    \draw[-] (1,2) -- (0,1);
    
    \draw [fill=white] (-1,2) circle [radius =.1];
    \draw [fill] (0,1) circle [radius =.1];
    \draw [fill=white] (1,2) circle [radius =.1];
    
    \draw[-] (-1,3) -- (0,2);
    \draw[-] (1,3) -- (0,2);
    
    \draw [fill=white] (-1,3) circle [radius =.1];
    \draw [fill=white] (0,2) circle [radius =.1];
    \draw [fill=white] (1,3) circle [radius =.1];
\end{tikzpicture}
\xrightarrow{\jrho}
\begin{tikzpicture}[scale = 0.9024, baseline={([yshift=-.8ex]current bounding box.center)}]
    \draw[-] (-1,1) -- (-1,3);
    \draw[-] (1,1) -- (1,3);
    \draw[-] (0,0) -- (0,2);

    \draw[-] (-1,1) -- (0,0);
    \draw[-] (1,1) -- (0,0);
    
    \draw [fill] (-1,1) circle [radius =.1];
    \draw [fill] (0,0) circle [radius =.1];
    \draw [fill] (1,1) circle [radius =.1];
    
    \draw[-] (-1,2) -- (0,1);
    \draw[-] (1,2) -- (0,1);
    
    \draw [fill] (-1,2) circle [radius =.1];
    \draw [fill] (0,1) circle [radius =.1];
    \draw [fill] (1,2) circle [radius =.1];
    
    \draw[-] (-1,3) -- (0,2);
    \draw[-] (1,3) -- (0,2);
    
    \draw [fill=white] (-1,3) circle [radius =.1];
    \draw [fill] (0,2) circle [radius =.1];
    \draw [fill=white] (1,3) circle [radius =.1];
\end{tikzpicture}\xrightarrow{\jrho}
\begin{tikzpicture}[scale = 0.9024, baseline={([yshift=-.8ex]current bounding box.center)}]
    \draw[-] (-1,1) -- (-1,3);
    \draw[-] (1,1) -- (1,3);
    \draw[-] (0,0) -- (0,2);

    \draw[-] (-1,1) -- (0,0);
    \draw[-] (1,1) -- (0,0);
    
    \draw [fill] (-1,1) circle [radius =.1];
    \draw [fill] (0,0) circle [radius =.1];
    \draw [fill] (1,1) circle [radius =.1];
    
    \draw[-] (-1,2) -- (0,1);
    \draw[-] (1,2) -- (0,1);
    
    \draw [fill] (-1,2) circle [radius =.1];
    \draw [fill] (0,1) circle [radius =.1];
    \draw [fill] (1,2) circle [radius =.1];
    
    \draw[-] (-1,3) -- (0,2);
    \draw[-] (1,3) -- (0,2);
    
    \draw [fill] (-1,3) circle [radius =.1];
    \draw [fill] (0,2) circle [radius =.1];
    \draw [fill] (1,3) circle [radius =.1];
\end{tikzpicture}\return{\jrho}
\]
}}
\noindent\fbox{\parbox{\textwidth}{
\[
\begin{tikzpicture}[scale = 0.9024, baseline={([yshift=-.8ex]current bounding box.center)}]
    \draw[-] (-1,1) -- (-1,3);
    \draw[-] (1,1) -- (1,3);
    \draw[-] (0,0) -- (0,2);

    \draw[-] (-1,1) -- (0,0);
    \draw[-] (1,1) -- (0,0);
    
    \draw [fill] (-1,1) circle [radius =.1];
    \draw [fill] (0,0) circle [radius =.1];
    \draw [fill=white] (1,1) circle [radius =.1];
    
    \draw[-] (-1,2) -- (0,1);
    \draw[-] (1,2) -- (0,1);
    
    \draw [fill=white] (-1,2) circle [radius =.1];
    \draw [fill=white] (0,1) circle [radius =.1];
    \draw [fill=white] (1,2) circle [radius =.1];
    
    \draw[-] (-1,3) -- (0,2);
    \draw[-] (1,3) -- (0,2);
    
    \draw [fill=white] (-1,3) circle [radius =.1];
    \draw [fill=white] (0,2) circle [radius =.1];
    \draw [fill=white] (1,3) circle [radius =.1];
\end{tikzpicture}
\xrightarrow{\jrho}
\begin{tikzpicture}[scale = 0.9024, baseline={([yshift=-.8ex]current bounding box.center)}]
    \draw[-] (-1,1) -- (-1,3);
    \draw[-] (1,1) -- (1,3);
    \draw[-] (0,0) -- (0,2);

    \draw[-] (-1,1) -- (0,0);
    \draw[-] (1,1) -- (0,0);
    
    \draw [fill=white] (-1,1) circle [radius =.1];
    \draw [fill] (0,0) circle [radius =.1];
    \draw [fill] (1,1) circle [radius =.1];
    
    \draw[-] (-1,2) -- (0,1);
    \draw[-] (1,2) -- (0,1);
    
    \draw [fill=white] (-1,2) circle [radius =.1];
    \draw [fill] (0,1) circle [radius =.1];
    \draw [fill=white] (1,2) circle [radius =.1];
    
    \draw[-] (-1,3) -- (0,2);
    \draw[-] (1,3) -- (0,2);
    
    \draw [fill=white] (-1,3) circle [radius =.1];
    \draw [fill=white] (0,2) circle [radius =.1];
    \draw [fill=white] (1,3) circle [radius =.1];
\end{tikzpicture}
\xrightarrow{\jrho}
\begin{tikzpicture}[scale = 0.9024, baseline={([yshift=-.8ex]current bounding box.center)}]
    \draw[-] (-1,1) -- (-1,3);
    \draw[-] (1,1) -- (1,3);
    \draw[-] (0,0) -- (0,2);

    \draw[-] (-1,1) -- (0,0);
    \draw[-] (1,1) -- (0,0);
    
    \draw [fill] (-1,1) circle [radius =.1];
    \draw [fill] (0,0) circle [radius =.1];
    \draw [fill] (1,1) circle [radius =.1];
    
    \draw[-] (-1,2) -- (0,1);
    \draw[-] (1,2) -- (0,1);
    
    \draw [fill=white] (-1,2) circle [radius =.1];
    \draw [fill] (0,1) circle [radius =.1];
    \draw [fill] (1,2) circle [radius =.1];
    
    \draw[-] (-1,3) -- (0,2);
    \draw[-] (1,3) -- (0,2);
    
    \draw [fill=white] (-1,3) circle [radius =.1];
    \draw [fill] (0,2) circle [radius =.1];
    \draw [fill=white] (1,3) circle [radius =.1];
\end{tikzpicture}
\xrightarrow{\jrho}
\begin{tikzpicture}[scale = 0.9024, baseline={([yshift=-.8ex]current bounding box.center)}]
    \draw[-] (-1,1) -- (-1,3);
    \draw[-] (1,1) -- (1,3);
    \draw[-] (0,0) -- (0,2);

    \draw[-] (-1,1) -- (0,0);
    \draw[-] (1,1) -- (0,0);
    
    \draw [fill] (-1,1) circle [radius =.1];
    \draw [fill] (0,0) circle [radius =.1];
    \draw [fill] (1,1) circle [radius =.1];
    
    \draw[-] (-1,2) -- (0,1);
    \draw[-] (1,2) -- (0,1);
    
    \draw [fill] (-1,2) circle [radius =.1];
    \draw [fill] (0,1) circle [radius =.1];
    \draw [fill] (1,2) circle [radius =.1];
    
    \draw[-] (-1,3) -- (0,2);
    \draw[-] (1,3) -- (0,2);
    
    \draw [fill=white] (-1,3) circle [radius =.1];
    \draw [fill] (0,2) circle [radius =.1];
    \draw [fill] (1,3) circle [radius =.1];
\end{tikzpicture}
\xrightarrow{\jrho}
\begin{tikzpicture}[scale = 0.9024, baseline={([yshift=-.8ex]current bounding box.center)}]
    \draw[-] (-1,1) -- (-1,3);
    \draw[-] (1,1) -- (1,3);
    \draw[-] (0,0) -- (0,2);

    \draw[-] (-1,1) -- (0,0);
    \draw[-] (1,1) -- (0,0);
    
    \draw [fill] (-1,1) circle [radius =.1];
    \draw [fill] (0,0) circle [radius =.1];
    \draw [fill=white] (1,1) circle [radius =.1];
    
    \draw[-] (-1,2) -- (0,1);
    \draw[-] (1,2) -- (0,1);
    
    \draw [fill] (-1,2) circle [radius =.1];
    \draw [fill] (0,1) circle [radius =.1];
    \draw [fill=white] (1,2) circle [radius =.1];
    
    \draw[-] (-1,3) -- (0,2);
    \draw[-] (1,3) -- (0,2);
    
    \draw [fill] (-1,3) circle [radius =.1];
    \draw [fill] (0,2) circle [radius =.1];
    \draw [fill=white] (1,3) circle [radius =.1];
\end{tikzpicture}\xrightarrow{\jrho}
\]
\[
\begin{tikzpicture}[scale = 0.9024, baseline={([yshift=-.8ex]current bounding box.center)}]
    \draw[-] (-1,1) -- (-1,3);
    \draw[-] (1,1) -- (1,3);
    \draw[-] (0,0) -- (0,2);

    \draw[-] (-1,1) -- (0,0);
    \draw[-] (1,1) -- (0,0);
    
    \draw [fill=white] (-1,1) circle [radius =.1];
    \draw [fill] (0,0) circle [radius =.1];
    \draw [fill] (1,1) circle [radius =.1];
    
    \draw[-] (-1,2) -- (0,1);
    \draw[-] (1,2) -- (0,1);
    
    \draw [fill=white] (-1,2) circle [radius =.1];
    \draw [fill=white] (0,1) circle [radius =.1];
    \draw [fill=white] (1,2) circle [radius =.1];
    
    \draw[-] (-1,3) -- (0,2);
    \draw[-] (1,3) -- (0,2);
    
    \draw [fill=white] (-1,3) circle [radius =.1];
    \draw [fill=white] (0,2) circle [radius =.1];
    \draw [fill=white] (1,3) circle [radius =.1];
\end{tikzpicture}
\xrightarrow{\jrho}
\begin{tikzpicture}[scale = 0.9024, baseline={([yshift=-.8ex]current bounding box.center)}]
    \draw[-] (-1,1) -- (-1,3);
    \draw[-] (1,1) -- (1,3);
    \draw[-] (0,0) -- (0,2);

    \draw[-] (-1,1) -- (0,0);
    \draw[-] (1,1) -- (0,0);
    
    \draw [fill] (-1,1) circle [radius =.1];
    \draw [fill] (0,0) circle [radius =.1];
    \draw [fill=white] (1,1) circle [radius =.1];
    
    \draw[-] (-1,2) -- (0,1);
    \draw[-] (1,2) -- (0,1);
    
    \draw [fill=white] (-1,2) circle [radius =.1];
    \draw [fill] (0,1) circle [radius =.1];
    \draw [fill=white] (1,2) circle [radius =.1];
    
    \draw[-] (-1,3) -- (0,2);
    \draw[-] (1,3) -- (0,2);
    
    \draw [fill=white] (-1,3) circle [radius =.1];
    \draw [fill=white] (0,2) circle [radius =.1];
    \draw [fill=white] (1,3) circle [radius =.1];
\end{tikzpicture}
\xrightarrow{\jrho}
\begin{tikzpicture}[scale = 0.9024, baseline={([yshift=-.8ex]current bounding box.center)}]
    \draw[-] (-1,1) -- (-1,3);
    \draw[-] (1,1) -- (1,3);
    \draw[-] (0,0) -- (0,2);

    \draw[-] (-1,1) -- (0,0);
    \draw[-] (1,1) -- (0,0);
    
    \draw [fill] (-1,1) circle [radius =.1];
    \draw [fill] (0,0) circle [radius =.1];
    \draw [fill] (1,1) circle [radius =.1];
    
    \draw[-] (-1,2) -- (0,1);
    \draw[-] (1,2) -- (0,1);
    
    \draw [fill] (-1,2) circle [radius =.1];
    \draw [fill] (0,1) circle [radius =.1];
    \draw [fill=white] (1,2) circle [radius =.1];
    
    \draw[-] (-1,3) -- (0,2);
    \draw[-] (1,3) -- (0,2);
    
    \draw [fill=white] (-1,3) circle [radius =.1];
    \draw [fill] (0,2) circle [radius =.1];
    \draw [fill=white] (1,3) circle [radius =.1];
\end{tikzpicture}
\xrightarrow{\jrho}
\begin{tikzpicture}[scale = 0.9024, baseline={([yshift=-.8ex]current bounding box.center)}]
    \draw[-] (-1,1) -- (-1,3);
    \draw[-] (1,1) -- (1,3);
    \draw[-] (0,0) -- (0,2);

    \draw[-] (-1,1) -- (0,0);
    \draw[-] (1,1) -- (0,0);
    
    \draw [fill] (-1,1) circle [radius =.1];
    \draw [fill] (0,0) circle [radius =.1];
    \draw [fill] (1,1) circle [radius =.1];
    
    \draw[-] (-1,2) -- (0,1);
    \draw[-] (1,2) -- (0,1);
    
    \draw [fill] (-1,2) circle [radius =.1];
    \draw [fill] (0,1) circle [radius =.1];
    \draw [fill] (1,2) circle [radius =.1];
    
    \draw[-] (-1,3) -- (0,2);
    \draw[-] (1,3) -- (0,2);
    
    \draw [fill] (-1,3) circle [radius =.1];
    \draw [fill] (0,2) circle [radius =.1];
    \draw [fill=white] (1,3) circle [radius =.1];
\end{tikzpicture}
\xrightarrow{\jrho}
\begin{tikzpicture}[scale = 0.9024, baseline={([yshift=-.8ex]current bounding box.center)}]
    \draw[-] (-1,1) -- (-1,3);
    \draw[-] (1,1) -- (1,3);
    \draw[-] (0,0) -- (0,2);

    \draw[-] (-1,1) -- (0,0);
    \draw[-] (1,1) -- (0,0);
    
    \draw [fill=white] (-1,1) circle [radius =.1];
    \draw [fill] (0,0) circle [radius =.1];
    \draw [fill] (1,1) circle [radius =.1];
    
    \draw[-] (-1,2) -- (0,1);
    \draw[-] (1,2) -- (0,1);
    
    \draw [fill=white] (-1,2) circle [radius =.1];
    \draw [fill] (0,1) circle [radius =.1];
    \draw [fill] (1,2) circle [radius =.1];
    
    \draw[-] (-1,3) -- (0,2);
    \draw[-] (1,3) -- (0,2);
    
    \draw [fill=white] (-1,3) circle [radius =.1];
    \draw [fill] (0,2) circle [radius =.1];
    \draw [fill] (1,3) circle [radius =.1];
\end{tikzpicture}\return{\jrho}
\]
}}
\caption{The two orbits of rowmotion on order ideals of $\VV \times [3]$}
\label{fig:V3orbits}
\end{figure}
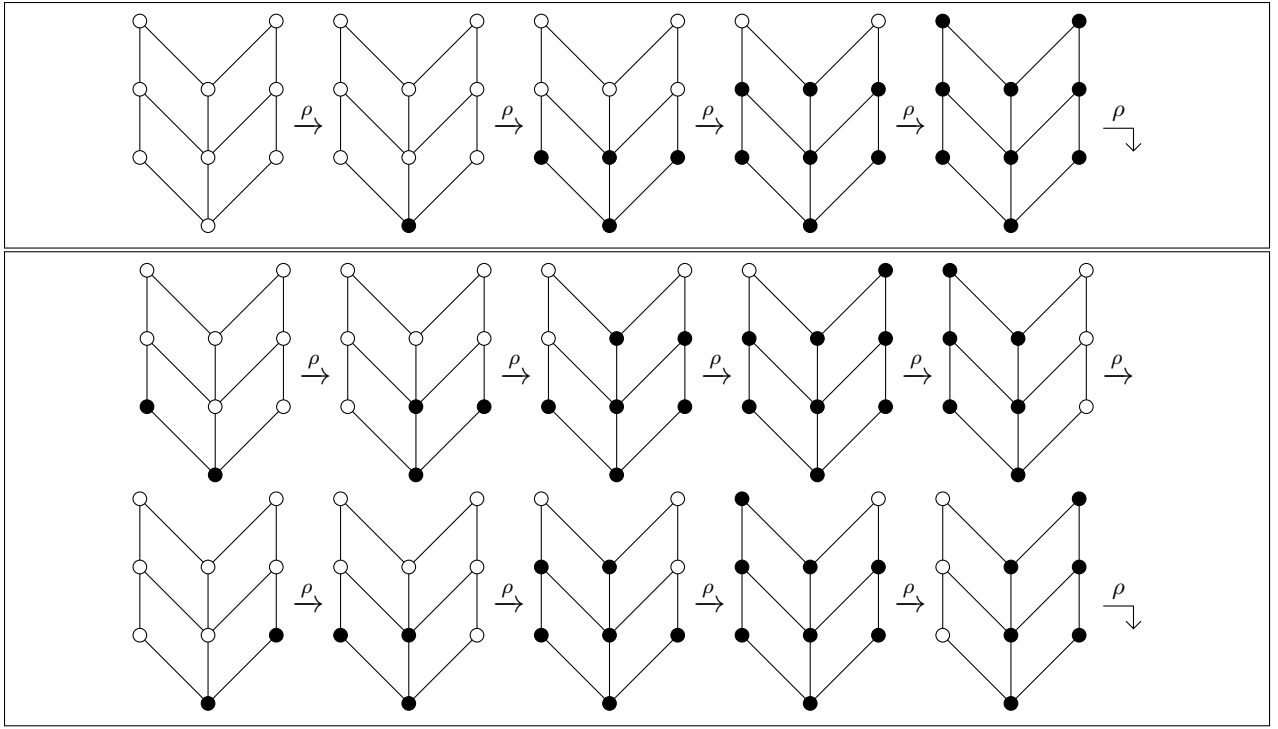

\begin{ex}\label{ex:V3orbits}
Figure~\ref{fig:V3orbits} shows the two $\jrho$-orbits, of sizes 5 and 10, on
$\mathcal{J}(\VV_4)$, confirming Theorem~\ref{thm:Vorder}.  It is also easy to check the
homomesies of Theorem~\ref{thm:Vhom}.    
\end{ex}

\begin{ex} \label{exorbit} Here is $\jrho$-orbit on $\mathcal{J}(\VV_4)$ of length 4, which
divides the order $2(4+2)=12$. Across the orbit the total number of elements at rank 1 in the side
fibers is 6, minus the two at the top of the center fiber, for an average of
$\frac{6-2}{4}=1 = \frac{2(4-1)}{4+2}$, agreeing with Theorem~\ref{thm:Vhom}(2).
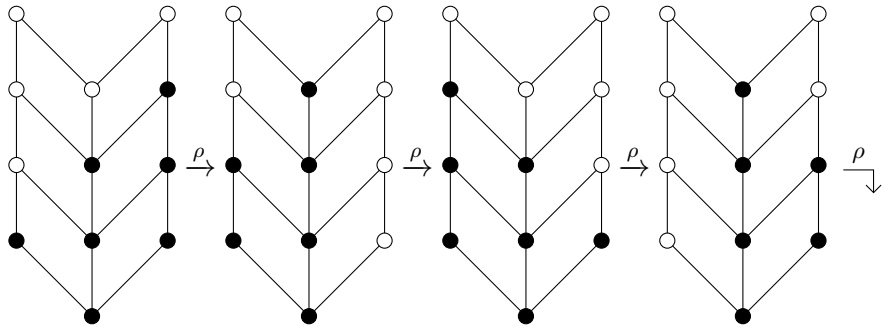
\[
\begin{tikzpicture}[baseline={([yshift=-.8ex]current bounding box.center)}]
    \draw[-] (-1,1) -- (-1,4);
    \draw[-] (1,1) -- (1,4);
    \draw[-] (0,0) -- (0,3);

    \draw[-] (-1,1) -- (0,0);
    \draw[-] (1,1) -- (0,0);
    
    \draw [fill] (-1,1) circle [radius =.1];
    \draw [fill] (0,0) circle [radius =.1];
    \draw [fill] (1,1) circle [radius =.1];
    
    \draw[-] (-1,2) -- (0,1);
    \draw[-] (1,2) -- (0,1);
    
    \draw [fill=white] (-1,2) circle [radius =.1];
    \draw [fill] (0,1) circle [radius =.1];
    \draw [fill] (1,2) circle [radius =.1];
    
    \draw[-] (-1,3) -- (0,2);
    \draw[-] (1,3) -- (0,2);
    
    \draw [fill=white] (-1,3) circle [radius =.1];
    \draw [fill] (0,2) circle [radius =.1];
    \draw [fill] (1,3) circle [radius =.1];
    
    \draw[-] (-1,4) -- (0,3);
    \draw[-] (1,4) -- (0,3);
    
    \draw [fill=white] (-1,4) circle [radius =.1];
    \draw [fill=white] (0,3) circle [radius =.1];
    \draw [fill=white] (1,4) circle [radius =.1];
\end{tikzpicture}
\xrightarrow{\jrho}
\begin{tikzpicture}[baseline={([yshift=-.8ex]current bounding box.center)}]
    \draw[-] (-1,1) -- (-1,4);
    \draw[-] (1,1) -- (1,4);
    \draw[-] (0,0) -- (0,3);

    \draw[-] (-1,1) -- (0,0);
    \draw[-] (1,1) -- (0,0);
    
    \draw [fill] (-1,1) circle [radius =.1];
    \draw [fill] (0,0) circle [radius =.1];
    \draw [fill=white] (1,1) circle [radius =.1];
    
    \draw[-] (-1,2) -- (0,1);
    \draw[-] (1,2) -- (0,1);
    
    \draw [fill] (-1,2) circle [radius =.1];
    \draw [fill] (0,1) circle [radius =.1];
    \draw [fill=white] (1,2) circle [radius =.1];
    
    \draw[-] (-1,3) -- (0,2);
    \draw[-] (1,3) -- (0,2);
    
    \draw [fill=white] (-1,3) circle [radius =.1];
    \draw [fill] (0,2) circle [radius =.1];
    \draw [fill=white] (1,3) circle [radius =.1];
    
    \draw[-] (-1,4) -- (0,3);
    \draw[-] (1,4) -- (0,3);
    
    \draw [fill=white] (-1,4) circle [radius =.1];
    \draw [fill] (0,3) circle [radius =.1];
    \draw [fill=white] (1,4) circle [radius =.1];
\end{tikzpicture}
\xrightarrow{\jrho}
\begin{tikzpicture}[baseline={([yshift=-.8ex]current bounding box.center)}]
    \draw[-] (-1,1) -- (-1,4);
    \draw[-] (1,1) -- (1,4);
    \draw[-] (0,0) -- (0,3);

    \draw[-] (-1,1) -- (0,0);
    \draw[-] (1,1) -- (0,0);
    
    \draw [fill] (-1,1) circle [radius =.1];
    \draw [fill] (0,0) circle [radius =.1];
    \draw [fill] (1,1) circle [radius =.1];
    
    \draw[-] (-1,2) -- (0,1);
    \draw[-] (1,2) -- (0,1);
    
    \draw [fill] (-1,2) circle [radius =.1];
    \draw [fill] (0,1) circle [radius =.1];
    \draw [fill=white] (1,2) circle [radius =.1];
    
    \draw[-] (-1,3) -- (0,2);
    \draw[-] (1,3) -- (0,2);
    
    \draw [fill] (-1,3) circle [radius =.1];
    \draw [fill] (0,2) circle [radius =.1];
    \draw [fill=white] (1,3) circle [radius =.1];
    
    \draw[-] (-1,4) -- (0,3);
    \draw[-] (1,4) -- (0,3);
    
    \draw [fill=white] (-1,4) circle [radius =.1];
    \draw [fill=white] (0,3) circle [radius =.1];
    \draw [fill=white] (1,4) circle [radius =.1];
\end{tikzpicture}
\xrightarrow{\jrho}
\begin{tikzpicture}[baseline={([yshift=-.8ex]current bounding box.center)}]
    \draw[-] (-1,1) -- (-1,4);
    \draw[-] (1,1) -- (1,4);
    \draw[-] (0,0) -- (0,3);

    \draw[-] (-1,1) -- (0,0);
    \draw[-] (1,1) -- (0,0);
    
    \draw [fill=white] (-1,1) circle [radius =.1];
    \draw [fill] (0,0) circle [radius =.1];
    \draw [fill] (1,1) circle [radius =.1];
    
    \draw[-] (-1,2) -- (0,1);
    \draw[-] (1,2) -- (0,1);
    
    \draw [fill=white] (-1,2) circle [radius =.1];
    \draw [fill] (0,1) circle [radius =.1];
    \draw [fill] (1,2) circle [radius =.1];
    
    \draw[-] (-1,3) -- (0,2);
    \draw[-] (1,3) -- (0,2);
    
    \draw [fill=white] (-1,3) circle [radius =.1];
    \draw [fill] (0,2) circle [radius =.1];
    \draw [fill=white] (1,3) circle [radius =.1];
    
    \draw[-] (-1,4) -- (0,3);
    \draw[-] (1,4) -- (0,3);
    
    \draw [fill=white] (-1,4) circle [radius =.1];
    \draw [fill] (0,3) circle [radius =.1];
    \draw [fill=white] (1,4) circle [radius =.1];
\end{tikzpicture}\return{\jrho}
\]
\end{ex}

To prove these theorems we utilize our equivariant bijection (Theorem~\ref{thm:eqwhirl}) from $\mathcal{J}(\VV_k)$ 
to $\mathcal{F}_k(\VV)$, 
then represent the latter by triples $f=(\ell,c,r)$ with $\ell \leq  c$ and $r\leq c$.  

This bijection $\phi$ sends an order ideal $I$ to a triple $(\ell,c,r)$, counting the number
of elements of the order ideal in the left, center, and right fibers respectively.  

\begin{ex} Here is the orbit of $\mathcal{F}_4(\VV)$ corresponding
to Example~\ref{exorbit}.  
\[(1,3,3)\xrightarrow{w} (2,4,0)\xrightarrow{w}(3,3,1)\xrightarrow{w} (0,4,2) \return{w} \]
\end{ex}

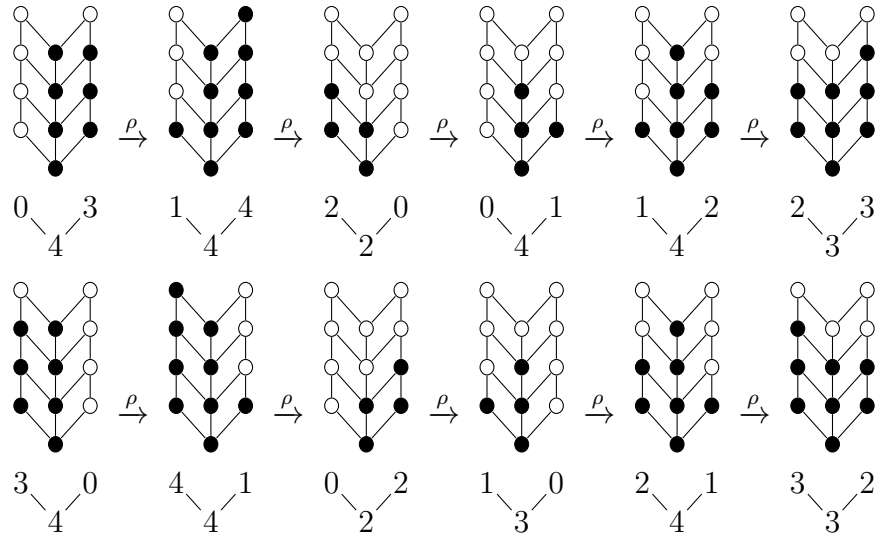
\begin{figure}
\[
\begin{tikzpicture}[scale = 0.51,xscale = 0.90, baseline={([yshift=-.8ex]current bounding box.center)}]
\node at (-1,-1) {$0$};
\node at (0,-2) {$4$};
\node at (1,-1) {$3$};
\draw [shorten >= 2mm,shorten <= 2mm] (0,-2) -- (-1,-1);
\draw [shorten >= 2mm,shorten <= 2mm] (0,-2) -- (1,-1);
%\node[draw=none,fill=none] at (0,-1) {$(0,4,3)$};

\draw [shorten >= 1.111mm,shorten <= 1.111mm] (0,0) -- (-1,1);
\draw [shorten >= 1.111mm,shorten <= 1.111mm] (0,0) -- (1,1);
\draw (-1,1) circle [radius=.2];
\draw [fill] (0,0) circle [radius=.2];
\draw [fill] (1,1) circle [radius=.2];
\draw [shorten >= 1.111mm,shorten <= 1.111mm] (0,1) -- (-1,2);
\draw [shorten >= 1.111mm,shorten <= 1.111mm] (0,1) -- (1,2);
\draw (-1,2) circle [radius=.2];
\draw [fill] (0,1) circle [radius=.2];
\draw [fill] (1,2) circle [radius=.2];
\draw [shorten >= 1.111mm,shorten <= 1.111mm] (0,2) -- (-1,3);
\draw [shorten >= 1.111mm,shorten <= 1.111mm] (0,2) -- (1,3);
\draw (-1,3) circle [radius=.2];
\draw [fill] (0,2) circle [radius=.2];
\draw [fill] (1,3) circle [radius=.2];
\draw [shorten >= 1.111mm,shorten <= 1.111mm] (0,3) -- (-1,4);
\draw [shorten >= 1.111mm,shorten <= 1.111mm] (0,3) -- (1,4);
\draw (-1,4) circle [radius=.2];
\draw [fill] (0,3) circle [radius=.2];
\draw (1,4) circle [radius=.2];
\draw [shorten >= 1.111mm,shorten <= 1.111mm] (0,0) -- (0,1);
\draw [shorten >= 1.111mm,shorten <= 1.111mm] (0,1) -- (0,2);
\draw [shorten >= 1.111mm,shorten <= 1.111mm] (0,2) -- (0,3);
\draw [shorten >= 1.111mm,shorten <= 1.111mm] (1,1) -- (1,2);
\draw [shorten >= 1.111mm,shorten <= 1.111mm] (1,2) -- (1,3);
\draw [shorten >= 1.111mm,shorten <= 1.111mm] (1,3) -- (1,4);
\draw [shorten >= 1.111mm,shorten <= 1.111mm] (-1,1) -- (-1,2);
\draw [shorten >= 1.111mm,shorten <= 1.111mm] (-1,2) -- (-1,3);
\draw [shorten >= 1.111mm,shorten <= 1.111mm] (-1,3) -- (-1,4);
\end{tikzpicture}
\xrightarrow{\rho}
\begin{tikzpicture}[scale = 0.51,xscale = 0.90, baseline={([yshift=-.8ex]current bounding box.center)}]
\node at (-1,-1) {$1$};
\node at (0,-2) {$4$};
\node at (1,-1) {$4$};
\draw [shorten >= 2mm,shorten <= 2mm] (0,-2) -- (-1,-1);
\draw [shorten >= 2mm,shorten <= 2mm] (0,-2) -- (1,-1);
%\node[draw=none,fill=none] at (0,-1) {$(1,4,4)$};
\draw [shorten >= 1.111mm,shorten <= 1.111mm] (0,0) -- (-1,1);
\draw [shorten >= 1.111mm,shorten <= 1.111mm] (0,0) -- (1,1);
\draw [fill] (-1,1) circle [radius=.2];
\draw [fill] (0,0) circle [radius=.2];
\draw [fill] (1,1) circle [radius=.2];
\draw [shorten >= 1.111mm,shorten <= 1.111mm] (0,1) -- (-1,2);
\draw [shorten >= 1.111mm,shorten <= 1.111mm] (0,1) -- (1,2);
\draw (-1,2) circle [radius=.2];
\draw [fill] (0,1) circle [radius=.2];
\draw [fill] (1,2) circle [radius=.2];
\draw [shorten >= 1.111mm,shorten <= 1.111mm] (0,2) -- (-1,3);
\draw [shorten >= 1.111mm,shorten <= 1.111mm] (0,2) -- (1,3);
\draw (-1,3) circle [radius=.2];
\draw [fill] (0,2) circle [radius=.2];
\draw [fill] (1,3) circle [radius=.2];
\draw [shorten >= 1.111mm,shorten <= 1.111mm] (0,3) -- (-1,4);
\draw [shorten >= 1.111mm,shorten <= 1.111mm] (0,3) -- (1,4);
\draw (-1,4) circle [radius=.2];
\draw [fill] (0,3) circle [radius=.2];
\draw [fill] (1,4) circle [radius=.2];
\draw [shorten >= 1.111mm,shorten <= 1.111mm] (0,0) -- (0,1);
\draw [shorten >= 1.111mm,shorten <= 1.111mm] (0,1) -- (0,2);
\draw [shorten >= 1.111mm,shorten <= 1.111mm] (0,2) -- (0,3);
\draw [shorten >= 1.111mm,shorten <= 1.111mm] (1,1) -- (1,2);
\draw [shorten >= 1.111mm,shorten <= 1.111mm] (1,2) -- (1,3);
\draw [shorten >= 1.111mm,shorten <= 1.111mm] (1,3) -- (1,4);
\draw [shorten >= 1.111mm,shorten <= 1.111mm] (-1,1) -- (-1,2);
\draw [shorten >= 1.111mm,shorten <= 1.111mm] (-1,2) -- (-1,3);
\draw [shorten >= 1.111mm,shorten <= 1.111mm] (-1,3) -- (-1,4);
\end{tikzpicture}
\xrightarrow{\rho}
\begin{tikzpicture}[scale = 0.51,xscale = 0.90, baseline={([yshift=-.8ex]current bounding box.center)}]
\node at (-1,-1) {$2$};
\node at (0,-2) {$2$};
\node at (1,-1) {$0$};
\draw [shorten >= 2mm,shorten <= 2mm] (0,-2) -- (-1,-1);
\draw [shorten >= 2mm,shorten <= 2mm] (0,-2) -- (1,-1);
%\node[draw=none,fill=none] at (0,-1) {$(2,2,0)$};
\draw [shorten >= 1.111mm,shorten <= 1.111mm] (0,0) -- (-1,1);
\draw [shorten >= 1.111mm,shorten <= 1.111mm] (0,0) -- (1,1);
\draw [fill] (-1,1) circle [radius=.2];
\draw [fill] (0,0) circle [radius=.2];
\draw (1,1) circle [radius=.2];
\draw [shorten >= 1.111mm,shorten <= 1.111mm] (0,1) -- (-1,2);
\draw [shorten >= 1.111mm,shorten <= 1.111mm] (0,1) -- (1,2);
\draw [fill] (-1,2) circle [radius=.2];
\draw [fill] (0,1) circle [radius=.2];
\draw (1,2) circle [radius=.2];
\draw [shorten >= 1.111mm,shorten <= 1.111mm] (0,2) -- (-1,3);
\draw [shorten >= 1.111mm,shorten <= 1.111mm] (0,2) -- (1,3);
\draw (-1,3) circle [radius=.2];
\draw (0,2) circle [radius=.2];
\draw (1,3) circle [radius=.2];
\draw [shorten >= 1.111mm,shorten <= 1.111mm] (0,3) -- (-1,4);
\draw [shorten >= 1.111mm,shorten <= 1.111mm] (0,3) -- (1,4);
\draw (-1,4) circle [radius=.2];
\draw (0,3) circle [radius=.2];
\draw (1,4) circle [radius=.2];
\draw [shorten >= 1.111mm,shorten <= 1.111mm] (0,0) -- (0,1);
\draw [shorten >= 1.111mm,shorten <= 1.111mm] (0,1) -- (0,2);
\draw [shorten >= 1.111mm,shorten <= 1.111mm] (0,2) -- (0,3);
\draw [shorten >= 1.111mm,shorten <= 1.111mm] (1,1) -- (1,2);
\draw [shorten >= 1.111mm,shorten <= 1.111mm] (1,2) -- (1,3);
\draw [shorten >= 1.111mm,shorten <= 1.111mm] (1,3) -- (1,4);
\draw [shorten >= 1.111mm,shorten <= 1.111mm] (-1,1) -- (-1,2);
\draw [shorten >= 1.111mm,shorten <= 1.111mm] (-1,2) -- (-1,3);
\draw [shorten >= 1.111mm,shorten <= 1.111mm] (-1,3) -- (-1,4);
\end{tikzpicture}
\xrightarrow{\rho}
\begin{tikzpicture}[scale = 0.51,xscale = 0.90, baseline={([yshift=-.8ex]current bounding box.center)}]
\node at (-1,-1) {$0$};
\node at (0,-2) {$4$};
\node at (1,-1) {$1$};
\draw [shorten >= 2mm,shorten <= 2mm] (0,-2) -- (-1,-1);
\draw [shorten >= 2mm,shorten <= 2mm] (0,-2) -- (1,-1);
%\node[draw=none,fill=none] at (0,-1) {$(0,4,1)$};
\draw [shorten >= 1.111mm,shorten <= 1.111mm] (0,0) -- (-1,1);
\draw [shorten >= 1.111mm,shorten <= 1.111mm] (0,0) -- (1,1);
\draw (-1,1) circle [radius=.2];
\draw [fill] (0,0) circle [radius=.2];
\draw [fill] (1,1) circle [radius=.2];
\draw [shorten >= 1.111mm,shorten <= 1.111mm] (0,1) -- (-1,2);
\draw [shorten >= 1.111mm,shorten <= 1.111mm] (0,1) -- (1,2);
\draw (-1,2) circle [radius=.2];
\draw [fill] (0,1) circle [radius=.2];
\draw (1,2) circle [radius=.2];
\draw [shorten >= 1.111mm,shorten <= 1.111mm] (0,2) -- (-1,3);
\draw [shorten >= 1.111mm,shorten <= 1.111mm] (0,2) -- (1,3);
\draw (-1,3) circle [radius=.2];
\draw [fill] (0,2) circle [radius=.2];
\draw (1,3) circle [radius=.2];
\draw [shorten >= 1.111mm,shorten <= 1.111mm] (0,3) -- (-1,4);
\draw [shorten >= 1.111mm,shorten <= 1.111mm] (0,3) -- (1,4);
\draw (-1,4) circle [radius=.2];
\draw (0,3) circle [radius=.2];
\draw (1,4) circle [radius=.2];
\draw [shorten >= 1.111mm,shorten <= 1.111mm] (0,0) -- (0,1);
\draw [shorten >= 1.111mm,shorten <= 1.111mm] (0,1) -- (0,2);
\draw [shorten >= 1.111mm,shorten <= 1.111mm] (0,2) -- (0,3);
\draw [shorten >= 1.111mm,shorten <= 1.111mm] (1,1) -- (1,2);
\draw [shorten >= 1.111mm,shorten <= 1.111mm] (1,2) -- (1,3);
\draw [shorten >= 1.111mm,shorten <= 1.111mm] (1,3) -- (1,4);
\draw [shorten >= 1.111mm,shorten <= 1.111mm] (-1,1) -- (-1,2);
\draw [shorten >= 1.111mm,shorten <= 1.111mm] (-1,2) -- (-1,3);
\draw [shorten >= 1.111mm,shorten <= 1.111mm] (-1,3) -- (-1,4);
\end{tikzpicture}
\xrightarrow{\rho}
\begin{tikzpicture}[scale = 0.51,xscale = 0.90, baseline={([yshift=-.8ex]current bounding box.center)}]
\node at (-1,-1) {$1$};
\node at (0,-2) {$4$};
\node at (1,-1) {$2$};
\draw [shorten >= 2mm,shorten <= 2mm] (0,-2) -- (-1,-1);
\draw [shorten >= 2mm,shorten <= 2mm] (0,-2) -- (1,-1);
%\node[draw=none,fill=none] at (0,-1) {$(1,4,2)$};
\draw [shorten >= 1.111mm,shorten <= 1.111mm] (0,0) -- (-1,1);
\draw [shorten >= 1.111mm,shorten <= 1.111mm] (0,0) -- (1,1);
\draw [fill] (-1,1) circle [radius=.2];
\draw [fill] (0,0) circle [radius=.2];
\draw [fill] (1,1) circle [radius=.2];
\draw [shorten >= 1.111mm,shorten <= 1.111mm] (0,1) -- (-1,2);
\draw [shorten >= 1.111mm,shorten <= 1.111mm] (0,1) -- (1,2);
\draw (-1,2) circle [radius=.2];
\draw [fill] (0,1) circle [radius=.2];
\draw [fill] (1,2) circle [radius=.2];
\draw [shorten >= 1.111mm,shorten <= 1.111mm] (0,2) -- (-1,3);
\draw [shorten >= 1.111mm,shorten <= 1.111mm] (0,2) -- (1,3);
\draw (-1,3) circle [radius=.2];
\draw [fill] (0,2) circle [radius=.2];
\draw (1,3) circle [radius=.2];
\draw [shorten >= 1.111mm,shorten <= 1.111mm] (0,3) -- (-1,4);
\draw [shorten >= 1.111mm,shorten <= 1.111mm] (0,3) -- (1,4);
\draw (-1,4) circle [radius=.2];
\draw [fill] (0,3) circle [radius=.2];
\draw (1,4) circle [radius=.2];
\draw [shorten >= 1.111mm,shorten <= 1.111mm] (0,0) -- (0,1);
\draw [shorten >= 1.111mm,shorten <= 1.111mm] (0,1) -- (0,2);
\draw [shorten >= 1.111mm,shorten <= 1.111mm] (0,2) -- (0,3);
\draw [shorten >= 1.111mm,shorten <= 1.111mm] (1,1) -- (1,2);
\draw [shorten >= 1.111mm,shorten <= 1.111mm] (1,2) -- (1,3);
\draw [shorten >= 1.111mm,shorten <= 1.111mm] (1,3) -- (1,4);
\draw [shorten >= 1.111mm,shorten <= 1.111mm] (-1,1) -- (-1,2);
\draw [shorten >= 1.111mm,shorten <= 1.111mm] (-1,2) -- (-1,3);
\draw [shorten >= 1.111mm,shorten <= 1.111mm] (-1,3) -- (-1,4);
\end{tikzpicture}
\xrightarrow{\rho}
\begin{tikzpicture}[scale = 0.51,xscale = 0.90, baseline={([yshift=-.8ex]current bounding box.center)}]
\node at (-1,-1) {$2$};
\node at (0,-2) {$3$};
\node at (1,-1) {$3$};
\draw [shorten >= 2mm,shorten <= 2mm] (0,-2) -- (-1,-1);
\draw [shorten >= 2mm,shorten <= 2mm] (0,-2) -- (1,-1);
%\node[draw=none,fill=none] at (0,-1) {$(2,3,3)$};
\draw [shorten >= 1.111mm,shorten <= 1.111mm] (0,0) -- (-1,1);
\draw [shorten >= 1.111mm,shorten <= 1.111mm] (0,0) -- (1,1);
\draw [fill] (-1,1) circle [radius=.2];
\draw [fill] (0,0) circle [radius=.2];
\draw [fill] (1,1) circle [radius=.2];
\draw [shorten >= 1.111mm,shorten <= 1.111mm] (0,1) -- (-1,2);
\draw [shorten >= 1.111mm,shorten <= 1.111mm] (0,1) -- (1,2);
\draw [fill] (-1,2) circle [radius=.2];
\draw [fill] (0,1) circle [radius=.2];
\draw [fill] (1,2) circle [radius=.2];
\draw [shorten >= 1.111mm,shorten <= 1.111mm] (0,2) -- (-1,3);
\draw [shorten >= 1.111mm,shorten <= 1.111mm] (0,2) -- (1,3);
\draw (-1,3) circle [radius=.2];
\draw [fill] (0,2) circle [radius=.2];
\draw [fill] (1,3) circle [radius=.2];
\draw [shorten >= 1.111mm,shorten <= 1.111mm] (0,3) -- (-1,4);
\draw [shorten >= 1.111mm,shorten <= 1.111mm] (0,3) -- (1,4);
\draw (-1,4) circle [radius=.2];
\draw (0,3) circle [radius=.2];
\draw (1,4) circle [radius=.2];
\draw [shorten >= 1.111mm,shorten <= 1.111mm] (0,0) -- (0,1);
\draw [shorten >= 1.111mm,shorten <= 1.111mm] (0,1) -- (0,2);
\draw [shorten >= 1.111mm,shorten <= 1.111mm] (0,2) -- (0,3);
\draw [shorten >= 1.111mm,shorten <= 1.111mm] (1,1) -- (1,2);
\draw [shorten >= 1.111mm,shorten <= 1.111mm] (1,2) -- (1,3);
\draw [shorten >= 1.111mm,shorten <= 1.111mm] (1,3) -- (1,4);
\draw [shorten >= 1.111mm,shorten <= 1.111mm] (-1,1) -- (-1,2);
\draw [shorten >= 1.111mm,shorten <= 1.111mm] (-1,2) -- (-1,3);
\draw [shorten >= 1.111mm,shorten <= 1.111mm] (-1,3) -- (-1,4);
\end{tikzpicture}
\]
\[
\begin{tikzpicture}[scale = 0.51,xscale = 0.90, baseline={([yshift=-.8ex]current bounding box.center)}]
\node at (-1,-1) {$3$};
\node at (0,-2) {$4$};
\node at (1,-1) {$0$};
\draw [shorten >= 2mm,shorten <= 2mm] (0,-2) -- (-1,-1);
\draw [shorten >= 2mm,shorten <= 2mm] (0,-2) -- (1,-1);
%\node[draw=none,fill=none] at (0,-1) {$(3,4,0)$};
\draw [shorten >= 1.111mm,shorten <= 1.111mm] (0,0) -- (-1,1);
\draw [shorten >= 1.111mm,shorten <= 1.111mm] (0,0) -- (1,1);
\draw [fill] (-1,1) circle [radius=.2];
\draw [fill] (0,0) circle [radius=.2];
\draw (1,1) circle [radius=.2];
\draw [shorten >= 1.111mm,shorten <= 1.111mm] (0,1) -- (-1,2);
\draw [shorten >= 1.111mm,shorten <= 1.111mm] (0,1) -- (1,2);
\draw [fill] (-1,2) circle [radius=.2];
\draw [fill] (0,1) circle [radius=.2];
\draw(1,2) circle [radius=.2];
\draw [shorten >= 1.111mm,shorten <= 1.111mm] (0,2) -- (-1,3);
\draw [shorten >= 1.111mm,shorten <= 1.111mm] (0,2) -- (1,3);
\draw [fill] (-1,3) circle [radius=.2];
\draw [fill] (0,2) circle [radius=.2];
\draw(1,3) circle [radius=.2];
\draw [shorten >= 1.111mm,shorten <= 1.111mm] (0,3) -- (-1,4);
\draw [shorten >= 1.111mm,shorten <= 1.111mm] (0,3) -- (1,4);
\draw(-1,4) circle [radius=.2];
\draw [fill] (0,3) circle [radius=.2];
\draw(1,4) circle [radius=.2];
\draw [shorten >= 1.111mm,shorten <= 1.111mm] (0,0) -- (0,1);
\draw [shorten >= 1.111mm,shorten <= 1.111mm] (0,1) -- (0,2);
\draw [shorten >= 1.111mm,shorten <= 1.111mm] (0,2) -- (0,3);
\draw [shorten >= 1.111mm,shorten <= 1.111mm] (1,1) -- (1,2);
\draw [shorten >= 1.111mm,shorten <= 1.111mm] (1,2) -- (1,3);
\draw [shorten >= 1.111mm,shorten <= 1.111mm] (1,3) -- (1,4);
\draw [shorten >= 1.111mm,shorten <= 1.111mm] (-1,1) -- (-1,2);
\draw [shorten >= 1.111mm,shorten <= 1.111mm] (-1,2) -- (-1,3);
\draw [shorten >= 1.111mm,shorten <= 1.111mm] (-1,3) -- (-1,4);
\end{tikzpicture}
\xrightarrow{\rho}
\begin{tikzpicture}[scale = 0.51,xscale = 0.90, baseline={([yshift=-.8ex]current bounding box.center)}]
\node at (-1,-1) {$4$};
\node at (0,-2) {$4$};
\node at (1,-1) {$1$};
\draw [shorten >= 2mm,shorten <= 2mm] (0,-2) -- (-1,-1);
\draw [shorten >= 2mm,shorten <= 2mm] (0,-2) -- (1,-1);
%\node[draw=none,fill=none] at (0,-1) {$(4,4,1)$};
\draw [shorten >= 1.111mm,shorten <= 1.111mm] (0,0) -- (-1,1);
\draw [shorten >= 1.111mm,shorten <= 1.111mm] (0,0) -- (1,1);
\draw [fill] (-1,1) circle [radius=.2];
\draw [fill] (0,0) circle [radius=.2];
\draw [fill] (1,1) circle [radius=.2];
\draw [shorten >= 1.111mm,shorten <= 1.111mm] (0,1) -- (-1,2);
\draw [shorten >= 1.111mm,shorten <= 1.111mm] (0,1) -- (1,2);
\draw [fill] (-1,2) circle [radius=.2];
\draw [fill] (0,1) circle [radius=.2];
\draw(1,2) circle [radius=.2];
\draw [shorten >= 1.111mm,shorten <= 1.111mm] (0,2) -- (-1,3);
\draw [shorten >= 1.111mm,shorten <= 1.111mm] (0,2) -- (1,3);
\draw [fill] (-1,3) circle [radius=.2];
\draw [fill] (0,2) circle [radius=.2];
\draw(1,3) circle [radius=.2];
\draw [shorten >= 1.111mm,shorten <= 1.111mm] (0,3) -- (-1,4);
\draw [shorten >= 1.111mm,shorten <= 1.111mm] (0,3) -- (1,4);
\draw [fill] (-1,4) circle [radius=.2];
\draw [fill] (0,3) circle [radius=.2];
\draw (1,4) circle [radius=.2];
\draw [shorten >= 1.111mm,shorten <= 1.111mm] (0,0) -- (0,1);
\draw [shorten >= 1.111mm,shorten <= 1.111mm] (0,1) -- (0,2);
\draw [shorten >= 1.111mm,shorten <= 1.111mm] (0,2) -- (0,3);
\draw [shorten >= 1.111mm,shorten <= 1.111mm] (1,1) -- (1,2);
\draw [shorten >= 1.111mm,shorten <= 1.111mm] (1,2) -- (1,3);
\draw [shorten >= 1.111mm,shorten <= 1.111mm] (1,3) -- (1,4);
\draw [shorten >= 1.111mm,shorten <= 1.111mm] (-1,1) -- (-1,2);
\draw [shorten >= 1.111mm,shorten <= 1.111mm] (-1,2) -- (-1,3);
\draw [shorten >= 1.111mm,shorten <= 1.111mm] (-1,3) -- (-1,4);
\end{tikzpicture}
\xrightarrow{\rho}
\begin{tikzpicture}[scale = 0.51,xscale = 0.90, baseline={([yshift=-.8ex]current bounding box.center)}]
\node at (-1,-1) {$0$};
\node at (0,-2) {$2$};
\node at (1,-1) {$2$};
\draw [shorten >= 2mm,shorten <= 2mm] (0,-2) -- (-1,-1);
\draw [shorten >= 2mm,shorten <= 2mm] (0,-2) -- (1,-1);
%\node[draw=none,fill=none] at (0,-1) {$(0,2,2)$};
\draw [shorten >= 1.111mm,shorten <= 1.111mm] (0,0) -- (-1,1);
\draw [shorten >= 1.111mm,shorten <= 1.111mm] (0,0) -- (1,1);
\draw (-1,1) circle [radius=.2];
\draw [fill] (0,0) circle [radius=.2];
\draw [fill] (1,1) circle [radius=.2];
\draw [shorten >= 1.111mm,shorten <= 1.111mm] (0,1) -- (-1,2);
\draw [shorten >= 1.111mm,shorten <= 1.111mm] (0,1) -- (1,2);
\draw (-1,2) circle [radius=.2];
\draw [fill] (0,1) circle [radius=.2];
\draw [fill] (1,2) circle [radius=.2];
\draw [shorten >= 1.111mm,shorten <= 1.111mm] (0,2) -- (-1,3);
\draw [shorten >= 1.111mm,shorten <= 1.111mm] (0,2) -- (1,3);
\draw (-1,3) circle [radius=.2];
\draw (0,2) circle [radius=.2];
\draw (1,3) circle [radius=.2];
\draw [shorten >= 1.111mm,shorten <= 1.111mm] (0,3) -- (-1,4);
\draw [shorten >= 1.111mm,shorten <= 1.111mm] (0,3) -- (1,4);
\draw (-1,4) circle [radius=.2];
\draw (0,3) circle [radius=.2];
\draw (1,4) circle [radius=.2];
\draw [shorten >= 1.111mm,shorten <= 1.111mm] (0,0) -- (0,1);
\draw [shorten >= 1.111mm,shorten <= 1.111mm] (0,1) -- (0,2);
\draw [shorten >= 1.111mm,shorten <= 1.111mm] (0,2) -- (0,3);
\draw [shorten >= 1.111mm,shorten <= 1.111mm] (1,1) -- (1,2);
\draw [shorten >= 1.111mm,shorten <= 1.111mm] (1,2) -- (1,3);
\draw [shorten >= 1.111mm,shorten <= 1.111mm] (1,3) -- (1,4);
\draw [shorten >= 1.111mm,shorten <= 1.111mm] (-1,1) -- (-1,2);
\draw [shorten >= 1.111mm,shorten <= 1.111mm] (-1,2) -- (-1,3);
\draw [shorten >= 1.111mm,shorten <= 1.111mm] (-1,3) -- (-1,4);
\end{tikzpicture}
\xrightarrow{\rho}
\begin{tikzpicture}[scale = 0.51,xscale = 0.90, baseline={([yshift=-.8ex]current bounding box.center)}]
\node at (-1,-1) {$1$};
\node at (0,-2) {$3$};
\node at (1,-1) {$0$};
\draw [shorten >= 2mm,shorten <= 2mm] (0,-2) -- (-1,-1);
\draw [shorten >= 2mm,shorten <= 2mm] (0,-2) -- (1,-1);
%\node[draw=none,fill=none] at (0,-1) {$(1,3,0)$};
\draw [shorten >= 1.111mm,shorten <= 1.111mm] (0,0) -- (-1,1);
\draw [shorten >= 1.111mm,shorten <= 1.111mm] (0,0) -- (1,1);
\draw [fill] (-1,1) circle [radius=.2];
\draw [fill] (0,0) circle [radius=.2];
\draw (1,1) circle [radius=.2];
\draw [shorten >= 1.111mm,shorten <= 1.111mm] (0,1) -- (-1,2);
\draw [shorten >= 1.111mm,shorten <= 1.111mm] (0,1) -- (1,2);
\draw (-1,2) circle [radius=.2];
\draw [fill] (0,1) circle [radius=.2];
\draw (1,2) circle [radius=.2];
\draw [shorten >= 1.111mm,shorten <= 1.111mm] (0,2) -- (-1,3);
\draw [shorten >= 1.111mm,shorten <= 1.111mm] (0,2) -- (1,3);
\draw (-1,3) circle [radius=.2];
\draw [fill] (0,2) circle [radius=.2];
\draw (1,3) circle [radius=.2];
\draw [shorten >= 1.111mm,shorten <= 1.111mm] (0,3) -- (-1,4);
\draw [shorten >= 1.111mm,shorten <= 1.111mm] (0,3) -- (1,4);
\draw (-1,4) circle [radius=.2];
\draw (0,3) circle [radius=.2];
\draw (1,4) circle [radius=.2];
\draw [shorten >= 1.111mm,shorten <= 1.111mm] (0,0) -- (0,1);
\draw [shorten >= 1.111mm,shorten <= 1.111mm] (0,1) -- (0,2);
\draw [shorten >= 1.111mm,shorten <= 1.111mm] (0,2) -- (0,3);
\draw [shorten >= 1.111mm,shorten <= 1.111mm] (1,1) -- (1,2);
\draw [shorten >= 1.111mm,shorten <= 1.111mm] (1,2) -- (1,3);
\draw [shorten >= 1.111mm,shorten <= 1.111mm] (1,3) -- (1,4);
\draw [shorten >= 1.111mm,shorten <= 1.111mm] (-1,1) -- (-1,2);
\draw [shorten >= 1.111mm,shorten <= 1.111mm] (-1,2) -- (-1,3);
\draw [shorten >= 1.111mm,shorten <= 1.111mm] (-1,3) -- (-1,4);
\end{tikzpicture}
\xrightarrow{\rho}
\begin{tikzpicture}[scale = 0.51,xscale = 0.90, baseline={([yshift=-.8ex]current bounding box.center)}]
\node at (-1,-1) {$2$};
\node at (0,-2) {$4$};
\node at (1,-1) {$1$};
\draw [shorten >= 2mm,shorten <= 2mm] (0,-2) -- (-1,-1);
\draw [shorten >= 2mm,shorten <= 2mm] (0,-2) -- (1,-1);
%\node[draw=none,fill=none] at (0,-1) {$(2,4,1)$};
\draw [shorten >= 1.111mm,shorten <= 1.111mm] (0,0) -- (-1,1);
\draw [shorten >= 1.111mm,shorten <= 1.111mm] (0,0) -- (1,1);
\draw [fill] (-1,1) circle [radius=.2];
\draw [fill] (0,0) circle [radius=.2];
\draw [fill] (1,1) circle [radius=.2];
\draw [shorten >= 1.111mm,shorten <= 1.111mm] (0,1) -- (-1,2);
\draw [shorten >= 1.111mm,shorten <= 1.111mm] (0,1) -- (1,2);
\draw [fill] (-1,2) circle [radius=.2];
\draw [fill] (0,1) circle [radius=.2];
\draw (1,2) circle [radius=.2];
\draw [shorten >= 1.111mm,shorten <= 1.111mm] (0,2) -- (-1,3);
\draw [shorten >= 1.111mm,shorten <= 1.111mm] (0,2) -- (1,3);
\draw (-1,3) circle [radius=.2];
\draw [fill] (0,2) circle [radius=.2];
\draw (1,3) circle [radius=.2];
\draw [shorten >= 1.111mm,shorten <= 1.111mm] (0,3) -- (-1,4);
\draw [shorten >= 1.111mm,shorten <= 1.111mm] (0,3) -- (1,4);
\draw (-1,4) circle [radius=.2];
\draw [fill] (0,3) circle [radius=.2];
\draw (1,4) circle [radius=.2];
\draw [shorten >= 1.111mm,shorten <= 1.111mm] (0,0) -- (0,1);
\draw [shorten >= 1.111mm,shorten <= 1.111mm] (0,1) -- (0,2);
\draw [shorten >= 1.111mm,shorten <= 1.111mm] (0,2) -- (0,3);
\draw [shorten >= 1.111mm,shorten <= 1.111mm] (1,1) -- (1,2);
\draw [shorten >= 1.111mm,shorten <= 1.111mm] (1,2) -- (1,3);
\draw [shorten >= 1.111mm,shorten <= 1.111mm] (1,3) -- (1,4);
\draw [shorten >= 1.111mm,shorten <= 1.111mm] (-1,1) -- (-1,2);
\draw [shorten >= 1.111mm,shorten <= 1.111mm] (-1,2) -- (-1,3);
\draw [shorten >= 1.111mm,shorten <= 1.111mm] (-1,3) -- (-1,4);
\end{tikzpicture}
\xrightarrow{\rho}
\begin{tikzpicture}[scale = 0.51,xscale = 0.90, baseline={([yshift=-.8ex]current bounding box.center)}]
\node at (-1,-1) {$3$};
\node at (0,-2) {$3$};
\node at (1,-1) {$2$};
\draw [shorten >= 2mm,shorten <= 2mm] (0,-2) -- (-1,-1);
\draw [shorten >= 2mm,shorten <= 2mm] (0,-2) -- (1,-1);
%\node[draw=none,fill=none] at (0,-1) {$(3,3,2)$};
\draw [shorten >= 1.111mm,shorten <= 1.111mm] (0,0) -- (-1,1);
\draw [shorten >= 1.111mm,shorten <= 1.111mm] (0,0) -- (1,1);
\draw [fill] (-1,1) circle [radius=.2];
\draw [fill] (0,0) circle [radius=.2];
\draw [fill] (1,1) circle [radius=.2];
\draw [shorten >= 1.111mm,shorten <= 1.111mm] (0,1) -- (-1,2);
\draw [shorten >= 1.111mm,shorten <= 1.111mm] (0,1) -- (1,2);
\draw [fill] (-1,2) circle [radius=.2];
\draw [fill] (0,1) circle [radius=.2];
\draw [fill] (1,2) circle [radius=.2];
\draw [shorten >= 1.111mm,shorten <= 1.111mm] (0,2) -- (-1,3);
\draw [shorten >= 1.111mm,shorten <= 1.111mm] (0,2) -- (1,3);
\draw [fill] (-1,3) circle [radius=.2];
\draw [fill] (0,2) circle [radius=.2];
\draw (1,3) circle [radius=.2];
\draw [shorten >= 1.111mm,shorten <= 1.111mm] (0,3) -- (-1,4);
\draw [shorten >= 1.111mm,shorten <= 1.111mm] (0,3) -- (1,4);
\draw (-1,4) circle [radius=.2];
\draw (0,3) circle [radius=.2];
\draw (1,4) circle [radius=.2];
\draw [shorten >= 1.111mm,shorten <= 1.111mm] (0,0) -- (0,1);
\draw [shorten >= 1.111mm,shorten <= 1.111mm] (0,1) -- (0,2);
\draw [shorten >= 1.111mm,shorten <= 1.111mm] (0,2) -- (0,3);
\draw [shorten >= 1.111mm,shorten <= 1.111mm] (1,1) -- (1,2);
\draw [shorten >= 1.111mm,shorten <= 1.111mm] (1,2) -- (1,3);
\draw [shorten >= 1.111mm,shorten <= 1.111mm] (1,3) -- (1,4);
\draw [shorten >= 1.111mm,shorten <= 1.111mm] (-1,1) -- (-1,2);
\draw [shorten >= 1.111mm,shorten <= 1.111mm] (-1,2) -- (-1,3);
\draw [shorten >= 1.111mm,shorten <= 1.111mm] (-1,3) -- (-1,4);
\end{tikzpicture}
\]
\caption{An orbit of rowmotion on $P=\VV \times [4]$, with corresponding $P$-partitions (which
are being whirled), as in Theorem~\ref{thm:eqwhirl}.}  
\label{fig:rowWhirl}
\end{figure}
Figure~\ref{fig:rowWhirl} shows an orbit of rowmotion on order ideals of $P=\VV \times [4]$.
Below each ideal is the corresponding $P$-partition in $\calF_{4}(P)$, and these are being whirled
equivariantly with $\rho $.  

\begin{prop}\label{VCounting}
The number of order ideals of $\VV_k$ is given by 
$\displaystyle \sum_{i=1}^{k+1} i^2 = \frac{(k+1)(k+2)(2k+3)}{6}$.
%$|\mathcal{J}(\VV_k)| = \frac{k(k+1)(2k+1)}{6}$. 
\end{prop}
\begin{proof}
Every order ideal of $\VV_k$ is associated with a triple of labels of the elements of $\VV$.
Let $j \in [0, k]$ denote the label of $c$.
The possible labels for $\ell$ and $r$ are precisely the elements of $[0, j]$, yielding a total of $(j+1)^2$ labelings.
Summing over all possible values of $j$ produces the desired result.
%We first show that the number of triples $(a_1,a_2,a_3)\in [0,k+1]^3$ with $a_1,a_3<a_2$
%(strict inequality) is given by $\frac{k(k+1)(2k+1)}{6}$. 
%If $a_1=a_3< a_2$ we have ${k+1\choose 2} = \frac{(k+1)(k)}{2}$ such triples. If $a_1<a_3<a_2$ or $a_3<a_1<a_2$, then we have
%$ 2\cdot {k+1\choose 3} =\frac{(k+1)k(k-1)}{3}$ such triples.  These sum to 
%$\frac{k(k+1)(2k+1)}{6}$. These triples are in bijection with those defining order ideals of
%$\VV_{k}$ by simply subtracting 1 from the middle component, making the strict inequality
%with $a_{2}$ weak.   
\end{proof}

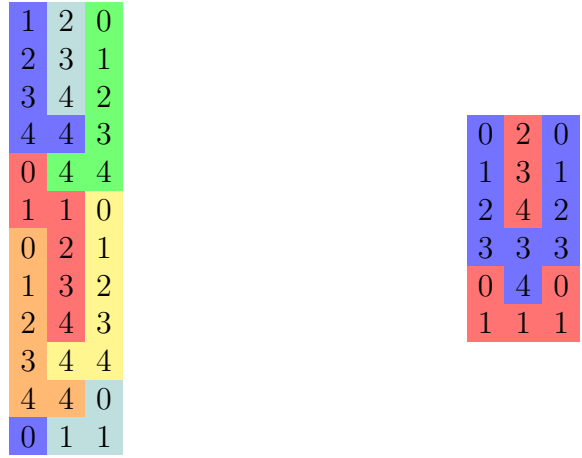
\begin{figure}
\hfill
\begin{tikzpicture}[baseline={([yshift=-.8ex]current bounding box.center)}]
    \tikzColoredSquare{0}{0}{1}{white!45!blue}
    \tikzColoredSquare{0}{-1}{2}{white!45!blue}
    \tikzColoredSquare{0}{-2}{3}{white!45!blue}
    \tikzColoredSquare{0}{-3}{4}{white!45!blue}
    \tikzColoredSquare{1}{-3}{4}{white!45!blue}
    \tikzColoredSquare{2}{0}{0}{white!45!green}
    \tikzColoredSquare{2}{-1}{1}{white!45!green}
    \tikzColoredSquare{2}{-2}{2}{white!45!green}
    \tikzColoredSquare{2}{-3}{3}{white!45!green}
    \tikzColoredSquare{2}{-4}{4}{white!45!green}
    \tikzColoredSquare{1}{-4}{4}{white!45!green}
    \tikzColoredSquare{0}{-4}{0}{white!45!red}
    \tikzColoredSquare{0}{-5}{1}{white!45!red}
    \tikzColoredSquare{1}{-5}{1}{white!45!red}
    \tikzColoredSquare{1}{-6}{2}{white!45!red}
    \tikzColoredSquare{1}{-7}{3}{white!45!red}
    \tikzColoredSquare{1}{-8}{4}{white!45!red}
    \tikzColoredSquare{2}{-5}{0}{white!45!yellow}
    \tikzColoredSquare{2}{-6}{1}{white!45!yellow}
    \tikzColoredSquare{2}{-7}{2}{white!45!yellow}
    \tikzColoredSquare{2}{-8}{3}{white!45!yellow}
    \tikzColoredSquare{2}{-9}{4}{white!45!yellow}
    \tikzColoredSquare{1}{-9}{4}{white!45!yellow}
    \tikzColoredSquare{0}{-6}{0}{white!45!orange}
    \tikzColoredSquare{0}{-7}{1}{white!45!orange}
    \tikzColoredSquare{0}{-8}{2}{white!45!orange}
    \tikzColoredSquare{0}{-9}{3}{white!45!orange}
    \tikzColoredSquare{0}{-10}{4}{white!45!orange}
    \tikzColoredSquare{1}{-10}{4}{white!45!orange}
    \tikzColoredSquare{2}{-10}{0}{white!75!teal}
    \tikzColoredSquare{2}{-11}{1}{white!75!teal}
    \tikzColoredSquare{1}{-11}{1}{white!75!teal}
    \tikzColoredSquare{1}{0}{2}{white!75!teal}
    \tikzColoredSquare{1}{-1}{3}{white!75!teal}
    \tikzColoredSquare{1}{-2}{4}{white!75!teal}
    \tikzColoredSquare{0}{-11}{0}{white!45!blue}
\end{tikzpicture}
\hfill
\begin{tikzpicture}[baseline={([yshift=-.8ex]current bounding box.center)}]

    \def\m{0}
    \tikzColoredSquare{0}{\m}{0}{white!45!blue}
    \tikzColoredSquare{1}{\m}{2}{white!45!red}
    \tikzColoredSquare{2}{\m}{0}{white!45!blue}
    
    \def\m{-1}
    \tikzColoredSquare{0}{\m}{1}{white!45!blue}
    \tikzColoredSquare{1}{\m}{3}{white!45!red}
    \tikzColoredSquare{2}{\m}{1}{white!45!blue}
    
    \def\m{-2}
    \tikzColoredSquare{0}{\m}{2}{white!45!blue}
    \tikzColoredSquare{1}{\m}{4}{white!45!red}
    \tikzColoredSquare{2}{\m}{2}{white!45!blue}
    
    \def\m{-3}
    \tikzColoredSquare{0}{\m}{3}{white!45!blue}
    \tikzColoredSquare{1}{\m}{3}{white!45!blue}
    \tikzColoredSquare{2}{\m}{3}{white!45!blue}
    
    \def\m{-4}
    \tikzColoredSquare{0}{\m}{0}{white!45!red}
    \tikzColoredSquare{1}{\m}{4}{white!45!blue}
    \tikzColoredSquare{2}{\m}{0}{white!45!red}
    
    \def\m{-5}
    \tikzColoredSquare{0}{\m}{1}{white!45!red}
    \tikzColoredSquare{1}{\m}{4}{white!45!blue}
    \tikzColoredSquare{2}{\m}{1}{white!45!red}
    
    \def\m{-5}
    \tikzColoredSquare{0}{\m}{1}{white!45!red}
    \tikzColoredSquare{1}{\m}{1}{white!45!red}
    \tikzColoredSquare{2}{\m}{1}{white!45!red}
\end{tikzpicture}
\hfill
\caption{Two orbit boards
of $\mathcal{F}_4(\VV)$, one with six whorms and one with two whorms.}\label{fig:veetwoorbits}
\label{fig:whormboard}
\end{figure}

\subsection{Periodicity and homomesy via center-seeking whorms}
To show that the order of $\jrho$ on $\mathcal{J}({\VV }_k)$ is $2(k+2)$ we end up proving
something stronger, namely that $\jrho^{k+2}(I)$ is the reflection of $I$ across the
center chain. Our method is to investigate the whorms that arise from repeatedly
whirling a $k$-bounded $P$-partition. 

Recall from Definition~\ref{def:whorms} that, given a whirling orbit board,
$\mathcal{R}=\{f,w(f), w^2(f),\dots\}$ of $w$ on $\mathcal{F}_k(\VV)$, a whorm $\xi$ is a
maximal set of whorm-connected elements. Figure~\ref{fig:veetwoorbits} shows two orbit boards
of $\mathcal{F}_4(\VV)$, one with six whorms and one with two whorms. Notice that each
whorm in the second orbit has two ``starting" positions.

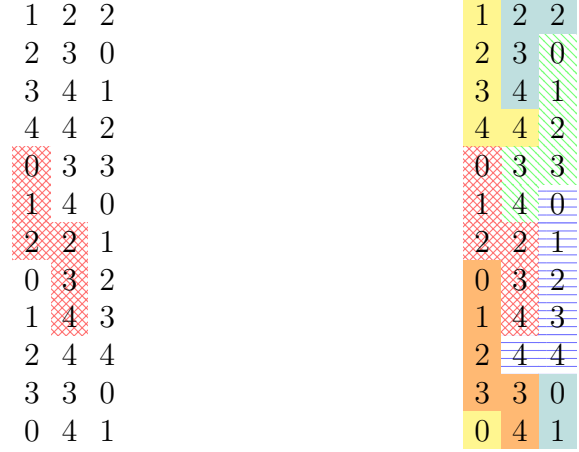
\begin{figure}
\hfill
\begin{tikzpicture}[baseline={([yshift=-.8ex]current bounding box.center)}]
    \tikzColoredSquare{0}{0}{1}{white}
    \tikzColoredSquare{0}{-1}{2}{white}
    \tikzColoredSquare{0}{-2}{3}{white}
    \tikzColoredSquare{0}{-3}{4}{white}
    \tikzColoredSquare{1}{-3}{4}{white}
    \tikzColoredSquare{2}{-1}{0}{white}
    \tikzColoredSquare{2}{-2}{1}{white}
    \tikzColoredSquare{2}{-3}{2}{white}
    \tikzColoredSquare{2}{-4}{3}{white}
    \tikzColoredSquare{1}{-4}{3}{white}
    \tikzColoredSquare{1}{-5}{4}{white}
    \tikzColoredSquarech{0}{-4}{0}{white!45!red}
    \tikzColoredSquarech{0}{-5}{1}{white!45!red}
    \tikzColoredSquarech{0}{-6}{2}{white!45!red}
    \tikzColoredSquarech{1}{-6}{2}{white!45!red}
    \tikzColoredSquarech{1}{-7}{3}{white!45!red}
    \tikzColoredSquarech{1}{-8}{4}{white!45!red}
    \tikzColoredSquare{2}{-5}{0}{white}
    \tikzColoredSquare{2}{-6}{1}{white}
    \tikzColoredSquare{2}{-7}{2}{white}
    \tikzColoredSquare{2}{-8}{3}{white}
    \tikzColoredSquare{2}{-9}{4}{white}
    \tikzColoredSquare{1}{-9}{4}{white}
    \tikzColoredSquare{0}{-7}{0}{white}
    \tikzColoredSquare{0}{-8}{1}{white}
    \tikzColoredSquare{0}{-9}{2}{white}
    \tikzColoredSquare{0}{-10}{3}{white}
    \tikzColoredSquare{1}{-10}{3}{white}
    \tikzColoredSquare{1}{-11}{4}{white}
    \tikzColoredSquare{2}{-10}{0}{white}
    \tikzColoredSquare{2}{-11}{1}{white}
    \tikzColoredSquare{2}{0}{2}{white}
    \tikzColoredSquare{1}{0}{2}{white}
    \tikzColoredSquare{1}{-1}{3}{white}
    \tikzColoredSquare{1}{-2}{4}{white}
    \tikzColoredSquare{0}{-11}{0}{white}
\end{tikzpicture}\hfill
\begin{tikzpicture}[baseline={([yshift=-.8ex]current bounding box.center)}]
    \tikzColoredSquare{0}{0}{1}{white!45!yellow}
    \tikzColoredSquare{0}{-1}{2}{white!45!yellow}
    \tikzColoredSquare{0}{-2}{3}{white!45!yellow}
    \tikzColoredSquare{0}{-3}{4}{white!45!yellow}
    \tikzColoredSquare{1}{-3}{4}{white!45!yellow}
    \tikzColoredSquarenwl{2}{-1}{0}{white!45!green}
    \tikzColoredSquarenwl{2}{-2}{1}{white!45!green}
    \tikzColoredSquarenwl{2}{-3}{2}{white!45!green}
    \tikzColoredSquarenwl{2}{-4}{3}{white!45!green}
    \tikzColoredSquarenwl{1}{-4}{3}{white!45!green}
    \tikzColoredSquarenwl{1}{-5}{4}{white!45!green}
    \tikzColoredSquarech{0}{-4}{0}{white!45!red}
    \tikzColoredSquarech{0}{-5}{1}{white!45!red}
    \tikzColoredSquarech{0}{-6}{2}{white!45!red}
    \tikzColoredSquarech{1}{-6}{2}{white!45!red}
    \tikzColoredSquarech{1}{-7}{3}{white!45!red}
    \tikzColoredSquarech{1}{-8}{4}{white!45!red}
    \tikzColoredSquarehl{2}{-5}{0}{white!45!blue}
    \tikzColoredSquarehl{2}{-6}{1}{white!45!blue}
    \tikzColoredSquarehl{2}{-7}{2}{white!45!blue}
    \tikzColoredSquarehl{2}{-8}{3}{white!45!blue}
    \tikzColoredSquarehl{2}{-9}{4}{white!45!blue}
    \tikzColoredSquarehl{1}{-9}{4}{white!45!blue}
    \tikzColoredSquare{0}{-7}{0}{white!45!orange}
    \tikzColoredSquare{0}{-8}{1}{white!45!orange}
    \tikzColoredSquare{0}{-9}{2}{white!45!orange}
    \tikzColoredSquare{0}{-10}{3}{white!45!orange}
    \tikzColoredSquare{1}{-10}{3}{white!45!orange}
    \tikzColoredSquare{1}{-11}{4}{white!45!orange}
    \tikzColoredSquare{2}{-10}{0}{white!75!teal}
    \tikzColoredSquare{2}{-11}{1}{white!75!teal}
    \tikzColoredSquare{2}{0}{2}{white!75!teal}
    \tikzColoredSquare{1}{0}{2}{white!75!teal}
    \tikzColoredSquare{1}{-1}{3}{white!75!teal}
    \tikzColoredSquare{1}{-2}{4}{white!75!teal}
    \tikzColoredSquare{0}{-11}{0}{white!45!yellow}
\end{tikzpicture}
\hfill
    \caption{The same orbit board twice, the one on the left highlights a single whorm, the one on the right highlights all whorms.}\label{fig:onewhorm}
\end{figure}

Each whorm in an orbit board of $\VV\times [k]$ starts on the left, or the right, or both left and
right; we call the former \emph{one-tailed} (specifically \emph{right-tailed} or
\emph{left-tailed}) and the latter \emph{two-tailed}. We will also use the terms
\emph{left-whorm} and \emph{right-whorm} for short in the one-tailed case.

Since these whorms move down the orbit board at every step, except for one move to the center, we consider
them as  a sequence of function values in
the orbit board which start at $0$ and end at $k$, where one value is repeated when moving
into the center.  We call these \emph{center-seeking whorms}.
(Since an orbit board is actually a cylinder, we have a ``can of whorms'' to deal with.)

In the left orbit of  Figure~\ref{fig:onewhorm} we isolate one
example of a left whorm: 
$$\xi = \{(\ell,(0,3,3)),(\ell,(1,4,0)),(\ell,(2,2,1)),(c,(2,2,1)), (c,(0,3,2)), (c,(1,4,3))\},$$
visualized within an orbit board of $\mathcal{F}_4(\VV)$.  

It is easy to see that an orbit board is tiled either entirely by one-tailed whorms or
entirely by two-tailed whorms. In fact, if $f(\ell)=f(r)$ in any row, then $w^df(\ell) =
w^df(r)$ for all $d\in \mathbb N$ and all $f\in\mathcal F_k(\VV)$. This is because the
labels of $\ell$ and $r$ represent the result of whirling at \emph{incomporable} elements of
the poset $\VV$. Furthermore, this forces the whorms to go to the middle from both sides in
the same row, which makes the elements of the two-tails whirl-connected, so happens exactly
when we have a two-tailed whorms. In the situation where $f(\ell ) \neq f(r)$ in any row,
this also persists, so none of the whorms can be two-tailed. 

We first observe that all one-tailed whorms have $k+2$ elements, since each contains the $k+1$
elements $0,...,k$, exactly one of which is doubled. 
Define $\tl(\xi): = 1+\min \{f(c): (c,f)\in \xi\}$, the number of elements in an
outer column (the ``tail length'') and $\hl(\xi): = k+2 - \tl(\xi)$, the number
of elements in the center column (the ``head length''). 
For the red whorm in the orbit on the left of Figure~\ref{fig:onewhorm}, $\tl(\xi) = 3$ and $\hl(\xi) = 3$. 
For the blue two-tailed worm on the right of Figure~\ref{fig:whormboard}, $\tl(\xi) =
4$ and $\hl(\xi) = 2$. Note that in a one-tailed orbit board, each whorm has the same
cardinality, namely $k+2$; but in a two-tailed orbit board, the whorms can have different
cardinalities (e.g., 10 and 8 in Figure~\ref{fig:whormboard}).  

\begin{ex} \label{ex:bp}
The right orbit board in Figure~\ref{fig:onewhorm} is the previous example with all the whorms
colored. The number of elements in the left column of the  yellow, red, and
orange whorms are $5,3,4$ respectively,  and the orbit board is of length $12$. 
\end{ex}

It follows that the order of whirling divides the sum of  $\tl(\xi)$ over all whorms
$\xi\in S$. In the setting of $\mathcal{F}_k(\VV)$, as long as we know
$\tl(\xi)$ and whether $f(\ell) = 0$, $f(r) =0$, or both, then we can recover the
entire whorm. 

The center of any row intersects with eactly one whorm. A natural cyclic ordering on whorms in an orbit board arises from scanning down its center fiber.
\begin{definition}
We will place a circular order on the whorms. Let $\xi_1$ and $\xi_2$ be whorms
in an orbit board of $\mathcal{F}_k(\VV)$. If there exists $(c, f) \in \xi_1$ with
$f(c) = k$ such that $(c, w(f))\in \xi_2$, then we say $\xi_2$ is \emph{in front of}
$\xi_1$. 
We call a sequence of whorms \emph{consecutive} if each is in front of the next.  
\end{definition} 
\begin{ex}
In Figure~\ref{fig:onewhorm} the blue (horizontal lines) whorm is in front of the red
(crosshatch) whorm,
which is in front of the green (northwest lines) whorm.  
\end{ex}

It is not hard to show that a within a one-tailed orbit board, consecutive whorms alternate starting from the left
and from the right. For when a whorm goes into the center from one side, say the left
without loss of generality, then the next value in the left column will be 0, while in the
right column it will be greater than 0. This persists until the next whorm moves to the
center (possibly during the same step).  At that point, by definition of whirling, the
value in the center column will change from $k$ to $1+\max\{f(\ell ), f(r) \} = 1+f(r)$, meaning that
the values in the right column and in the center will be equal.  See the one-tailed examples
in Figures~\ref{fig:whormboard}~and~\ref{fig:onewhorm}.  We use this in the next lemma,
which is key to proving our main results. 

\begin{lem}\label{lem:bp_consecutive}
Given an orbit board $\mathcal{R}$ of $w$ on $\mathcal{F}_k(\VV)$, let
$\xi_1,\xi_2,$ and $\xi_3$ be three consecutive whorms (not necessarily
all distinct), that 
is, $\xi_3$ is in front of $\xi_2$ which is in front of $\xi_1$ in
$\mathcal{R}$. 
\begin{enumerate}

\item If $\mathcal{R}$ is tiled entirely by one-tailed whorms, then 
$$\tl(\xi_1)+\tl(\xi_2)+\tl(\xi_3)=2(k+2).$$

\item Otherwise, if $\mathcal{R}$ is tiled entirely by two-tailed whorms, then 
$$\tl(\xi_1)+\tl(\xi_2) = k+2.$$

\end{enumerate}
\end{lem}

\begin{proof} 
(1) For the case of one-tailed whorms the key observation is that $\tl (\xi_{3}) =
\hl(\xi_1) + \hl(\xi_2)$ for consecutive whorms.  This holds because the tail of
$\xi_{3}$ begins one row below where the head of $\xi_{1}$ starts, and ends one
row below where $\xi_{2}$ ends.  But since $\tl (\xi ) = k+2 - \hl (\xi )$
for any $\xi $ we get
\[
\tl(\xi_3) = \hl(\xi_1) + \hl(\xi_2) = (k+2-\tl (\xi_{1}))+ (k+2-\tl (\xi_{2})),
\]
which implies the result. See Figure~\ref{fig:whormboard} for examples.  

(2) For the case of two-tailed whorms a similar argument shows that $\tl (\xi_{2})
=\hl(\xi_1)$ for two consecutive whorms. The result then follows as before.  
\end{proof}

\noindent In fact, the entire orbit board can be reconstructed simply from knowing the
values of $\tl(\xi_1)$ and $\tl(\xi_2)$ for two consecutive
whorms in the one-tailed case, and from a single $\tl(\xi_1)$ in the two-tailed case.  

\begin{ex}
In Figure~\ref{fig:onewhorm} we have $k=4$, $\tl(\text{green})=4$, $\tl(\text{red})=3 $, and $\tl(\text{blue})=5$, which sum to $12=2(4+2)$.  
\end{ex}

\begin{lem}\label{lem:srotate}
Given an orbit board with one-tailed whorms, let $\xi_1$, $\xi_2$,
$\xi_3$, and $\xi_4$ be consecutive, then $$\tl(\xi_4) = \tl(\xi_1).$$
Otherwise, if the orbit board contains two-tailed whorms, then $\tl(\xi_1) = \tl(\xi_3)$.
\end{lem}

\begin{proof} First assume the orbit board is one-tailed.
By Lemma~\ref{lem:bp_consecutive}, the difference of the equations,
\[\tl(\xi_1)+\tl(\xi_2)+\tl(\xi_3)=2(k+2),\text{ and,}\]
\[\tl(\xi_2)+\tl(\xi_3)+\tl(\xi_4)=2(k+2),\]
give the result. The proof follows similarly when the orbit board has two-tailed whorms. 
\end{proof}

\begin{lem}\label{lem:maxWhorms}
Let $\mathcal{R}$ be an orbit board  of $w$ on $\mathcal{F}_k(\VV)$.  
\begin{enumerate}

\item If $\mathcal{R}$ is tiled entirely by one-tailed whorms, then 
there are at most six distinct whorms.  

\item Otherwise, if $\mathcal{R}$ is tiled entirely by two-tailed whorms, then 
there are at most two distinct whorms.  
\end{enumerate}
\end{lem}
\begin{proof}
For orbit boards with one-tailed whorms, we are \emph{not} claiming the board starts to
repeat after three steps, since  $\xi_4$ will start on the opposite side from 
$\xi_1$. If we keep applying the previous Lemma to even more consecutive whorms, we
see $\tl(\xi_5)=\tl(\xi_2)$ and $\tl(\xi_6) = \tl(\xi_3)$. Finally
we get $\tl(\xi_k)=\tl(\xi_{k-6})$ for $k\geq 7$.  Since a whorm is completely
characterized by its tail length and sidedness, this means that the orbit board must repeat
by the row in which $\xi_{7}$ begins.  The proof in the two-tailed case is similar. 
\end{proof}

\begin{figure}
\hfill
\begin{tikzpicture}[baseline={([yshift=-.8ex]current bounding box.center)}]

    \def\m{0}
    \tikzColoredSquare{0}{\m}{0}{white}
    \tikzColoredSquare{1}{\m}{1}{white}
    \tikzColoredSquare{2}{\m}{1}{white}
    
    \def\m{-1}
    \tikzColoredSquare{0}{\m}{1}{white}
    \tikzColoredSquare{1}{\m}{2}{white}
    \tikzColoredSquare{2}{\m}{0}{white}
    
    \def\m{-2}
    \tikzColoredSquare{0}{\m}{2}{white}
    \tikzColoredSquare{1}{\m}{2}{white}
    \tikzColoredSquare{2}{\m}{1}{white}
    
    \def\m{-3}
    \tikzColoredSquare{0}{\m}{0}{white}
    \tikzColoredSquare{1}{\m}{2}{white}
    \tikzColoredSquare{2}{\m}{2}{white}
    
    \def\m{-4}
    \tikzColoredSquare{0}{\m}{1}{white}
    \tikzColoredSquare{1}{\m}{1}{white}
    \tikzColoredSquare{2}{\m}{0}{white}
    
    \def\m{-5}
    \tikzColoredSquare{0}{\m}{0}{white}
    \tikzColoredSquare{1}{\m}{2}{white}
    \tikzColoredSquare{2}{\m}{1}{white}
    
    \def\m{-6}
    \tikzColoredSquare{0}{\m}{1}{white}
    \tikzColoredSquare{1}{\m}{2}{white}
    \tikzColoredSquare{2}{\m}{2}{white}

    \def\m{-7}
    \tikzColoredSquare{0}{\m}{2}{white}
    \tikzColoredSquare{1}{\m}{2}{white}
    \tikzColoredSquare{2}{\m}{0}{white}
\end{tikzpicture}
\hfill
\begin{tikzpicture}[baseline={([yshift=-.8ex]current bounding box.center)}]
    \def\m{0}
    \tikzColoredSquare{0}{\m}{0}{white}
    \tikzColoredSquare{1}{\m}{2}{white}
    \tikzColoredSquare{2}{\m}{2}{white}
    
    \def\m{-1}
    \tikzColoredSquare{0}{\m}{1}{white}
    \tikzColoredSquare{1}{\m}{3}{white}
    \tikzColoredSquare{2}{\m}{0}{white}
    
    \def\m{-2}
    \tikzColoredSquare{0}{\m}{2}{white}
    \tikzColoredSquare{1}{\m}{2}{white}
    \tikzColoredSquare{2}{\m}{1}{white}
    
    \def\m{-3}
    \tikzColoredSquare{0}{\m}{0}{white}
    \tikzColoredSquare{1}{\m}{3}{white}
    \tikzColoredSquare{2}{\m}{2}{white}
    
    \def\m{-4}
    \tikzColoredSquare{0}{\m}{1}{white}
    \tikzColoredSquare{1}{\m}{3}{white}
    \tikzColoredSquare{2}{\m}{3}{white}
    
    \def\m{-5}
    \tikzColoredSquare{0}{\m}{2}{white}
    \tikzColoredSquare{1}{\m}{2}{white}
    \tikzColoredSquare{2}{\m}{0}{white}
    
    \def\m{-6}
    \tikzColoredSquare{0}{\m}{0}{white}
    \tikzColoredSquare{1}{\m}{3}{white}
    \tikzColoredSquare{2}{\m}{1}{white}

    \def\m{-7}
    \tikzColoredSquare{0}{\m}{1}{white}
    \tikzColoredSquare{1}{\m}{2}{white}
    \tikzColoredSquare{2}{\m}{2}{white}

    \def\m{-8}
    \tikzColoredSquare{0}{\m}{2}{white}
    \tikzColoredSquare{1}{\m}{3}{white}
    \tikzColoredSquare{2}{\m}{0}{white}

    \def\m{-9}
    \tikzColoredSquare{0}{\m}{3}{white}
    \tikzColoredSquare{1}{\m}{3}{white}
    \tikzColoredSquare{2}{\m}{1}{white}
\end{tikzpicture}
\hfill
\caption{Two orbits of whirling, the left one of length 8 for $k=2$ and the right one of length
10 for $k=3$, each, confirming the whirling order in these cases.}\label{fig:smallk}
\end{figure}

\begin{figure}
\hfill
\begin{tikzpicture}[baseline={([yshift=-.8ex]current bounding box.center)}]

    \def\m{0}
    \tikzColoredSquare{0}{\m}{0}{white}
    \tikzColoredSquare{1}{\m}{0}{white}
    \tikzColoredSquare{2}{\m}{0}{white}
    
    \def\m{-1}
    \tikzColoredSquare{0}{\m}{0}{white}
    \tikzColoredSquare{1}{\m}{1}{white}
    \tikzColoredSquare{2}{\m}{0}{white}
    
    \def\m{-2}
    \tikzColoredSquare{0}{\m}{1}{white}
    \tikzColoredSquare{1}{\m}{1}{white}
    \tikzColoredSquare{2}{\m}{1}{white}
\end{tikzpicture}
\hfill
\begin{tikzpicture}[baseline={([yshift=-.8ex]current bounding box.center)}]
    \def\m{0}
    \tikzColoredSquare{0}{\m}{0}{white}
    \tikzColoredSquare{1}{\m}{1}{white}
    \tikzColoredSquare{2}{\m}{1}{white}
    
    \def\m{-1}
    \tikzColoredSquare{0}{\m}{1}{white}
    \tikzColoredSquare{1}{\m}{1}{white}
    \tikzColoredSquare{2}{\m}{0}{white}
\end{tikzpicture}
\hfill
\caption{Both orbits of whirling when $k=1$, confirming the order is 6 in this case.}\label{fig:orderidealv}
\end{figure}
\begin{thm}\label{thm:halfwaysymm}
Let $(x,y,z)\in \mathcal{F}_k(\VV)$, then $w^{k+2}(x,y,z) = (z,y,x)$.
\end{thm}

\begin{proof}
Let $\mathcal{R}$ be an orbit board of $w$. The proof involves a careful analysis of when
whorms move into the center of the orbit board.  For the one-tailed case, we use
Lemma~\ref{lem:bp_consecutive} and that $\tl (\xi ) = k+2 - \hl (\xi )$ for any
$\xi $ to conclude that
$$\hl(\xi_1)+\hl(\xi_2)+\hl(\xi_3)=(k+2).$$
At this point, $\xi_4$ arrives from the side opposite $\xi_{1}$, but with the
same tail length.  Therefore these two whorms must be reflections of one another (around the vertical middle of the
$\mathcal{R}$, though translated vertically as well).  This forces the values in row
$k+2+1$ to be the mirror reflection of those in row 1 (where $\xi_{1}$ starts). 

In the two-tailed case, that is, when $x=z$, by Lemma~\ref{lem:srotate} we know $\tl(\xi_1) = \tl(\xi_3)$.  
This implies $\xi_1$ and $\xi_3$ are the same whorm, since they are two-tailed whorms. Therefore, the orbit board
repeats after $\tl(\xi_1)+\tl(\xi_2) = k+2$ steps, by Lemma~\ref{lem:bp_consecutive}.  
\end{proof}

\begin{cor}\label{cor:veeorder}
The order of $w$ on $\mathcal{F}_k(\VV)$ divides $2(k+2)$.
\end{cor}

We are now ready to finish off the proofs of our  main results for this section. 

\medskip 
\begin{proof}[Proof of Theorem~\ref{thm:Vorder}]
\noindent Easy computations shows that order of $w$
is \textit{exactly} $2(k+2)$ for $k=1, 2$, or 3. (See
Figures~\ref{fig:smallk}--\ref{fig:orderidealv}.)   Now for $k\geq 4$ we can pick any three
distinct values  
$b_{1}, b_{2}, b_{3}$ that sum to $k+2$ and construct an orbit board from six consecutive whorms with 
three distinct tail lengths and $\tl (\xi_{i})=b_{i} = \tl (\xi_{i+3})$ for
$i\in [3]$.  These whorms will all be distinct, so we get no repeats until after $2(k+2)$
steps.  It is also straightforward to directly construct an order ideal in $\VV \times [k]$
whose orbit is of this size.  Hence, The order of $w$ on $\mathcal{F}_k(\VV)$ is
\emph{exactly} $2(k+2)$. 
\end{proof}

In our proofs the following definition will be handy.  
\begin{definition}\label{def:superorbit}
A \emph{super orbit} of an action is an orbit repeated one or more times.  The length of a
super orbit is always a multiple of the orbit size.   Similarly a \emph{super orbit board}
will be the array representing a super orbit. 
\end{definition}

\begin{proof}[Proof of Theorem~\ref{thm:Vhom}]
The symmetry homomesy (1) is immediate from Theorem~\ref{thm:halfwaysymm}.  

For (2) if $I$ is
an order ideal of $\VV\times [k]$, then there is a $P$-partition $f=(x,y,z)\in\mathcal{F}$
in bijection with $I$ given by Theorem~\ref{thm:eqwhirl}. We know that 
\begin{itemize}

\item $\ell_1\in I \iff x>0$, 

\item $r_1\in I \iff z>0$, and 

\item  $c_k\in I \iff y=k$. 

\end{itemize}

In the one-tailed case, let $\xi_1,\xi_2,\dots,\xi_6$ be six
consecutive whorms partitioning the (possibly super) orbit board by Lemma~\ref{lem:maxWhorms}.  This gives a
total of $4(k+2)$ entries in the left and right columns, six of which are $0$,
leaving $4(k+2)-6$ nonzero entries, each of which contributes to either $\chi_{r_1}$ or
$\chi_{\ell_1}$. The number of entries equal to $k$ in the center column is 6 because there
are 6 whorms that end at $k$ in the center. Thus, we get the average
$$\frac{\sum_{f\in\mathcal{R}}\chi_{\ell_1}(f) + \chi_{r_1}(f) - \chi_{c_k}(f)}{\# \mathcal{R}}=\frac{4(k+2)-12}{2(k+2)} = \frac{2(k-1)}{k+2}.$$
The two-tailed case follows analogously.  
\end{proof} 
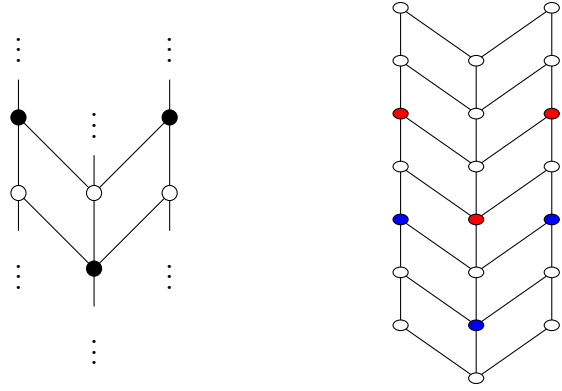
\begin{figure}
\begin{tikzpicture}[baseline={([yshift=-.8ex]current bounding box.center)}]
    \draw[-] (-1,1.5) -- (-1,3.5);
    \draw[-] (1,1.5) -- (1,3.5);
    \draw[-] (0,0.5) -- (0,2.5);

    \node at (-1,1) {$\vdots$};
    \node at (0,0) {$\vdots$};
    \node at (1,1) {$\vdots$};
    
    \draw[-] (-1,2) -- (0,1);
    \draw[-] (1,2) -- (0,1);
    
    \draw [fill=white] (-1,2) circle [radius =.1];
    \draw [fill] (0,1) circle [radius =.1];
    \draw [fill=white] (1,2) circle [radius =.1];
    
    \draw[-] (-1,3) -- (0,2);
    \draw[-] (1,3) -- (0,2);
    
    \draw [fill] (-1,3) circle [radius =.1];
    \draw [fill=white] (0,2) circle [radius =.1];
    \draw [fill] (1,3) circle [radius =.1];
    
    \node at (-1,4) {$\vdots$};
    \node at (0,3) {$\vdots$};
    \node at (1,4)  {$\vdots$};
\end{tikzpicture}\ \ \ \ \ \ \ \ \ \ \ \ \ \ \ \ \ \ \ \ 
\begin{tikzpicture}[yscale = 0.7, baseline={([yshift=-.8ex]current bounding box.center)}]
    \draw[-] (-1,1) -- (-1,7);
    \draw[-] (1,1) -- (1,7);
    \draw[-] (0,0) -- (0,6);

    \draw[-] (-1,1) -- (0,0);
    \draw[-] (1,1) -- (0,0);
    
    \draw [fill=white] (-1,1) circle [radius =.1];
    \draw [fill=white] (0,0) circle [radius =.1];
    \draw [fill=white] (1,1) circle [radius =.1];
    
    \draw[-] (-1,2) -- (0,1);
    \draw[-] (1,2) -- (0,1);
    
    \draw [fill=white] (-1,2) circle [radius =.1];
    \draw [fill=blue] (0,1) circle [radius =.1];
    \draw [fill=white] (1,2) circle [radius =.1];
    
    \draw[-] (-1,3) -- (0,2);
    \draw[-] (1,3) -- (0,2);
    
    \draw [fill=blue] (-1,3) circle [radius =.1];
    \draw [fill=white] (0,2) circle [radius =.1];
    \draw [fill=blue] (1,3) circle [radius =.1];
    
    \draw[-] (-1,4) -- (0,3);
    \draw[-] (1,4) -- (0,3);
    
    \draw [fill=white] (-1,4) circle [radius =.1];
    \draw [fill=red] (0,3) circle [radius =.1];
    \draw [fill=white] (1,4) circle [radius =.1];
    
    \draw[-] (-1,5) -- (0,4);
    \draw[-] (1,5) -- (0,4);
    
    \draw [fill=red] (-1,5) circle [radius =.1];
    \draw [fill=white] (0,4) circle [radius =.1];
    \draw [fill=red] (1,5) circle [radius =.1];
    
    \draw[-] (-1,6) -- (0,5);
    \draw[-] (1,6) -- (0,5);
    
    \draw [fill=white] (-1,6) circle [radius =.1];
    \draw [fill=white] (0,5) circle [radius =.1];
    \draw [fill=white] (1,6) circle [radius =.1];
    
    \draw[-] (-1,7) -- (0,6);
    \draw[-] (1,7) -- (0,6);
    
    \draw [fill=white] (-1,7) circle [radius =.1];
    \draw [fill=white] (0,6) circle [radius =.1];
    \draw [fill=white] (1,7) circle [radius =.1];
\end{tikzpicture}
\caption{On the left, a section of $\VV_k$ with elements of a flux capacitor configuration, $\ell_i,r_i$ and $c_{i-1}$, shaded black. On the right, a section of  $\VV_7$ with elements counted by $F_3-F_5$, shaded black.}\label{fig:flux}
\end{figure}

Using the whorm decomposition of the orbit board also gives a very clean proof of the
following ``flux capacitor'' homomesy for rowmotion on $\mathcal{J}(\VV_k)$.   
Let $F_i = \chi_{\ell_i}+\chi_{r_i} + \chi_{c_{i-1}}$, as shown in Figure~\ref{fig:flux}.

\begin{lem}\label{lem:fluxpairs}
Let $\mathcal{R}$ be a \textbf{one-tailed} orbit board for the action of rowmotion on order ideals of
$\mathcal{J}(\VV_k)$.  The the average value of the difference of successive flux-capacitor
indicator functions, $F_{i} - F_{i+1}$, is $\frac{3}{k+2}$ for $i\in [2,k-1]$.
\end{lem}

\begin{proof}
From the whirling orbit-board point of view, the statistic $F_{i} - F_{i+1}$ is nonzero
exactly when $i$ appears in the left or right column or $i-1$
appears in the center column.   (Otherwise, e.g., both left-column indicator functions are zero,
or both are one, canceling out when we take the difference $F_{i} - F_{i+1}$.) 

In the one-tailed case, if
$i$ appears in the tail of a given whorm then $i-1$ will not appear in the head, and
vice versa. So each of the six whorms in a super orbit board (Lemma~\ref{lem:maxWhorms})
will contribute exactly 1 to the statistic $F_{i} - F_{i+1}$, depending on whether $i$
appears in a side column or $i-1$ appears in the center column.  This gives a total
contribution of 6 to the statistic.  Divide this by the size of the super orbit, $2(k+2)$,
to get the result.  
\end{proof}

Unfortunately in the two-tailed (symmetric) case, this lemma fails to hold, as the following
counterexample shows.  Earlier versions of this work incorrectly asserted that the above
lemma held in \textbf{all} cases, which would have led to a stronger version of
Theorem~\ref{thm:flux} which omitted the ``symmetrically placed'' hypothesis.  

\begin{figure}
\hfill
\begin{tikzpicture}[baseline={([yshift=-.8ex]current bounding box.center)}]

    \def\m{0}
    \tikzColoredSquare{0}{\m}{6}{white!45!red}
    \tikzColoredSquare{1}{\m}{6}{white!45!red}
    \tikzColoredSquare{2}{\m}{6}{white!45!red}

    \def\m{-1}
    \tikzColoredSquare{0}{\m}{0}{white!45!blue}
    \tikzColoredSquare{1}{\m}{7}{white!45!red}
    \tikzColoredSquare{2}{\m}{0}{white!45!blue}
    
    \def\m{-2}
    \tikzColoredSquare{0}{\m}{1}{white!45!blue}
    \tikzColoredSquare{1}{\m}{8}{white!45!red}
    \tikzColoredSquare{2}{\m}{1}{white!45!blue}
    
    \def\m{-3}
    \tikzColoredSquare{0}{\m}{2}{white!45!blue}
    \tikzColoredSquare{1}{\m}{9}{white!45!red}
    \tikzColoredSquare{2}{\m}{2}{white!45!blue}
    
    \def\m{-4}
    \tikzColoredSquare{0}{\m}{3}{white!45!blue}
    \tikzColoredSquare{1}{\m}{3}{white!45!blue}
    \tikzColoredSquare{2}{\m}{3}{white!45!blue}

    \def\m{-5}
    \tikzColoredSquare{0}{\m}{0}{white!45!red}
    \tikzColoredSquare{1}{\m}{4}{white!45!blue}
    \tikzColoredSquare{2}{\m}{0}{white!45!red}
    
    \def\m{-6}
    \tikzColoredSquare{0}{\m}{1}{white!45!red}
    \tikzColoredSquare{1}{\m}{5}{white!45!blue}
    \tikzColoredSquare{2}{\m}{1}{white!45!red}
    
    \def\m{-7}
    \tikzColoredSquare{0}{\m}{2}{white!45!red}
    \tikzColoredSquare{1}{\m}{6}{white!45!blue}
    \tikzColoredSquare{2}{\m}{2}{white!45!red}
    
    \def\m{-8}
    \tikzColoredSquare{0}{\m}{3}{white!45!red}
    \tikzColoredSquare{1}{\m}{7}{white!45!blue}
    \tikzColoredSquare{2}{\m}{3}{white!45!red}
    
    \def\m{-9}
    \tikzColoredSquare{0}{\m}{4}{white!45!red}
    \tikzColoredSquare{1}{\m}{8}{white!45!blue}
    \tikzColoredSquare{2}{\m}{4}{white!45!red}
    
    \def\m{-10}
    \tikzColoredSquare{0}{\m}{5}{white!45!red}
    \tikzColoredSquare{1}{\m}{9}{white!45!blue}
    \tikzColoredSquare{2}{\m}{5}{white!45!red}
\end{tikzpicture}
\hfill
\caption{An orbit board of whirling in $\calF_{9}(\VV)$, tiled by 2 two-tailed worms}\label{fig:2tailCEG}

\end{figure}

\begin{ex}\label{ex:2tailedCEG}
Consider the orbit board from whirling $\calF_{9}(\VV )$ starting at $666$, as shown in
Figure~\ref{fig:2tailCEG}. 
we see that
$F_{2}-F_{3}$ has average $\frac{4}{k+2}$ (since there are four 2s in outer column and no 1s in
the center column).  Meanwhile, $F_{7}-F_{8}$ has average $\frac{2}{k+2}$ (since there are zero
7s in outer column and two 6s in the middle).  
\end{ex}

\begin{thm}\label{thm:flux}
For $k>1$. Let $F_i = \chi_{\ell_i}+\chi_{r_i} + \chi_{c_{i-1}}$ be the indicator function
of the flux capacitor at height $i$. Under the action of rowmotion on order ideals of
$\mathcal{J}(\VV_k)$, the difference of \textbf{symmetrically-placed} flux-capacitor
statistics, $F_{i} - F_{k+2-i}$, is $\frac{3(k+2-2i)}{k+2}$-mesic.
\end{thm}

\begin{proof}
For the one-tailed case, we can use a simple telescoping sum to get a more general result.  
For any $i<j$ we may write
\[
F_{i}-F_{j} = (F_{i}-F_{i+1})+(F_{i+1}-F_{i+2})+\dots + (F_{j-1}-F_{j}). 
\]
Grouping consecutive pairs and using Lemma~\ref{lem:fluxpairs}, we get the average of
$F_i-F_j$ is $\displaystyle \frac{3(j-i)}{k+2}$.  This specializes to the stated homomesy when
$j=k+2-i$ for symmetrically-placed flux-capacitor statistics.  A similar argument works for
$j>i$, when the homomesic average will be non-positive.    

For the two-tailed case we need a separate argument, since Lemma~\ref{lem:fluxpairs} no
longer applies.  Fix $j=k+2-i$, and assume again for now that $i<j$.  
Within any two-tailed whorm, the multiset of labels is the set of all integers in
$[0,k]$ with integers in $[0,\tl (\xi)-2]$ appearing twice, and $\tl (\xi)-1$
appearing thrice (where the whorm moves to the center).  Call the two whorms tiling the
length $k+2$ (super) orbit board $\xi_{1}$ and $\xi_{2}$.  By
Lemma~\ref{lem:bp_consecutive} we have  
\[
(\tl(\xi_{1})-1) + (\tl( \xi_{2})-1) =k. 
\]
For example, in the length 11 orbit board in Figure~\ref{fig:2tailCEG}, the labels on the red
and blue whorms are 6 and 3, which sum to $k=9$.  

Now assume that $\tl(\xi_{1})-1 =m$, so $\tl( \xi_{2})-1 = k-m$.  From the
whirling whorms point of view, the contribution of a given label $\lambda$ in the tiling to
the statistic $F_{i} - F_{k+2-i}$ is given as follows:
\begin{itemize}

\item If $\lambda$ lies in a side column and $i\leq \lambda \leq k+1-i$, then that label
contributes $+1$ to the statistic.  

\item If $\lambda$ lies in the center column and $i-1\leq \lambda \leq k-i$, then that label
contributes $+1$ to the statistic.  

\end{itemize}
Outside of the stated intervals, the statistics cancel each other out, so $\lambda$
contributes zero. This all follows directly from the definition of $F_{i}$.

\begin{ex}
Here are the two whorms from the orbit shown in Figure~\ref{fig:2tailCEG} written
horizontally in one line. The tails are highlighted in blue and the heads highlighted in
green. The labels that contribute to $F_2-F_9$ are circled. Here $m=6$ and $k-m=3$, both of
which lie between $i=2$ and $k+2-i = 9$.  
\[
\begin{tikzpicture}[baseline={([yshift=-.8ex]current bounding box.center)}]

    \def\m{0}
    \tikzColoredSquare{-1}{\m}{\xi_1}{white}
    \tikzColoredSquare{0}{\m}{0}{white!65!blue}
    \tikzColoredSquare{1}{\m}{1}{white!65!blue}
    \tikzColoredSquare{2}{\m}{\raisebox{.5pt}{\textcircled{\raisebox{-.9pt} {2}}}}{white!65!blue}
    \tikzColoredSquare{3}{\m}{\raisebox{.5pt}{\textcircled{\raisebox{-.9pt} {3}}}}{white!65!blue}
    \tikzColoredSquare{4}{\m}{\raisebox{.5pt}{\textcircled{\raisebox{-.9pt} {4}}}}{white!65!blue}
    \tikzColoredSquare{5}{\m}{\raisebox{.5pt}{\textcircled{\raisebox{-.9pt} {5}}}}{white!65!blue}
    \tikzColoredSquare{6}{\m}{\raisebox{.5pt}{\textcircled{\raisebox{-.9pt} {6}}}}{white!65!blue}
    \tikzColoredSquare{7}{\m}{\raisebox{.5pt}{\textcircled{\raisebox{-.9pt} {6}}}}{white!65!green}
    \tikzColoredSquare{8}{\m}{\raisebox{.5pt}{\textcircled{\raisebox{-.9pt} {7}}}}{white!65!green}
    \tikzColoredSquare{9}{\m}{8}{white!65!green}
    \tikzColoredSquare{10}{\m}{9}{white!65!green}

    \def\m{-2}
    \tikzColoredSquare{-1}{\m}{\xi_2}{white}
    \tikzColoredSquare{0}{\m}{0}{white!65!blue}
    \tikzColoredSquare{1}{\m}{1}{white!65!blue}
    \tikzColoredSquare{2}{\m}{\raisebox{.5pt}{\textcircled{\raisebox{-.9pt} {2}}}}{white!65!blue}
    \tikzColoredSquare{3}{\m}{\raisebox{.5pt}{\textcircled{\raisebox{-.9pt} {3}}}}{white!65!blue}
    \tikzColoredSquare{4}{\m}{\raisebox{.5pt}{\textcircled{\raisebox{-.9pt} {3}}}}{white!65!green}
    \tikzColoredSquare{5}{\m}{\raisebox{.5pt}{\textcircled{\raisebox{-.9pt} {4}}}}{white!65!green}
    \tikzColoredSquare{6}{\m}{\raisebox{.5pt}{\textcircled{\raisebox{-.9pt} {5}}}}{white!65!green}
    \tikzColoredSquare{7}{\m}{\raisebox{.5pt}{\textcircled{\raisebox{-.9pt} {6}}}}{white!65!green}
    \tikzColoredSquare{8}{\m}{\raisebox{.5pt}{\textcircled{\raisebox{-.9pt} {7}}}}{white!65!green}
    \tikzColoredSquare{9}{\m}{8}{white!65!green}
    \tikzColoredSquare{10}{\m}{9}{white!65!green}
    
\end{tikzpicture}
\]
The number of circled labels in the whorms will give us the total value of the statistic
$F_2-F_9$ provided the circled labels in the tail section (blue part) of the whorm are
counted twice since the whorm is two-tailed. Divide by the length of the orbit board, $11$
to get $\frac{21}{11} = \frac{3(k+2-2i)}{k+2}$. 
\end{ex}

In general for $i\leq m\leq k+2-i$ the ``symmetrically-placed'' condition forces the number of contributing labels
in the \emph{head} of $\xi_{1}$ to equal the number of contributing labels in the \emph{tail} of
$\xi_{2}$.  Hence, the total number of $+1$'s in the orbit will be thrice the
number of contributing labels in the tail of $\xi_{1}$ plus thrice the number of
contributing labels in the tail of $\xi_{2}$.  But this is just thrice the total number
of contributing labels, since the number in the tail of $\xi_{2}$ equals the number in
the head of $\xi_{1}$.  Hence, we get a total contribution of $3((k+2-i) -i)
=3(k+2-2i)$.  The homomesy follows by dividing by the length, $k+2$, of the orbit board. A
symmetrical argument, where we count contributing $-1$'s, handles the case where
$i>j=k+2-i$.  
\end{proof}

\section{Dynamics for rowmotion on $\BB_n \times [k]$} \label{sec:claw}
In this section we extend our results for the chain of V's poset to a ``chain of claws''
poset, defined below. The equivariant bijection (Theorem~\ref{thm:eqwhirl}) and techniques
from Section~\ref{sec:vee} extend with only limited difficulty to this new setting.  

\subsection{The chain of claws poset}\label{ss:chainclaws}

\begin{definition}\label{def:chainclaws}
We define the \emph{claw poset} $\BB_n = \{ b_1,\dots,b_n,\widehat{0}\}$ where each $b_{i}$ covers $\widehat{0}$.  For example, the Hasse diagram of $\BB_4$ would be \begin{tikzpicture}[scale=0.6,baseline={([yshift=-.8ex]current bounding box.center)}]
    \draw[-] (-1,1) -- (1.66,0);
    \draw[-] (-.33,1) -- (1.66,0);
    \draw[-] (.33,1) -- (1.66,0);
    \draw[-] (1,1) -- (1.66,0);
    
    \draw [fill] (1.66,0) circle [radius =.1];
    \draw [fill] (-1,1) circle [radius =.1];
    \draw [fill] (-.33,1) circle [radius =.1];
    \draw [fill] (1,1) circle [radius =.1];
    \draw [fill] (.33,1) circle [radius =.1];
\end{tikzpicture}.
The \textit{chain of claws poset} is defined to be $\BB_{n}\times [k]$.  
\end{definition}

Using the established equivariant bijection (Theorem~\ref{thm:eqwhirl}) between
$\mathcal{J}(\BB_n\times [k])$ and $k$-bounded $P$-partitions $\mathcal{F}_k(\BB_n)$ that
sends rowmotion to whirling, we can prove homomesy and periodicity results similar to those
for $\BB_2=\VV$.  
Now instead of \emph{triples} of numbers, we will consider orbit boards of
\emph{$(n+1)$-tuples} on $[0,k]$,  $\left(f(b_{1}),f(b_{2}),\dots
,f(b_{n}),f(\widehat{0})\right)$, satisfying $f(b_{i})\leq 
f(\widehat{0})$ for each $i\in [n]$.  Notice that for the case $n=2$, we have moved the column
formerly in the middle (corresponding to $\widehat{0}$) to be the \textbf{rightmost} one.  
\begin{figure}
    
\begin{tikzpicture}[baseline={([yshift=-.8ex]current bounding box.center)}]

    \tikzColoredSquare{-3.25}{1}{\text{columns:} }{{rgb,255:red,255; green,255; blue,255}}
    \tikzColoredSquare{-1}{1}{b_1 }{{rgb,255:red,255; green,255; blue,255}}
    \tikzColoredSquare{0}{1}{b_2}{{rgb,255:red,255; green,255; blue,255}}
    \tikzColoredSquare{1}{1}{b_3}{{rgb,255:red,255; green,255; blue,255}}
    \tikzColoredSquare{2}{1}{b_4}{{rgb,255:red,255; green,255; blue,255}}
    \tikzColoredSquare{3}{1}{b_5}{{rgb,255:red,255; green,255; blue,255}}
    \tikzColoredSquare{4}{1}{\widehat 0}{{rgb,255:red,255; green,255; blue,255}}

    \tikzColoredSquare{4}{0}{3}{{rgb,255:red,240; green,50; blue,230}}
    
    \tikzColoredSquare{-1}{0}{1}{{rgb,255:red,60; green,180; blue,75}}
    \tikzColoredSquare{0}{0}{1}{{rgb,255:red,60; green,180; blue,75}}
    \tikzColoredSquare{1}{0}{0}{{rgb,255:red,245; green,130; blue,48}}
    \tikzColoredSquare{2}{0}{1}{{rgb,255:red,60; green,180; blue,75}}
    \tikzColoredSquare{3}{0}{2}{{rgb,255:red,170; green,255; blue,195}}
    
    \tikzColoredSquare{4}{-1}{3}{{rgb,255:red,170; green,255; blue,195}}
    
    \tikzColoredSquare{-1}{-1}{2}{{rgb,255:red,60; green,180; blue,75}}
    \tikzColoredSquare{0}{-1}{2}{{rgb,255:red,60; green,180; blue,75}}
    \tikzColoredSquare{1}{-1}{1}{{rgb,255:red,245; green,130; blue,48}}
    \tikzColoredSquare{2}{-1}{2}{{rgb,255:red,60; green,180; blue,75}}
    \tikzColoredSquare{3}{-1}{3}{{rgb,255:red,170; green,255; blue,195}}

    \tikzColoredSquare{4}{-2}{3}{{rgb,255:red,60; green,180; blue,75}}
    
    \tikzColoredSquare{-1}{-2}{3}{{rgb,255:red,60; green,180; blue,75}}
    \tikzColoredSquare{0}{-2}{3}{{rgb,255:red,60; green,180; blue,75}}
    \tikzColoredSquare{1}{-2}{2}{{rgb,255:red,245; green,130; blue,48}}
    \tikzColoredSquare{2}{-2}{3}{{rgb,255:red,60; green,180; blue,75}}
    \tikzColoredSquare{3}{-2}{0}{{rgb,255:red,128; green,128; blue,128}}
    
    \tikzColoredSquare{4}{-3}{3}{{rgb,255:red,245; green,130; blue,48}}
    
    \tikzColoredSquare{-1}{-3}{0}{{rgb,255:red,0; green,130; blue,200}}
    \tikzColoredSquare{0}{-3}{0}{{rgb,255:red,0; green,130; blue,200}}
    \tikzColoredSquare{1}{-3}{3}{{rgb,255:red,245; green,130; blue,48}}
    \tikzColoredSquare{2}{-3}{0}{{rgb,255:red,0; green,130; blue,200}}
    \tikzColoredSquare{3}{-3}{1}{{rgb,255:red,128; green,128; blue,128}}
    
    \tikzColoredSquare{4}{-4}{2}{{rgb,255:red,128; green,128; blue,128}}
    
    \tikzColoredSquare{-1}{-4}{1}{{rgb,255:red,0; green,130; blue,200}}
    \tikzColoredSquare{0}{-4}{1}{{rgb,255:red,0; green,130; blue,200}}
    \tikzColoredSquare{1}{-4}{0}{{rgb,255:red,220; green,190; blue,255}}
    \tikzColoredSquare{2}{-4}{1}{{rgb,255:red,0; green,130; blue,200}}
    \tikzColoredSquare{3}{-4}{2}{{rgb,255:red,128; green,128; blue,128}}
    
    \tikzColoredSquare{4}{-5}{3}{{rgb,255:red,128; green,128; blue,128}}
    
    \tikzColoredSquare{-1}{-5}{2}{{rgb,255:red,0; green,130; blue,200}}
    \tikzColoredSquare{0}{-5}{2}{{rgb,255:red,0; green,130; blue,200}}
    \tikzColoredSquare{1}{-5}{1}{{rgb,255:red,220; green,190; blue,255}}
    \tikzColoredSquare{2}{-5}{2}{{rgb,255:red,0; green,130; blue,200}}
    \tikzColoredSquare{3}{-5}{0}{{rgb,255:red,70; green,240; blue,240}}
    
    \tikzColoredSquare{4}{-6}{3}{{rgb,255:red,0; green,130; blue,200}}
    
    \tikzColoredSquare{-1}{-6}{3}{{rgb,255:red,0; green,130; blue,200}}
    \tikzColoredSquare{0}{-6}{3}{{rgb,255:red,0; green,130; blue,200}}
    \tikzColoredSquare{1}{-6}{2}{{rgb,255:red,220; green,190; blue,255}}
    \tikzColoredSquare{2}{-6}{3}{{rgb,255:red,0; green,130; blue,200}}
    \tikzColoredSquare{3}{-6}{1}{{rgb,255:red,70; green,240; blue,240}}
    
    \tikzColoredSquare{4}{-7}{3}{{rgb,255:red,220; green,190; blue,255}}
        
    \tikzColoredSquare{-1}{-7}{0}{{rgb,255:red,230; green,25; blue,75}}
    \tikzColoredSquare{0}{-7}{0}{{rgb,255:red,230; green,25; blue,75}}
    \tikzColoredSquare{1}{-7}{3}{{rgb,255:red,220; green,190; blue,255}}
    \tikzColoredSquare{2}{-7}{0}{{rgb,255:red,230; green,25; blue,75}}
    \tikzColoredSquare{3}{-7}{2}{{rgb,255:red,70; green,240; blue,240}}
    
    \tikzColoredSquare{4}{-8}{3}{{rgb,255:red,70; green,240; blue,240}}
        
    \tikzColoredSquare{-1}{-8}{1}{{rgb,255:red,230; green,25; blue,75}}
    \tikzColoredSquare{0}{-8}{1}{{rgb,255:red,230; green,25; blue,75}}
    \tikzColoredSquare{1}{-8}{0}{{rgb,255:red,255; green,255; blue,25}}
    \tikzColoredSquare{2}{-8}{1}{{rgb,255:red,230; green,25; blue,75}}
    \tikzColoredSquare{3}{-8}{3}{{rgb,255:red,70; green,240; blue,240}}
    
    \tikzColoredSquare{4}{-9}{2}{{rgb,255:red,230; green,25; blue,75}}
    
    \tikzColoredSquare{-1}{-9}{2}{{rgb,255:red,230; green,25; blue,75}}
    \tikzColoredSquare{0}{-9}{2}{{rgb,255:red,230; green,25; blue,75}}
    \tikzColoredSquare{1}{-9}{1}{{rgb,255:red,255; green,255; blue,25}}
    \tikzColoredSquare{2}{-9}{2}{{rgb,255:red,230; green,25; blue,75}}
    \tikzColoredSquare{3}{-9}{0}{{rgb,255:red,170; green,110; blue,40}}
    
    \tikzColoredSquare{4}{-10}{3}{{rgb,255:red,230; green,25; blue,75}}
        
    \tikzColoredSquare{-1}{-10}{0}{{rgb,255:red,250; green,190; blue,212}}
    \tikzColoredSquare{0}{-10}{0}{{rgb,255:red,250; green,190; blue,212}}
    \tikzColoredSquare{1}{-10}{2}{{rgb,255:red,255; green,255; blue,25}}
    \tikzColoredSquare{2}{-10}{0}{{rgb,255:red,250; green,190; blue,212}}
    \tikzColoredSquare{3}{-10}{1}{{rgb,255:red,170; green,110; blue,40}}
    
    \tikzColoredSquare{4}{-11}{3}{{rgb,255:red,255; green,255; blue,25}}
        
    \tikzColoredSquare{-1}{-11}{1}{{rgb,255:red,250; green,190; blue,212}}
    \tikzColoredSquare{0}{-11}{1}{{rgb,255:red,250; green,190; blue,212}}
    \tikzColoredSquare{1}{-11}{3}{{rgb,255:red,255; green,255; blue,25}}
    \tikzColoredSquare{2}{-11}{1}{{rgb,255:red,250; green,190; blue,212}}
    \tikzColoredSquare{3}{-11}{2}{{rgb,255:red,170; green,110; blue,40}}
    
    \tikzColoredSquare{4}{-12}{3}{{rgb,255:red,170; green,110; blue,40}}
        
    \tikzColoredSquare{-1}{-12}{2}{{rgb,255:red,250; green,190; blue,212}}
    \tikzColoredSquare{0}{-12}{2}{{rgb,255:red,250; green,190; blue,212}}
    \tikzColoredSquare{1}{-12}{0}{{rgb,255:red,240; green,50; blue,230}}
    \tikzColoredSquare{2}{-12}{2}{{rgb,255:red,250; green,190; blue,212}}
    \tikzColoredSquare{3}{-12}{3}{{rgb,255:red,170; green,110; blue,40}}
    
    \tikzColoredSquare{4}{-13}{3}{{rgb,255:red,250; green,190; blue,212}}
    
    \tikzColoredSquare{-1}{-13}{3}{{rgb,255:red,250; green,190; blue,212}}
    \tikzColoredSquare{0}{-13}{3}{{rgb,255:red,250; green,190; blue,212}}
    \tikzColoredSquare{1}{-13}{1}{{rgb,255:red,240; green,50; blue,230}}
    \tikzColoredSquare{2}{-13}{3}{{rgb,255:red,250; green,190; blue,212}}
    \tikzColoredSquare{3}{-13}{0}{{rgb,255:red,170; green,255; blue,195}}
    
    \tikzColoredSquare{4}{-14}{2}{{rgb,255:red,240; green,50; blue,230}}
        
    \tikzColoredSquare{-1}{-14}{0}{{rgb,255:red,60; green,180; blue,75}}
    \tikzColoredSquare{0}{-14}{0}{{rgb,255:red,60; green,180; blue,75}}
    \tikzColoredSquare{1}{-14}{2}{{rgb,255:red,240; green,50; blue,230}}
    \tikzColoredSquare{2}{-14}{0}{{rgb,255:red,60; green,180; blue,75}}
    \tikzColoredSquare{3}{-14}{1}{{rgb,255:red,170; green,255; blue,195}}

\end{tikzpicture}
\caption{An orbit board of whirling on $\BB_5\times[3]$.}\label{fig:Borbit}
\end{figure}
The following proposition is an analogue of Proposition~\ref{VCounting}.

\begin{prop}
The number of order ideals of $\BB_n\times[k]$ is given by 
$\displaystyle \sum_{i=1}^{k+1} i^n$.
%$|\mathcal{J}(\VV_k)| = \frac{k(k+1)(2k+1)}{6}$. 
\end{prop}
\begin{proof}
Every order ideal of $\BB_n\times[k]$ is associated with a tuple of labels of the elements of $\BB_n$.
Let $j \in [0, k]$ denote the label of $\widehat{0}$.
The possible labels for $b_1,\dots,b_{n-1},$ and $b_n$ are precisely the elements of $[0, j]$, yielding a total of $(j+1)^n$ labelings.
Summing over all possible values of $j$ produces the desired result.
\end{proof}

The orbit board in Figure~\ref{fig:Borbit} has been decomposed (canonically) into a mixture of one-tailed
and multi-tailed whorms.  Since now we list the labels at $\widehat{0}$ in the rightmost
column, the whorms are ``(right) edge-seeking'' rather than ``center-seeking''.  So each
whorm starts with one or more tails lying within the first $k$ columns, and then jumps to the
rightmost column at some point.  

Note that the first, second, and fourth columns in Figure~\ref{fig:Borbit} are identical.  This is no
accident.  If two (or more) entries among the first $n$ in a given row are the \textit{same}, then
those positions (columns) remain the same throughout the \emph{entire} orbit board.  This is
because the entries $b_1,\dots, b_n$ represent the result of whirling at \emph{incomporable}
elements of the poset $\BB_n$, so the whirls at those elements commute with one another.  
 Furthermore, these two entries must belong to the same whorm,
because each will be whorm-connected via $\widehat{0}$ exactly when their value matches the
value of the last entry. So orbit boards for whirling on $\BB_{n}$ can decompose into
multi-tailed whorms with up to $n$ tails.  (The $n$-tails case entails near-constant
rows and exactly two whorms decomposing the super-orbit board.) 

These observations will allow us to generalize our periodicity and homomesy
results from $\VV$ to $\BB_n$. First we define a map that will be equivalent to whirling
$k+2$ times.  This generalizes the ``reflection'' $(x,y,z)\mapsto (z,y,x)$, which is the result
of whirling $k+2$ times in the $\VV_{k}$ case. 

\begin{definition}\label{def:whirlA}
For any $A \subseteq [0,k]$, define the family of order-reversing maps $\mathcal{F}_k^{A}(\BB_n)
= \{f: f\in\mathcal{F}_k(\BB_n)\text{ and }f(b_j) \in A\text{ for all } j \in[n]\}$. Then
set $\overline{w}_{A}: \mathcal{F}_k^{A}(\BB_n)$ to be the map which whirls
(cyclically increments \textit{within the subset $A$}) each label on the non-$\widehat{0}$
elements of $\BB_{n}$.  By definition of whirling (repeatedly incrementing until a valid
labeling is reached), this preserves the condition that $f(b_{i})\leq f(\widehat{0})$ for all $i\in
[n]$, so is well-defined.  Since there are no relations between the leaves of $\BB_{n}$, one
can update the labels in any order.  
\end{definition}

This definition becomes particularly useful, when $A$ is precisely the set of current labels
at nonzero elements of $f$.  
\begin{definition}\label{def:whirlC}

Given $f\in \mathcal{F}_k(\BB_n)$, set $A(f) = \{a: f(b_j) = a \text{ for some } j\in[n]\}$,
i.e., $A(f)$ is the set of values that the $P$-partition $f$ attains on the
non-$\widehat{0}$ elements of $\BB_n$. Set $\alpha = \#A$ and $\alpha(f) = \# A(f)$.
For any $f,g\in\mathcal{F}_k(\BB_n)$,
if $ g=w^j(f)$ for some $j\in\NN$, then $\alpha(f)=\alpha(g)$. (Whirling preserves the
number of \textit{distinct} labels.) So we may sometimes write
just $\alpha$ when an orbit is fixed.  For the rest of this section, we will consider the map
$\overline{w}_{A(f)}$, where whirling takes place within the set of possible labels given by
the current function.  (This set remains constant throughout an orbit of
$\overline{w}_{A(f)}$.)  
\end{definition}

\begin{ex}
Consider $f=(1,3,3,0,4,1,6)\in \mathcal{F}_9(\BB_6)$. We see $A(f) = \{0,1,3,4\}$ so 
$$\overline{w}_{A(f)}(1,3,3,0,4,1,6) = (3,4,4,1,0,3,6).$$
The last entry remains unchanged, and the earlier entries are increasing cyclically within
the set $A(f)=\{0,1,3,4\}$, with $\alpha =4$.  As another example, the reader can compare
rows five apart in the orbit board of Figure~\ref{fig:Borbit}.  
\end{ex}

In the special case of $\VV$ ($=\BB_2$) our set $A=A(f)$ within any
orbit will have at most two elements, hence $\overline{w}_{A(f)}$ will just toggle
between those two values at the left and the right. This means that $\overline{w}_{A(f)}$ is
the same as reflecting values across the center of the orbit board, which we already saw was the
effect of $w^{k+2}$. Our next theorem will generalize this to the case of $\BB_n$, except
now the center column has moved to be the rightmost one.  

\subsection{Periodicity and homomesy via edge-seeking whorms}\label{ss:clawmain}

Our goal now is to generalize the periodicity and homomesy results for $\VV \times
[k]$ to $\BB_{n} \times [k]$.  Our first result shows that the sum of the tail-lengths
of the first $\alpha (f) + 1$ whorms (counting the length only once for multi-tailed whorms)
is constant.  If $f\in \mathcal{F}^k(\BB_n)$ satisfies $f(\widehat{0})\not\in A(f)$, then $f$ will
contain entries from $\alpha+1$ distinct whorms. Otherwise, $f$ will intersect exactly
$\alpha $ whorms.  

Define $\tl(\xi) = 1+ \text{min}\{f(\hat 0): (\hat 0,f)\in \xi\}$ and
$\hl(\xi): = k+2 - \tl(\xi)$.  If there exists $(\hat 0, f) \in \xi_1$
with $f(\hat 0) = k$ such that $(\hat 0, w(f))\in \xi_2$, then we say $\xi_2$ is
\emph{in front of} $\xi_1$. (In other words, if their heads are consecutive in the
rightmost column of the orbit board.)  In Figure $\ref{fig:Borbit}$, the yellow whorm is in
front of the red whorm. We call a sequence of whorms \emph{consecutive} if each is in front
of the next.  The next result generalizes~Lemma~\ref{lem:bp_consecutive}.  

\begin{lem}\label{lem:Bconwhorm}
Fix a whirling orbit of $\mathcal{F}_k(\BB_n)$, decomposed into whorms, with $\alpha = \alpha (f)$ for any $f$ in the
orbit.
If $\xi_1,\dots,\xi_{\alpha+1}$ are $\alpha+1$ consecutive whorms, then
\begin{equation}\label{eq:tailsum}
\tl(\xi_1) + \cdots + \tl(\xi_{\alpha+1})  = \alpha(k+2).
\end{equation}
Using $\hl(\xi)+ \tl(\xi) = k+2$, this can be written equivalently as

\begin{equation}\label{eq:headsum}
\hl(\xi_1) + \cdots + \hl(\xi_{\alpha+1})  = k+2.
\end{equation}
or as
\begin{equation}\label{eq:headsortails}
\tl(\xi_{\alpha+1}) = \hl(\xi_1)+\dots+\hl(\xi_\alpha).
\end{equation}

\end{lem}

\begin{proof}
The proof follows similarly to that of Lemma~\ref{lem:bp_consecutive}. 
The key observation is that $\tl (\xi_{\alpha+1}) =\hl(\xi_1) +\cdots +\hl(\xi_\alpha)$ for consecutive whorms.  
This holds because the tail of $\xi_{\alpha+1}$ begins one row below where the head of $\xi_{1}$ starts, and ends one row below where $\xi_{\alpha}$ ends.  
But since $\tl (\xi ) = k+2 - \hl (\xi )$
for any $\xi $ we get
\[
\tl(\xi_{\alpha+1}) = \hl(\xi_1) + \cdots + \hl(\xi_\alpha) = (k+2-\tl (\xi_{1}))+\dots + (k+2-\tl (\xi_{\alpha})),
\]
which implies the result. See Figure~\ref{fig:Borbit} for examples.\end{proof}

\begin{cor}\label{cor:Bwhormrepeat}
Fix a whirling orbit of $\mathcal{F}_k(\BB_n)$, decomposed into whorms, with
$\alpha = \alpha (f)$ for any $f$ in the orbit.
If $\xi_1,\dots,\xi_{\alpha+2}$ are consecutive whorms, then 
$\tl(\xi_1) = \tl(\xi_{\alpha+2})$.  
\end{cor}

\begin{proof}
Applying Lemma~\ref{lem:Bconwhorm} to two consecutive consecutive collections of whorms, we obtain
$$ \tl(\xi_1) + \cdots + \tl(\xi_{\alpha+1})  = \alpha(k+2)= \tl(\xi_2) + \cdots + \tl(\xi_{\alpha+2}),$$
which reduces to $\tl(\xi_1) = \tl(\xi_{\alpha+2})$.
\end{proof}

The next theorem, which generalizes Theorem~\ref{thm:halfwaysymm}, gives a simple way of
describing the $(k+2)$nd iteration of whirling as the much simpler action of whirling (in any order) at each leaf
of $\BB_{n}$ within the set of allowable labels $A(f)$. 

\begin{figure}
    
\begin{tikzpicture}[baseline={([yshift=-.8ex]current bounding box.center)}]
    
    \tikzColoredSquare{-4.25}{1}{\text{columns:} }{{rgb,255:red,255; green,255; blue,255}}
    \tikzColoredSquare{-2}{1}{b_1 }{{rgb,255:red,255; green,255; blue,255}}
    \tikzColoredSquare{-1}{1}{\dots}{{rgb,255:red,255; green,255; blue,255}}
    \tikzColoredSquare{0}{1}{b_i}{{rgb,255:red,255; green,255; blue,255}}
    \tikzColoredSquare{1}{1}{\dots}{{rgb,255:red,255; green,255; blue,255}}
    \tikzColoredSquare{2}{1}{b_j}{{rgb,255:red,255; green,255; blue,255}}
    \tikzColoredSquare{3}{1}{\dots}{{rgb,255:red,255; green,255; blue,255}}
    \tikzColoredSquare{4}{1}{\widehat 0}{{rgb,255:red,255; green,255; blue,255}}
    
    \def\m{0}
    \tikzColoredSquare{-2}{\m}{ }{{rgb,255:red,200; green,200; blue,200}}
    \tikzColoredSquare{-1}{\m}{ }{{rgb,255:red,200; green,200; blue,200}}
    \tikzColoredSquare{0}{\m}{ }{{rgb,255:red,200; green,200; blue,200}}
    \tikzColoredSquare{1}{\m}{ }{{rgb,255:red,200; green,200; blue,200}}
    \tikzColoredSquare{2}{\m}{ }{{rgb,255:red,200; green,200; blue,200}}
    \tikzColoredSquare{3}{\m}{ }{{rgb,255:red,200; green,200; blue,200}}
    \tikzColoredSquare{4}{\m}{ }{{rgb,255:red,200; green,200; blue,200}}
    \def\m{-1}
    \tikzColoredSquare{-2}{\m}{ }{{rgb,255:red,200; green,200; blue,200}}
    \tikzColoredSquare{-1}{\m}{ }{{rgb,255:red,200; green,200; blue,200}}
    \tikzColoredSquare{0}{\m}{ }{{rgb,255:red,200; green,200; blue,200}}
    \tikzColoredSquare{1}{\m}{ }{{rgb,255:red,200; green,200; blue,200}}
    \tikzColoredSquare{2}{\m}{ }{{rgb,255:red,200; green,200; blue,200}}
    \tikzColoredSquare{3}{\m}{ }{{rgb,255:red,200; green,200; blue,200}}
    \tikzColoredSquare{4}{\m}{ }{{rgb,255:red,200; green,200; blue,200}}
    \def\m{-2}
    \tikzColoredSquare{-2}{\m}{ }{{rgb,255:red,200; green,200; blue,200}}
    \tikzColoredSquare{-1}{\m}{ }{{rgb,255:red,200; green,200; blue,200}}
    \tikzColoredSquare{0}{\m}{ }{{rgb,255:red,200; green,200; blue,200}}
    \tikzColoredSquare{1}{\m}{ }{{rgb,255:red,200; green,200; blue,200}}
    \tikzColoredSquare{2}{\m}{ }{{rgb,255:red,200; green,200; blue,200}}
    \tikzColoredSquare{3}{\m}{ }{{rgb,255:red,200; green,200; blue,200}}
    \tikzColoredSquare{4}{\m}{ }{{rgb,255:red,200; green,200; blue,200}}
    \def\m{-3}
    \tikzColoredSquare{-2}{\m}{ }{{rgb,255:red,200; green,200; blue,200}}
    \tikzColoredSquare{-1}{\m}{ }{{rgb,255:red,200; green,200; blue,200}}
    \tikzColoredSquare{0}{\m}{ }{{rgb,255:red,200; green,200; blue,200}}
    \tikzColoredSquare{1}{\m}{ }{{rgb,255:red,200; green,200; blue,200}}
    \tikzColoredSquare{2}{\m}{ }{{rgb,255:red,200; green,200; blue,200}}
    \tikzColoredSquare{3}{\m}{ }{{rgb,255:red,200; green,200; blue,200}}
    \tikzColoredSquare{4}{\m}{ }{{rgb,255:red,200; green,200; blue,200}}
    \def\m{-4}
    \tikzColoredSquare{-2}{\m}{ }{{rgb,255:red,200; green,200; blue,200}}
    \tikzColoredSquare{-1}{\m}{ }{{rgb,255:red,200; green,200; blue,200}}
    \tikzColoredSquare{0}{\m}{ }{{rgb,255:red,200; green,200; blue,200}}
    \tikzColoredSquare{1}{\m}{ }{{rgb,255:red,200; green,200; blue,200}}
    \tikzColoredSquare{2}{\m}{ }{{rgb,255:red,200; green,200; blue,200}}
    \tikzColoredSquare{3}{\m}{ }{{rgb,255:red,200; green,200; blue,200}}
    \tikzColoredSquare{4}{\m}{ }{{rgb,255:red,200; green,200; blue,200}}
    \def\m{-5}
    \tikzColoredSquare{-2}{\m}{ }{{rgb,255:red,200; green,200; blue,200}}
    \tikzColoredSquare{-1}{\m}{ }{{rgb,255:red,200; green,200; blue,200}}
    \tikzColoredSquare{0}{\m}{ }{{rgb,255:red,200; green,200; blue,200}}
    \tikzColoredSquare{1}{\m}{ }{{rgb,255:red,200; green,200; blue,200}}
    \tikzColoredSquare{2}{\m}{ }{{rgb,255:red,200; green,200; blue,200}}
    \tikzColoredSquare{3}{\m}{ }{{rgb,255:red,200; green,200; blue,200}}
    \tikzColoredSquare{4}{\m}{ }{{rgb,255:red,200; green,200; blue,200}}
    \def\m{-6}
    \tikzColoredSquare{-2}{\m}{ }{{rgb,255:red,200; green,200; blue,200}}
    \tikzColoredSquare{-1}{\m}{ }{{rgb,255:red,200; green,200; blue,200}}
    \tikzColoredSquare{0}{\m}{ }{{rgb,255:red,200; green,200; blue,200}}
    \tikzColoredSquare{1}{\m}{ }{{rgb,255:red,200; green,200; blue,200}}
    \tikzColoredSquare{2}{\m}{ }{{rgb,255:red,200; green,200; blue,200}}
    \tikzColoredSquare{3}{\m}{ }{{rgb,255:red,200; green,200; blue,200}}
    \tikzColoredSquare{4}{\m}{ }{{rgb,255:red,200; green,200; blue,200}}
    \def\m{-7}
    \tikzColoredSquare{-2}{\m}{ }{{rgb,255:red,200; green,200; blue,200}}
    \tikzColoredSquare{-1}{\m}{ }{{rgb,255:red,200; green,200; blue,200}}
    \tikzColoredSquare{0}{\m}{ }{{rgb,255:red,200; green,200; blue,200}}
    \tikzColoredSquare{1}{\m}{ }{{rgb,255:red,200; green,200; blue,200}}
    \tikzColoredSquare{2}{\m}{ }{{rgb,255:red,200; green,200; blue,200}}
    \tikzColoredSquare{3}{\m}{ }{{rgb,255:red,200; green,200; blue,200}}
    \tikzColoredSquare{4}{\m}{ }{{rgb,255:red,200; green,200; blue,200}}
    \def\m{-8}
    \tikzColoredSquare{-2}{\m}{ }{{rgb,255:red,255; green,255; blue,255}}
    \tikzColoredSquare{-1}{\m}{ }{{rgb,255:red,255; green,255; blue,255}}
    \tikzColoredSquare{0}{\m}{ }{{rgb,255:red,255; green,255; blue,255}}
    \tikzColoredSquare{1}{\m}{ }{{rgb,255:red,255; green,255; blue,255}}
    \tikzColoredSquare{2}{\m}{ }{{rgb,255:red,255; green,255; blue,255}}
    \tikzColoredSquare{3}{\m}{ }{{rgb,255:red,255; green,255; blue,255}}
    \tikzColoredSquare{4}{\m}{ }{{rgb,255:red,255; green,255; blue,255}}
    \def\m{-9}
    \tikzColoredSquare{-2}{\m}{ }{{rgb,255:red,200; green,200; blue,200}}
    \tikzColoredSquare{-1}{\m}{ }{{rgb,255:red,200; green,200; blue,200}}
    \tikzColoredSquare{0}{\m}{ }{{rgb,255:red,200; green,200; blue,200}}
    \tikzColoredSquare{1}{\m}{ }{{rgb,255:red,200; green,200; blue,200}}
    \tikzColoredSquare{2}{\m}{ }{{rgb,255:red,200; green,200; blue,200}}
    \tikzColoredSquare{3}{\m}{ }{{rgb,255:red,200; green,200; blue,200}}
    \tikzColoredSquare{4}{\m}{ }{{rgb,255:red,200; green,200; blue,200}}
    \def\m{-10}
    \tikzColoredSquare{-2}{\m}{ }{{rgb,255:red,200; green,200; blue,200}}
    \tikzColoredSquare{-1}{\m}{ }{{rgb,255:red,200; green,200; blue,200}}
    \tikzColoredSquare{0}{\m}{ }{{rgb,255:red,200; green,200; blue,200}}
    \tikzColoredSquare{1}{\m}{ }{{rgb,255:red,200; green,200; blue,200}}
    \tikzColoredSquare{2}{\m}{ }{{rgb,255:red,200; green,200; blue,200}}
    \tikzColoredSquare{3}{\m}{ }{{rgb,255:red,200; green,200; blue,200}}
    \tikzColoredSquare{4}{\m}{ }{{rgb,255:red,200; green,200; blue,200}}
    \def\m{-11}
    \tikzColoredSquare{-2}{\m}{ }{{rgb,255:red,200; green,200; blue,200}}
    \tikzColoredSquare{-1}{\m}{ }{{rgb,255:red,200; green,200; blue,200}}
    \tikzColoredSquare{0}{\m}{ }{{rgb,255:red,200; green,200; blue,200}}
    \tikzColoredSquare{1}{\m}{ }{{rgb,255:red,200; green,200; blue,200}}
    \tikzColoredSquare{2}{\m}{ }{{rgb,255:red,200; green,200; blue,200}}
    \tikzColoredSquare{3}{\m}{ }{{rgb,255:red,200; green,200; blue,200}}
    \tikzColoredSquare{4}{\m}{ }{{rgb,255:red,200; green,200; blue,200}}
    
    \tikzColoredSquare{2}{0}{\rotatebox{270}{$\hspace{.1cm}\xi_1$} }{{rgb,255:red,60; green,180; blue,75}}

    \tikzColoredSquare{0}{-1}{\rotatebox{270}{$\hspace{.1cm}\xi_2$} }{{rgb,255:red,245; green,130; blue,48}}
    \tikzColoredSquare{2}{-1}{ }{{rgb,255:red,60; green,180; blue,75}}

    \tikzColoredSquare{0}{-2}{ }{{rgb,255:red,245; green,130; blue,48}}
    \tikzColoredSquare{2}{-2}{ }{{rgb,255:red,60; green,180; blue,75}}

    \tikzColoredSquare{0}{-2}{ }{{rgb,255:red,245; green,130; blue,48}}
    \tikzColoredSquare{2}{-2}{t_1}{{rgb,255:red,60; green,180; blue,75}}
    \tikzColoredSquare{4}{-2}{t_1}{{rgb,255:red,60; green,180; blue,75}}
    
    \tikzColoredSquare{0}{-3}{ }{{rgb,255:red,245; green,130; blue,48}}
    \tikzColoredSquare{2}{-3}{ }{{rgb,255:red,0; green,130; blue,200}}
    \tikzColoredSquare{4}{-3}{ }{{rgb,255:red,60; green,180; blue,75}}
    
    \tikzColoredSquare{0}{-4}{ }{{rgb,255:red,245; green,130; blue,48}}
    \tikzColoredSquare{2}{-4}{ }{{rgb,255:red,0; green,130; blue,200}}
    \tikzColoredSquare{4}{-4}{ }{{rgb,255:red,60; green,180; blue,75}}
    
    \tikzColoredSquare{0}{-5}{ }{{rgb,255:red,245; green,130; blue,48}}
    \tikzColoredSquare{2}{-5}{ }{{rgb,255:red,0; green,130; blue,200}}
    \tikzColoredSquare{4}{-5}{ }{{rgb,255:red,245; green,130; blue,48}}
    
    \tikzColoredSquare{0}{-6}{ }{{rgb,255:red,230; green,25; blue,75}}
    \tikzColoredSquare{2}{-6}{ }{{rgb,255:red,0; green,130; blue,200}}
    \tikzColoredSquare{4}{-6}{ }{{rgb,255:red,245; green,130; blue,48}}
 
    \tikzColoredSquare{1}{-7.75}{\vdots}{{rgb,255:red,255; green,255; blue,255}}
    \tikzColoredSquare{0}{-7}{ }{{rgb,255:red,230; green,25; blue,75}}
    \tikzColoredSquare{1}{-7}{ }{{rgb,255:red,200; green,200; blue,200}}
    
   \tikzColoredSquare{2}{-7}{ }{{rgb,255:red,0; green,130; blue,200}}
    
    \tikzColoredSquare{0}{-9}{ }{{rgb,255:red,230; green,25; blue,75}}
    \tikzColoredSquare{2}{-9}{ }{{rgb,255:red,200; green,200; blue,200}}
    \tikzColoredSquare{4}{-9}{ }{{rgb,255:red,0; green,130; blue,200}}
    
    \tikzColoredSquare{0}{-10}{ }{{rgb,255:red,230; green,25; blue,75}}
    \tikzColoredSquare{2}{-10}{ }{{rgb,255:red,200; green,200; blue,200}}
    \tikzColoredSquare{4}{-10}{ }{{rgb,255:red,0; green,130; blue,200}}
    
    \tikzColoredSquare{0}{-11}{t_1}{{rgb,255:red,230; green,25; blue,75}}
    \tikzColoredSquare{4}{-11}{t_1}{{rgb,255:red,230; green,25; blue,75}}
    
    \tikzColoredSquare{2}{-3}{\rotatebox{270}{$\hspace{.5cm}\xi_{\alpha+1}$} }{{rgb,255:red,0; green,130; blue,200}}
    \tikzColoredSquare{0}{-6}{\rotatebox{270}{$\hspace{.5cm}\xi_{\alpha+2}$} }{{rgb,255:red,230; green,25; blue,75}}
    \draw[-] (2.30,-0.75) -- (2.70,-0.75);
    \draw[-] (2.5,-0.75) -- (2.5,-5.25);
    \draw[-] (2.30,-5.25) -- (2.70,-5.25);
    \tikzColoredSquare{6}{-6}{k+2}{{rgb,255:red,255; green,255; blue,255}}

    \draw[-] (-1.25,-5.75) -- (2.25,-5.75);
    \draw[-] (2.25,-5.75) -- (2.25,.25);
    \draw[-] (2.25,.25) -- (-1.25,.25);
    \draw[-] (-1.25,.25) -- (-1.25,-5.75);
    
\end{tikzpicture}
\caption{Sketch for the proof of Theorem~\ref{thm:worder}. Consecutive whorms are denoted by
$\xi_1,\dots,\xi_{\alpha+2}$. The last label in the tail of $\xi_1$ is $t_1$. Note that the tail
lengths of $\xi_{\alpha+2}$ and $\xi_1$, are the same, so the labels in {\color{rgb,255:red,60; green,180; blue,75} column $b_j$}
reproduce those in {\color{rgb,255:red,230; green,25; blue,75} column $b_i$} exactly $k+2$ rows further along the orbit.}\label{fig:worder} 
\end{figure}
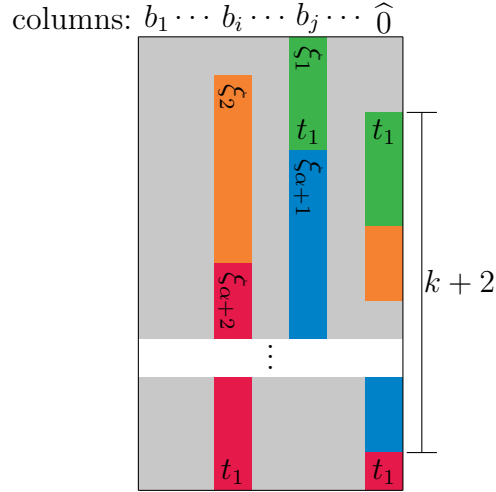

\begin{thm}\label{thm:worder}
Let $w$ be the whirling action on $k$-bounded $P$-partitions in $\mathcal{F}_k(\BB_n)$. 
For any $f\in \mathcal{F}_k(\BB_n)$ with $A=A(f)$ and $\alpha =\alpha (f) = \#A(f)$, 
we have $w^{k+2}f=\overline{w}_{A(f)}f$. Thus, $w^{\alpha (k+2)}f=f$.
\end{thm}

\begin{proof}
Let $\mathcal{R}$ be an orbit board of $w$, decomposed into whorms, with
$A=A(f)$ for some (row) $f$ in the orbit. Let $\xi_1,\dots,\xi_{\alpha}$ be
the consecutive whorms whose tails intersect $f$. Because the whorms are consecutive, in row $f$
the largest value in $A$ belong to $\xi_1$, the next largest to $\xi_2$, and so on. 
In this proof we will show that the labels in the tail(s) of $\xi_1$ will appear in
the tail(s) of $\xi_{\alpha +2}$ exactly $k+2$ rows laters. This process can be iterated to
include all the whorms, proving that the labels change as described by  
$\overline{w}_{A(f)}f$.

By Lemma~\ref{lem:Bconwhorm}~(3), the sum of the head lengths of the whorms $\xi_{1},
\xi_{2},\dotsc ,\xi_{\alpha +1}$ is exactly $k+2$.  Thus, the head
of $\xi_{\alpha +2}$ begins exactly $k+2$ rows below the row where the head  of
$\xi_{1}$ begins. 
Equivalently, the tail of
$\xi_{\alpha +2}$ \emph{ends} in the row exactly $k+2$ below 
where the tail of $\xi_{1}$ \emph{ends}.  
Since these two whorms have the \emph{same} tail length by Corollary~\ref{cor:Bwhormrepeat},
we get that each entry within each tail of $\xi_{\alpha +2}$ is the same as each entry
within each tail of $\xi_{1}$, cyclically shifted down $k+2$ rows and appearing in different
column(s). Note also that the tail column(s) of $\xi_{\alpha +2}$ are the \emph{same}
as the tail column(s) of $\xi_{2}$.  So in any of these ($\xi_{2}$) columns, the label $k+2$ rows
later has become the label from \textbf{the previous whorm} $\xi_{1}$ in row $f$, which is the cyclically next
larger value within $A(f)$. 
Thus, $w^{k+2}f=\overline{w}_{A(f)}f$, proving the theorem. 
\end{proof}

\begin{ex}
Here is the orbit board generated by $(0,0,2,5,0,6)\in \mathcal F_6(\BB_5)$, with $\alpha=3$
and divided into three sections of length $k+2$. In each section one of the whorms
$\xi_1,\xi_5,$ and $\xi_9$ is highlighted. Each whorm starts in
the first row and ends in the penultimate row of its section. The
tail-lengths of $\xi_1$ and $\xi_{\alpha+ 2}$ are the same, just as in
Corollary~\ref{cor:Bwhormrepeat}. Reading across from one section to the next (cyclically)
in each row, illustrates Theorem~\ref{thm:worder}:  $w^{k+2}=\overline{w}_{A}$. 
\[
\begin{tikzpicture}[baseline={([yshift=-.8ex]current bounding box.center)}]
	\def\m{0}
    \tikzColoredSquare{-1}{\m}{0}{{rgb,255:red,230; green,25; blue,75}}
    \tikzColoredSquare{ 0}{\m}{0}{{rgb,255:red,230; green,25; blue,75}}
    \tikzColoredSquare{ 1}{\m}{2}{{rgb,255:red,255; green,255; blue,255}}
    \tikzColoredSquare{ 2}{\m}{5}{{rgb,255:red,255; green,255; blue,255}}
    \tikzColoredSquare{ 3}{\m}{0}{{rgb,255:red,230; green,25; blue,75}}
    \tikzColoredSquare{ 4}{\m}{6}{{rgb,255:red,255; green,255; blue,255}}
	\def\m{-1}
    \tikzColoredSquare{-1}{\m}{1}{{rgb,255:red,230; green,25; blue,75}}
    \tikzColoredSquare{ 0}{\m}{1}{{rgb,255:red,230; green,25; blue,75}}
    \tikzColoredSquare{ 1}{\m}{3}{{rgb,255:red,255; green,255; blue,255}}
    \tikzColoredSquare{ 2}{\m}{6}{{rgb,255:red,255; green,255; blue,255}}
    \tikzColoredSquare{ 3}{\m}{1}{{rgb,255:red,230; green,25; blue,75}}
    \tikzColoredSquare{ 4}{\m}{6}{{rgb,255:red,255; green,255; blue,255}}
	\def\m{-2}
    \tikzColoredSquare{-1}{\m}{2}{{rgb,255:red,230; green,25; blue,75}}
    \tikzColoredSquare{ 0}{\m}{2}{{rgb,255:red,230; green,25; blue,75}}
    \tikzColoredSquare{ 1}{\m}{4}{{rgb,255:red,255; green,255; blue,255}}
    \tikzColoredSquare{ 2}{\m}{0}{{rgb,255:red,255; green,255; blue,255}}
    \tikzColoredSquare{ 3}{\m}{2}{{rgb,255:red,230; green,25; blue,75}}
    \tikzColoredSquare{ 4}{\m}{4}{{rgb,255:red,255; green,255; blue,255}}
	\def\m{-3}
    \tikzColoredSquare{-1}{\m}{3}{{rgb,255:red,230; green,25; blue,75}}
    \tikzColoredSquare{ 0}{\m}{3}{{rgb,255:red,230; green,25; blue,75}}
    \tikzColoredSquare{ 1}{\m}{0}{{rgb,255:red,255; green,255; blue,255}}
    \tikzColoredSquare{ 2}{\m}{1}{{rgb,255:red,255; green,255; blue,255}}
    \tikzColoredSquare{ 3}{\m}{3}{{rgb,255:red,230; green,25; blue,75}}
    \tikzColoredSquare{ 4}{\m}{5}{{rgb,255:red,255; green,255; blue,255}}
	\def\m{-4}
    \tikzColoredSquare{-1}{\m}{4}{{rgb,255:red,230; green,25; blue,75}}
    \tikzColoredSquare{ 0}{\m}{4}{{rgb,255:red,230; green,25; blue,75}}
    \tikzColoredSquare{ 1}{\m}{1}{{rgb,255:red,255; green,255; blue,255}}
    \tikzColoredSquare{ 2}{\m}{2}{{rgb,255:red,255; green,255; blue,255}}
    \tikzColoredSquare{ 3}{\m}{4}{{rgb,255:red,230; green,25; blue,75}}
    \tikzColoredSquare{ 4}{\m}{6}{{rgb,255:red,255; green,255; blue,255}}
	\def\m{-5}
    \tikzColoredSquare{-1}{\m}{5}{{rgb,255:red,230; green,25; blue,75}}
    \tikzColoredSquare{ 0}{\m}{5}{{rgb,255:red,230; green,25; blue,75}}
    \tikzColoredSquare{ 1}{\m}{2}{{rgb,255:red,255; green,255; blue,255}}
    \tikzColoredSquare{ 2}{\m}{3}{{rgb,255:red,255; green,255; blue,255}}
    \tikzColoredSquare{ 3}{\m}{5}{{rgb,255:red,230; green,25; blue,75}}
    \tikzColoredSquare{ 4}{\m}{5}{{rgb,255:red,230; green,25; blue,75}}
	\def\m{-6}
    \tikzColoredSquare{-1}{\m}{0}{{rgb,255:red,255; green,255; blue,255}}
    \tikzColoredSquare{ 0}{\m}{0}{{rgb,255:red,255; green,255; blue,255}}
    \tikzColoredSquare{ 1}{\m}{3}{{rgb,255:red,255; green,255; blue,255}}
    \tikzColoredSquare{ 2}{\m}{4}{{rgb,255:red,255; green,255; blue,255}}
    \tikzColoredSquare{ 3}{\m}{0}{{rgb,255:red,255; green,255; blue,255}}
    \tikzColoredSquare{ 4}{\m}{6}{{rgb,255:red,230; green,25; blue,75}}
	\def\m{-7}
    \tikzColoredSquare{-1}{\m}{1}{{rgb,255:red,255; green,255; blue,255}}
    \tikzColoredSquare{ 0}{\m}{1}{{rgb,255:red,255; green,255; blue,255}}
    \tikzColoredSquare{ 1}{\m}{4}{{rgb,255:red,255; green,255; blue,255}}
    \tikzColoredSquare{ 2}{\m}{5}{{rgb,255:red,255; green,255; blue,255}}
    \tikzColoredSquare{ 3}{\m}{1}{{rgb,255:red,255; green,255; blue,255}}
    \tikzColoredSquare{ 4}{\m}{5}{{rgb,255:red,255; green,255; blue,255}}
\end{tikzpicture}
\ \ \ \ 
\begin{tikzpicture}[baseline={([yshift=-.8ex]current bounding box.center)}]
	\def\m{0}
    \tikzColoredSquare{-1}{\m}{2}{{rgb,255:red,255; green,255; blue,255}}
    \tikzColoredSquare{ 0}{\m}{2}{{rgb,255:red,255; green,255; blue,255}}
    \tikzColoredSquare{ 1}{\m}{5}{{rgb,255:red,255; green,255; blue,255}}
    \tikzColoredSquare{ 2}{\m}{0}{{rgb,255:red,230; green,25; blue,75}}
    \tikzColoredSquare{ 3}{\m}{2}{{rgb,255:red,255; green,255; blue,255}}
    \tikzColoredSquare{ 4}{\m}{6}{{rgb,255:red,255; green,255; blue,255}}
	\def\m{-1}
    \tikzColoredSquare{-1}{\m}{3}{{rgb,255:red,255; green,255; blue,255}}
    \tikzColoredSquare{ 0}{\m}{3}{{rgb,255:red,255; green,255; blue,255}}
    \tikzColoredSquare{ 1}{\m}{6}{{rgb,255:red,255; green,255; blue,255}}
    \tikzColoredSquare{ 2}{\m}{1}{{rgb,255:red,230; green,25; blue,75}}
    \tikzColoredSquare{ 3}{\m}{3}{{rgb,255:red,255; green,255; blue,255}}
    \tikzColoredSquare{ 4}{\m}{6}{{rgb,255:red,255; green,255; blue,255}}
	\def\m{-2}
    \tikzColoredSquare{-1}{\m}{4}{{rgb,255:red,255; green,255; blue,255}}
    \tikzColoredSquare{ 0}{\m}{4}{{rgb,255:red,255; green,255; blue,255}}
    \tikzColoredSquare{ 1}{\m}{0}{{rgb,255:red,255; green,255; blue,255}}
    \tikzColoredSquare{ 2}{\m}{2}{{rgb,255:red,230; green,25; blue,75}}
    \tikzColoredSquare{ 3}{\m}{4}{{rgb,255:red,255; green,255; blue,255}}
    \tikzColoredSquare{ 4}{\m}{4}{{rgb,255:red,255; green,255; blue,255}}
	\def\m{-3}
    \tikzColoredSquare{-1}{\m}{0}{{rgb,255:red,255; green,255; blue,255}}
    \tikzColoredSquare{ 0}{\m}{0}{{rgb,255:red,255; green,255; blue,255}}
    \tikzColoredSquare{ 1}{\m}{1}{{rgb,255:red,255; green,255; blue,255}}
    \tikzColoredSquare{ 2}{\m}{3}{{rgb,255:red,230; green,25; blue,75}}
    \tikzColoredSquare{ 3}{\m}{0}{{rgb,255:red,255; green,255; blue,255}}
    \tikzColoredSquare{ 4}{\m}{5}{{rgb,255:red,255; green,255; blue,255}}
	\def\m{-4}
    \tikzColoredSquare{-1}{\m}{1}{{rgb,255:red,255; green,255; blue,255}}
    \tikzColoredSquare{ 0}{\m}{1}{{rgb,255:red,255; green,255; blue,255}}
    \tikzColoredSquare{ 1}{\m}{2}{{rgb,255:red,255; green,255; blue,255}}
    \tikzColoredSquare{ 2}{\m}{4}{{rgb,255:red,230; green,25; blue,75}}
    \tikzColoredSquare{ 3}{\m}{1}{{rgb,255:red,255; green,255; blue,255}}
    \tikzColoredSquare{ 4}{\m}{6}{{rgb,255:red,255; green,255; blue,255}}
	\def\m{-5}
    \tikzColoredSquare{-1}{\m}{2}{{rgb,255:red,255; green,255; blue,255}}
    \tikzColoredSquare{ 0}{\m}{2}{{rgb,255:red,255; green,255; blue,255}}
    \tikzColoredSquare{ 1}{\m}{3}{{rgb,255:red,255; green,255; blue,255}}
    \tikzColoredSquare{ 2}{\m}{5}{{rgb,255:red,230; green,25; blue,75}}
    \tikzColoredSquare{ 3}{\m}{2}{{rgb,255:red,255; green,255; blue,255}}
    \tikzColoredSquare{ 4}{\m}{5}{{rgb,255:red,230; green,25; blue,75}}
	\def\m{-6}
    \tikzColoredSquare{-1}{\m}{3}{{rgb,255:red,255; green,255; blue,255}}
    \tikzColoredSquare{ 0}{\m}{3}{{rgb,255:red,255; green,255; blue,255}}
    \tikzColoredSquare{ 1}{\m}{4}{{rgb,255:red,255; green,255; blue,255}}
    \tikzColoredSquare{ 2}{\m}{0}{{rgb,255:red,255; green,255; blue,255}}
    \tikzColoredSquare{ 3}{\m}{3}{{rgb,255:red,255; green,255; blue,255}}
    \tikzColoredSquare{ 4}{\m}{6}{{rgb,255:red,230; green,25; blue,75}}
	\def\m{-7}
    \tikzColoredSquare{-1}{\m}{4}{{rgb,255:red,255; green,255; blue,255}}
    \tikzColoredSquare{ 0}{\m}{4}{{rgb,255:red,255; green,255; blue,255}}
    \tikzColoredSquare{ 1}{\m}{5}{{rgb,255:red,255; green,255; blue,255}}
    \tikzColoredSquare{ 2}{\m}{1}{{rgb,255:red,255; green,255; blue,255}}
    \tikzColoredSquare{ 3}{\m}{4}{{rgb,255:red,255; green,255; blue,255}}
    \tikzColoredSquare{ 4}{\m}{5}{{rgb,255:red,255; green,255; blue,255}}
\end{tikzpicture}
\ \ \ \ 
\begin{tikzpicture}[baseline={([yshift=-.8ex]current bounding box.center)}]
	\def\m{0}
    \tikzColoredSquare{-1}{\m}{5}{{rgb,255:red,255; green,255; blue,255}}
    \tikzColoredSquare{ 0}{\m}{5}{{rgb,255:red,255; green,255; blue,255}}
    \tikzColoredSquare{ 1}{\m}{0}{{rgb,255:red,230; green,25; blue,75}}
    \tikzColoredSquare{ 2}{\m}{2}{{rgb,255:red,255; green,255; blue,255}}
    \tikzColoredSquare{ 3}{\m}{5}{{rgb,255:red,255; green,255; blue,255}}
    \tikzColoredSquare{ 4}{\m}{6}{{rgb,255:red,255; green,255; blue,255}}
	\def\m{-1}
    \tikzColoredSquare{-1}{\m}{6}{{rgb,255:red,255; green,255; blue,255}}
    \tikzColoredSquare{ 0}{\m}{6}{{rgb,255:red,255; green,255; blue,255}}
    \tikzColoredSquare{ 1}{\m}{1}{{rgb,255:red,230; green,25; blue,75}}
    \tikzColoredSquare{ 2}{\m}{3}{{rgb,255:red,255; green,255; blue,255}}
    \tikzColoredSquare{ 3}{\m}{6}{{rgb,255:red,255; green,255; blue,255}}
    \tikzColoredSquare{ 4}{\m}{6}{{rgb,255:red,255; green,255; blue,255}}
	\def\m{-2}
    \tikzColoredSquare{-1}{\m}{0}{{rgb,255:red,255; green,255; blue,255}}
    \tikzColoredSquare{ 0}{\m}{0}{{rgb,255:red,255; green,255; blue,255}}
    \tikzColoredSquare{ 1}{\m}{2}{{rgb,255:red,230; green,25; blue,75}}
    \tikzColoredSquare{ 2}{\m}{4}{{rgb,255:red,255; green,255; blue,255}}
    \tikzColoredSquare{ 3}{\m}{0}{{rgb,255:red,255; green,255; blue,255}}
    \tikzColoredSquare{ 4}{\m}{4}{{rgb,255:red,255; green,255; blue,255}}
	\def\m{-3}
    \tikzColoredSquare{-1}{\m}{1}{{rgb,255:red,255; green,255; blue,255}}
    \tikzColoredSquare{ 0}{\m}{1}{{rgb,255:red,255; green,255; blue,255}}
    \tikzColoredSquare{ 1}{\m}{3}{{rgb,255:red,230; green,25; blue,75}}
    \tikzColoredSquare{ 2}{\m}{0}{{rgb,255:red,255; green,255; blue,255}}
    \tikzColoredSquare{ 3}{\m}{1}{{rgb,255:red,255; green,255; blue,255}}
    \tikzColoredSquare{ 4}{\m}{5}{{rgb,255:red,255; green,255; blue,255}}
	\def\m{-4}
    \tikzColoredSquare{-1}{\m}{2}{{rgb,255:red,255; green,255; blue,255}}
    \tikzColoredSquare{ 0}{\m}{2}{{rgb,255:red,255; green,255; blue,255}}
    \tikzColoredSquare{ 1}{\m}{4}{{rgb,255:red,230; green,25; blue,75}}
    \tikzColoredSquare{ 2}{\m}{1}{{rgb,255:red,255; green,255; blue,255}}
    \tikzColoredSquare{ 3}{\m}{2}{{rgb,255:red,255; green,255; blue,255}}
    \tikzColoredSquare{ 4}{\m}{6}{{rgb,255:red,255; green,255; blue,255}}
	\def\m{-5}
    \tikzColoredSquare{-1}{\m}{3}{{rgb,255:red,255; green,255; blue,255}}
    \tikzColoredSquare{ 0}{\m}{3}{{rgb,255:red,255; green,255; blue,255}}
    \tikzColoredSquare{ 1}{\m}{5}{{rgb,255:red,230; green,25; blue,75}}
    \tikzColoredSquare{ 2}{\m}{2}{{rgb,255:red,255; green,255; blue,255}}
    \tikzColoredSquare{ 3}{\m}{3}{{rgb,255:red,255; green,255; blue,255}}
    \tikzColoredSquare{ 4}{\m}{5}{{rgb,255:red,230; green,25; blue,75}}
	\def\m{-6}
    \tikzColoredSquare{-1}{\m}{4}{{rgb,255:red,255; green,255; blue,255}}
    \tikzColoredSquare{ 0}{\m}{4}{{rgb,255:red,255; green,255; blue,255}}
    \tikzColoredSquare{ 1}{\m}{0}{{rgb,255:red,255; green,255; blue,255}}
    \tikzColoredSquare{ 2}{\m}{3}{{rgb,255:red,255; green,255; blue,255}}
    \tikzColoredSquare{ 3}{\m}{4}{{rgb,255:red,255; green,255; blue,255}}
    \tikzColoredSquare{ 4}{\m}{6}{{rgb,255:red,230; green,25; blue,75}}
	\def\m{-7}
    \tikzColoredSquare{-1}{\m}{5}{{rgb,255:red,255; green,255; blue,255}}
    \tikzColoredSquare{ 0}{\m}{5}{{rgb,255:red,255; green,255; blue,255}}
    \tikzColoredSquare{ 1}{\m}{1}{{rgb,255:red,255; green,255; blue,255}}
    \tikzColoredSquare{ 2}{\m}{4}{{rgb,255:red,255; green,255; blue,255}}
    \tikzColoredSquare{ 3}{\m}{5}{{rgb,255:red,255; green,255; blue,255}}
    \tikzColoredSquare{ 4}{\m}{5}{{rgb,255:red,255; green,255; blue,255}}
\end{tikzpicture}
\]
\end{ex}
\begin{porism}\label{por:whormsPerColumn}
In the above situation, each \textbf{column} of the super-orbit board of length $\alpha (k+2)$
intersects exactly $\alpha +1$ whorms.  
\end{porism}

\begin{proof}
Let $\mathcal{R}$ be the super-orbit board of $w$ with $\alpha (k+2)$ rows.  Decompose this 
into whorms, with $A=A(f)$ for some $f$ in the orbit. 
Let $\xi_1,\dots,\xi_{\alpha}$ be the consecutive whorms whose tails intersect $f$. 
By definition of \textit{consecutive}, the whorms $\xi_{i}$ and $\xi_{\alpha +i}$ occupy the same columns, 
so the whorms $\xi_1,\xi_{\alpha+1}, \xi_{2\alpha+1},\dots $ intersect this column. 
By Corollary~\ref{cor:Bwhormrepeat}, we have the equalities 
\begin{align*}
\tl(\xi_{2\alpha+1}) &= \tl(\xi_{\alpha}),\\
 \tl(\xi_{3\alpha+1}) &= \tl(\xi_{2\alpha}) = \tl(\xi_{\alpha-1}), \\
 \tl(\xi_{4\alpha+1}) &= \tl(\xi_{3\alpha}) = \tl(\xi_{2\alpha-1}) = \tl(\xi_{\alpha-2}), \\
 &\vdots\\
 \tl(\xi_{\alpha\alpha+1}) &= \tl(\xi_{(\alpha-1)\alpha}) = \tl(\xi_{(\alpha-2)\alpha-1}) =\dots= \tl(\xi_{2}).
\end{align*}
This gives the sum 
$$\tl(\xi_1) +\tl(\xi_{\alpha+1}) + \dots +\tl(\xi_{\alpha\alpha+1})
=\tl(\xi_1) +\tl(\xi_2)+ \dots  +\tl(\xi_{\alpha+1}) =\alpha (k+2)$$
by Lemma~\ref{lem:Bconwhorm}.  This is also the length of the super-orbit board; hence, the
total number of whorms intersecting any given column is $\alpha+1$, since the sum of tail
lengths of $\alpha +1$ worms is the total number of rows in the super-orbit board. 
\end{proof} 

\begin{figure}
\begin{tikzpicture}[baseline={([yshift=-.8ex]current bounding box.center)}]

    \def\m{0}
    \tikzColoredSquare{0}{\m}{3}{white}
    \tikzColoredSquare{1}{\m}{1}{white}
    \tikzColoredSquare{2}{\m}{3}{white}
    \tikzColoredSquare{3}{\m}{3}{white}
    \def\m{-1}
    \tikzColoredSquare{0}{\m}{0}{white}
    \tikzColoredSquare{1}{\m}{2}{white}
    \tikzColoredSquare{2}{\m}{0}{white}
    \tikzColoredSquare{3}{\m}{4}{white}
    \def\m{-2}
    \tikzColoredSquare{0}{\m}{1}{white}
    \tikzColoredSquare{1}{\m}{3}{white}
    \tikzColoredSquare{2}{\m}{1}{white}
    \tikzColoredSquare{3}{\m}{3}{white}
    \def\m{-3}
    \tikzColoredSquare{0}{\m}{2}{white}
    \tikzColoredSquare{1}{\m}{0}{white}
    \tikzColoredSquare{2}{\m}{2}{white}
    \tikzColoredSquare{3}{\m}{4}{white}
    
\end{tikzpicture}
    \caption{ A short orbit board of whirling on $\BB_3\times[4]$. In this example $\alpha = 2$.}\label{fig:shortorbit}
\end{figure}

By Porism~\ref{por:whormsPerColumn}, we will 
have at most $\alpha(\alpha+1)$ total distinct whorms in an orbit board.  (Indeed for equality to be
achieved, the orbit board has to have maximal length $\alpha (k+2)$ before repeating.) 
On the other hand, consider the orbit
board of $\mathcal{F}_4(\BB_3)$ in Figure~\ref{fig:shortorbit} with $\alpha = 2$. 
Here $w^4(f) = f$ so the orbit is only 4 rows long with 2 distinct whorms. 
We can extend this to a super orbit board with 12 rows and 6 whorms (in general with
$\alpha (k+2)$ rows and $\alpha(\alpha+1)$ whorms). 

\begin{thm}\label{thm:Bperiod}

Let $m=\min(k+1,n)$. The order of rowmotion on the chain of claws poset $\mathcal{J}(\BB_n\times [k])$ is $(k+2)\lcm (1,2,\dots ,m)$.
\end{thm}
\noindent Earlier versions of this work reported the weaker statement that the order divides $m!(k+2)$. 

\begin{proof}
By Theorem~\ref{thm:worder}, the whirling order of any $f\in \mathcal{F}_k(\BB_n)$ (i.e.,
the size of the orbit containing $f$) divides
$\alpha (f)(k+2)$, where $\alpha(f) \leq m=\min (k+1,n)$.  Thus the order of the entire map
divides the LCM of these values. 

To prove that this order is exact, we claim that for every $p \in [m]$, there exists an (actual) orbit of length $p(k+2)$, except when $p = k+1$, in which case the orbit has length $k+1$.
Our strategy is similar to the construction used to prove Theorem~\ref{thm:Vorder}. 
However, here we make use of Theorem~\ref{thm:worder}, that $w^{k+2}f=\overline{w}_{A(f)}f$,
i.e., that after $k+2$ steps, $w$ simply cyclically rotates the possible values of $A(f)$. 

Set $f_{1}:=(0,0,0,\dots, 0)$, and for each integer $p\in [2,m]$, set
$$f_p=(0,1,\dots,p-2,p-1, p-1, \dots, p-1,p-1).$$ 
Let $\mathcal{R}_p$ be the orbit board of $w$ containing $f_p$. 
Here $\alpha(f_p) = p$, and the action $\overline w_{A(f_p)}$ partitions $\mathcal{R}_p$
into cycles of length $p$.  
Thus $p$ divides the length of $\mathcal{R}_p$. See Figure~\ref{fig:lcm} for an example of
this construction. 

We will now show the following:
\begin{itemize}
\item If $p<k+1$, then the length of $\mathcal{R}_p$ is $p(k+2)$.
\item If $p=k+1$, then the length of $\mathcal{R}_p$ is $k+1$.
\end{itemize}
In the first case, it suffices to show there is no positive integer $q< k+2$ such that $w^q (f_p) = \overline w_{A(f_p)}f_p$. 
We demonstrate this by showing that $w^q (f_p)(\widehat 0) \neq f_p(\widehat 0)$ for all positive $q < k+2$, noting that $\overline w_{A(f_p)}$ leaves the rightmost label unchanged.

By definition, the rightmost label of $f_p$ is $p-1$, which is less than $k$ by assumption.
Since $\widehat{0}$ is the minimal element of the poset, by definition of whirling its
(rightmost) label will increase by one at 
each step until it reaches the maximum value $k$ at step $k - (p-1)$.  At this step, the
other values in the row (in some order) are $k-1,k-2,\dots , k-p$, which have all increased from their
initial values of $0,1,\dots ,p$ as the constraining label at $\widehat{0}$ has increased.  

From this point onward, each step will involve exactly one (possibly multi-tailed) whorm
moving to the rightmost column with value $k$. There are $p$ whorms, so the rightmost label
remains $k$ for $p$ additional steps, and thus cannot reset to $p-1$ for a total of
$k-(p-1)+p = k+1$ steps.  
Thus, $w^q (f_p) \neq \overline w_{A(f_p)} f_p$ for all $q < k+2$, and 
length of $\mathcal{R}_p$ is $p(k+2)$, finishing this case.

In the second case, we get a collapsing of the suborbit, and it only takes a single
whirl to get us to the $P$-partition with permuted labels.  That is, from the definition of
whirling when $p=k+1$, we get
$$w(f_p) = (1,2,\dots, p-1,0,\dots,0, p-1) = \overline w_{A(f_p)} (f_p). $$
Therefore, the length of $\mathcal{R}_p$ is exactly the order of $\overline w_{A(f)}$, which is $p=k+1$.

The orbit lengths of $\mathcal{R}_1,\dots,\mathcal{R}_m$ depend on the value of $m$. 
If $m = n$, the orbit lengths are $k+2, 2(k+2), \dots, (n-1)(k+2), n(k+2)$. 
If $m = k+1$, the orbit lengths are $k+2, 2(k+2), \dots, k(k+2), k+1$.
In either case, the LCM is $(k+2) \mathrm{LCM}(1,\dots,m)$.
\end{proof}

\begin{figure}
\begin{tikzpicture}[baseline={([yshift=-.8ex]current bounding box.center)}]

    \def\m{0}
    \tikzColoredSquare{0}{\m}{0}{white}
    \tikzColoredSquare{1}{\m}{0}{white}
    \tikzColoredSquare{2}{\m}{0}{white}
    \tikzColoredSquare{3}{\m}{0}{white}
    \tikzColoredSquare{4}{\m}{0}{white}
    \def\m{-1}
    \tikzColoredSquare{0}{\m}{0}{white}
    \tikzColoredSquare{1}{\m}{0}{white}
    \tikzColoredSquare{2}{\m}{0}{white}
    \tikzColoredSquare{3}{\m}{0}{white}
    \tikzColoredSquare{4}{\m}{1}{white}
    \def\m{-2}
    \tikzColoredSquare{0}{\m}{1}{white}
    \tikzColoredSquare{1}{\m}{1}{white}
    \tikzColoredSquare{2}{\m}{1}{white}
    \tikzColoredSquare{3}{\m}{1}{white}
    \tikzColoredSquare{4}{\m}{2}{white}
    \def\m{-3}
    \tikzColoredSquare{0}{\m}{2}{white}
    \tikzColoredSquare{1}{\m}{2}{white}
    \tikzColoredSquare{2}{\m}{2}{white}
    \tikzColoredSquare{3}{\m}{2}{white}
    \tikzColoredSquare{4}{\m}{3}{white}
    \def\m{-4}
    \tikzColoredSquare{0}{\m}{3}{white}
    \tikzColoredSquare{1}{\m}{3}{white}
    \tikzColoredSquare{2}{\m}{3}{white}
    \tikzColoredSquare{3}{\m}{3}{white}
    \tikzColoredSquare{4}{\m}{3}{white}
\end{tikzpicture}
\ \ \ \
\begin{tikzpicture}[baseline={([yshift=-.8ex]current bounding box.center)}]

    \def\m{0}
    \tikzColoredSquare{0}{\m}{0}{white}
    \tikzColoredSquare{1}{\m}{1}{white}
    \tikzColoredSquare{2}{\m}{1}{white}
    \tikzColoredSquare{3}{\m}{1}{white}
    \tikzColoredSquare{4}{\m}{1}{white}
    
    \def\m{-1}
    \tikzColoredSquare{0}{\m}{1}{white}
    \tikzColoredSquare{1}{\m}{0}{white}
    \tikzColoredSquare{2}{\m}{0}{white}
    \tikzColoredSquare{3}{\m}{0}{white}
    \tikzColoredSquare{4}{\m}{2}{white}
    
    \def\m{-2}
    \tikzColoredSquare{0}{\m}{2}{white}
    \tikzColoredSquare{1}{\m}{1}{white}
    \tikzColoredSquare{2}{\m}{1}{white}
    \tikzColoredSquare{3}{\m}{1}{white}
    \tikzColoredSquare{4}{\m}{3}{white}
    
    \def\m{-3}
    \tikzColoredSquare{0}{\m}{3}{white}
    \tikzColoredSquare{1}{\m}{2}{white}
    \tikzColoredSquare{2}{\m}{2}{white}
    \tikzColoredSquare{3}{\m}{2}{white}
    \tikzColoredSquare{4}{\m}{3}{white}
    
    \def\m{-4}
    \tikzColoredSquare{0}{\m}{0}{white}
    \tikzColoredSquare{1}{\m}{3}{white}
    \tikzColoredSquare{2}{\m}{3}{white}
    \tikzColoredSquare{3}{\m}{3}{white}
    \tikzColoredSquare{4}{\m}{3}{white}
    
    \def\m{-5}
    \tikzColoredSquare{0}{\m}{1}{white}
    \tikzColoredSquare{1}{\m}{0}{white}
    \tikzColoredSquare{2}{\m}{0}{white}
    \tikzColoredSquare{3}{\m}{0}{white}
    \tikzColoredSquare{4}{\m}{1}{white}
    
    \def\m{-6}
    \tikzColoredSquare{0}{\m}{0}{white}
    \tikzColoredSquare{1}{\m}{1}{white}
    \tikzColoredSquare{2}{\m}{1}{white}
    \tikzColoredSquare{3}{\m}{1}{white}
    \tikzColoredSquare{4}{\m}{2}{white}
    
    \def\m{-7}
    \tikzColoredSquare{0}{\m}{1}{white}
    \tikzColoredSquare{1}{\m}{2}{white}
    \tikzColoredSquare{2}{\m}{2}{white}
    \tikzColoredSquare{3}{\m}{2}{white}
    \tikzColoredSquare{4}{\m}{3}{white}
    
    \def\m{-8}
    \tikzColoredSquare{0}{\m}{2}{white}
    \tikzColoredSquare{1}{\m}{3}{white}
    \tikzColoredSquare{2}{\m}{3}{white}
    \tikzColoredSquare{3}{\m}{3}{white}
    \tikzColoredSquare{4}{\m}{3}{white}
    
    \def\m{-9}
    \tikzColoredSquare{0}{\m}{3}{white}
    \tikzColoredSquare{1}{\m}{0}{white}
    \tikzColoredSquare{2}{\m}{0}{white}
    \tikzColoredSquare{3}{\m}{0}{white}
    \tikzColoredSquare{4}{\m}{3}{white}
\end{tikzpicture}
\ \ \ \ 
\begin{tikzpicture}[baseline={([yshift=-.8ex]current bounding box.center)}]

    \def\m{0}
    \tikzColoredSquare{0}{\m}{0}{white}
    \tikzColoredSquare{1}{\m}{1}{white}
    \tikzColoredSquare{2}{\m}{2}{white}
    \tikzColoredSquare{3}{\m}{2}{white}
    \tikzColoredSquare{4}{\m}{2}{white}
    
    \def\m{-1}
    \tikzColoredSquare{0}{\m}{1}{white}
    \tikzColoredSquare{1}{\m}{2}{white}
    \tikzColoredSquare{2}{\m}{0}{white}
    \tikzColoredSquare{3}{\m}{0}{white}
    \tikzColoredSquare{4}{\m}{3}{white}
    
    \def\m{-2}
    \tikzColoredSquare{0}{\m}{2}{white}
    \tikzColoredSquare{1}{\m}{3}{white}
    \tikzColoredSquare{2}{\m}{1}{white}
    \tikzColoredSquare{3}{\m}{1}{white}
    \tikzColoredSquare{4}{\m}{3}{white}
    
    \def\m{-3}
    \tikzColoredSquare{0}{\m}{3}{white}
    \tikzColoredSquare{1}{\m}{0}{white}
    \tikzColoredSquare{2}{\m}{2}{white}
    \tikzColoredSquare{3}{\m}{2}{white}
    \tikzColoredSquare{4}{\m}{3}{white}
    \def\m{-4}
    \tikzColoredSquare{0}{\m}{0}{white}
    \tikzColoredSquare{1}{\m}{1}{white}
    \tikzColoredSquare{2}{\m}{3}{white}
    \tikzColoredSquare{3}{\m}{3}{white}
    \tikzColoredSquare{4}{\m}{3}{white}
    \def\m{-5}
    \tikzColoredSquare{0}{\m}{1}{white}
    \tikzColoredSquare{1}{\m}{2}{white}
    \tikzColoredSquare{2}{\m}{0}{white}
    \tikzColoredSquare{3}{\m}{0}{white}
    \tikzColoredSquare{4}{\m}{2}{white}
    \def\m{-6}
    \tikzColoredSquare{0}{\m}{2}{white}
    \tikzColoredSquare{1}{\m}{0}{white}
    \tikzColoredSquare{2}{\m}{1}{white}
    \tikzColoredSquare{3}{\m}{1}{white}
    \tikzColoredSquare{4}{\m}{3}{white}
    \def\m{-7}
    \tikzColoredSquare{0}{\m}{3}{white}
    \tikzColoredSquare{1}{\m}{1}{white}
    \tikzColoredSquare{2}{\m}{2}{white}
    \tikzColoredSquare{3}{\m}{2}{white}
    \tikzColoredSquare{4}{\m}{3}{white}
    \def\m{-8}
    \tikzColoredSquare{0}{\m}{0}{white}
    \tikzColoredSquare{1}{\m}{2}{white}
    \tikzColoredSquare{2}{\m}{3}{white}
    \tikzColoredSquare{3}{\m}{3}{white}
    \tikzColoredSquare{4}{\m}{3}{white}
    \def\m{-9}
    \tikzColoredSquare{0}{\m}{1}{white}
    \tikzColoredSquare{1}{\m}{3}{white}
    \tikzColoredSquare{2}{\m}{0}{white}
    \tikzColoredSquare{3}{\m}{0}{white}
    \tikzColoredSquare{4}{\m}{3}{white}
    \def\m{-10}
    \tikzColoredSquare{0}{\m}{2}{white}
    \tikzColoredSquare{1}{\m}{0}{white}
    \tikzColoredSquare{2}{\m}{1}{white}
    \tikzColoredSquare{3}{\m}{1}{white}
    \tikzColoredSquare{4}{\m}{2}{white}
    \def\m{-11}
    \tikzColoredSquare{0}{\m}{0}{white}
    \tikzColoredSquare{1}{\m}{1}{white}
    \tikzColoredSquare{2}{\m}{2}{white}
    \tikzColoredSquare{3}{\m}{2}{white}
    \tikzColoredSquare{4}{\m}{3}{white}
    \def\m{-12}
    \tikzColoredSquare{0}{\m}{1}{white}
    \tikzColoredSquare{1}{\m}{2}{white}
    \tikzColoredSquare{2}{\m}{3}{white}
    \tikzColoredSquare{3}{\m}{3}{white}
    \tikzColoredSquare{4}{\m}{3}{white}
    \def\m{-13}
    \tikzColoredSquare{0}{\m}{2}{white}
    \tikzColoredSquare{1}{\m}{3}{white}
    \tikzColoredSquare{2}{\m}{0}{white}
    \tikzColoredSquare{3}{\m}{0}{white}
    \tikzColoredSquare{4}{\m}{3}{white}
    \def\m{-14}
    \tikzColoredSquare{0}{\m}{3}{white}
    \tikzColoredSquare{1}{\m}{0}{white}
    \tikzColoredSquare{2}{\m}{1}{white}
    \tikzColoredSquare{3}{\m}{1}{white}
    \tikzColoredSquare{4}{\m}{3}{white}
\end{tikzpicture}
\ \ \ \
\begin{tikzpicture}[baseline={([yshift=-.8ex]current bounding box.center)}]

    \def\m{0}
    \tikzColoredSquare{0}{\m}{0}{white}
    \tikzColoredSquare{1}{\m}{1}{white}
    \tikzColoredSquare{2}{\m}{2}{white}
    \tikzColoredSquare{3}{\m}{3}{white}
    \tikzColoredSquare{4}{\m}{3}{white}
    
    \def\m{-1}
    \tikzColoredSquare{0}{\m}{1}{white}
    \tikzColoredSquare{1}{\m}{2}{white}
    \tikzColoredSquare{2}{\m}{3}{white}
    \tikzColoredSquare{3}{\m}{0}{white}
    \tikzColoredSquare{4}{\m}{3}{white}
    
    \def\m{-2}
    \tikzColoredSquare{0}{\m}{2}{white}
    \tikzColoredSquare{1}{\m}{3}{white}
    \tikzColoredSquare{2}{\m}{0}{white}
    \tikzColoredSquare{3}{\m}{1}{white}
    \tikzColoredSquare{4}{\m}{3}{white}
    
    \def\m{-3}
    \tikzColoredSquare{0}{\m}{3}{white}
    \tikzColoredSquare{1}{\m}{0}{white}
    \tikzColoredSquare{2}{\m}{1}{white}
    \tikzColoredSquare{3}{\m}{2}{white}
    \tikzColoredSquare{4}{\m}{3}{white}
\end{tikzpicture}
    \caption{Four orbit boards of whirling on $\BB_4\times[3]$ with $\alpha = 1,2,3,$ and
    $4$ respectively. The order of whirling here is $60= \lcm(1,2,3,4)\cdot 5$.}\label{fig:lcm}
\end{figure}

The following theorem is the analogue 
of the first homomesy in Theorem~\ref{thm:Vhom}; it follows directly from
Theorem~\ref{thm:worder} just as in the proof of Theorem~\ref{thm:Vhom}(1) for the case $n=2$.

\begin{thm}\label{thm:homet} Let $\chi_{(i,a)}$ denote the indicator function for $(b_{i},a)\in
\BB_n\times [k]$. Then for the action of rowmotion on $\mathcal{J}(\BB_n\times [k])$, 
the statistic $\chi_{(i,a)}-\chi_{(j,a)}$ is $0$-mesic for all $i,j\in[n]$ and $a\in[k]$. 
\end{thm}

\begin{remark}
The second homomesy in Theorem~\ref{thm:Vhom} fails to hold in general for $n>2$. 
The average of the statistic $\left(\sum_{i=1}^n\chi_{(i,1)}\right)-\chi_{(\widehat{0},k)}$
(analogous to Theorem~\ref{thm:Vhom}(2))
turns out to be \textbf{dependent} on $\alpha(f)$ (for any $f\in \mathcal{R}$) and can be
computed as 
\begin{equation}\label{eq:alpha-dependent}
\frac{n(\alpha) (k+2) - (n+\alpha)(\alpha + 1)}{\alpha(k+2)}.
\end{equation}
To see this consider the super orbit of length $\alpha (k+2)$ for some $f\in
\mathcal{F}_k(\BB_n)$ with $n\alpha(k+2)$ entries in the first $n$ columns. 
We know $\chi_{(i,1)}(I)=0$ if and only if the $f$ corresponding to $I$ satisfies $f(b_{i})
= 0$.  (See the proof of Lemma~\ref{lem:eqwhirl}.) But this is counted by the number of
whorm beginnings (considering each tail to be separate); by Porism~\ref{por:whormsPerColumn} this is $n(\alpha+1)$. 
Furthermore, $\chi_{(\widehat{0},k)}(I)=1$ if and
only if the $f$ corresponding to $I$ satisfies $f(\widehat{0}) = k$, which is counted by the
number of whorm endings, that is $\alpha(\alpha+1)$. Thus, the numerator counts 1 for
each of the $n\alpha (k+2) - n(\alpha +1)$ entries that are not at the start of a tail, and
$-1$ for the $\alpha (\alpha +1)$ whorm heads in the rightmost column.  Divide by the size
of the super-orbit board to get the average. In the $n=2$ case the only possibilities are
$\alpha=1$ or $\alpha =n$, for both of which (\ref{eq:alpha-dependent}) reduces to
$\frac{nk-2}{k+2} = \frac{2k-2}{k+2}$; this explains why we get the homomesy
Theorem~\ref{thm:Vhom}(2) for $\VV \times [k]$. 
\end{remark}

\begin{figure}
\begin{tikzpicture}[baseline={([yshift=-.8ex]current bounding box.center)}]

    \def\m{0}
    \tikzColoredSquare{0}{\m}{1}{white}
    \tikzColoredSquare{1}{\m}{3}{white}
    \tikzColoredSquare{2}{\m}{5}{white}
    \tikzColoredSquare{3}{\m}{5}{white}
    \tikzColoredSquare{4}{\m}{5}{white}
    
    \def\m{-1}
    \tikzColoredSquare{0}{\m}{2}{white}
    \tikzColoredSquare{1}{\m}{4}{white}
    \tikzColoredSquare{2}{\m}{0}{white}
    \tikzColoredSquare{3}{\m}{0}{white}
    \tikzColoredSquare{4}{\m}{6}{white}
    
    \def\m{-2}
    \tikzColoredSquare{0}{\m}{3}{white}
    \tikzColoredSquare{1}{\m}{5}{white}
    \tikzColoredSquare{2}{\m}{1}{white}
    \tikzColoredSquare{3}{\m}{1}{white}
    \tikzColoredSquare{4}{\m}{5}{white}
    
    \def\m{-3}
    \tikzColoredSquare{0}{\m}{4}{white}
    \tikzColoredSquare{1}{\m}{0}{white}
    \tikzColoredSquare{2}{\m}{2}{white}
    \tikzColoredSquare{3}{\m}{2}{white}
    \tikzColoredSquare{4}{\m}{6}{white}
    
    \def\m{-4}
    \tikzColoredSquare{0}{\m}{5}{white}
    \tikzColoredSquare{1}{\m}{1}{white}
    \tikzColoredSquare{2}{\m}{3}{white}
    \tikzColoredSquare{3}{\m}{3}{white}
    \tikzColoredSquare{4}{\m}{5}{white}
    
    \def\m{-5}
    \tikzColoredSquare{0}{\m}{0}{white}
    \tikzColoredSquare{1}{\m}{2}{white}
    \tikzColoredSquare{2}{\m}{4}{white}
    \tikzColoredSquare{3}{\m}{4}{white}
    \tikzColoredSquare{4}{\m}{6}{white}
\end{tikzpicture}
    \caption{An orbit board of whirling on $\BB_4\times[6]$. The average of the statistic $B_3-B_2$ is  $2/3$ which does not agree with the statement of the generalized flux capacitor in Remark~\ref{rem:Claw_Flux}. Notice if the labels in repeated column are ignored, then the average of the statistic $B_3-B_2$ is $1/2$ which agrees with the generalization.}\label{fig:Claw_Flux_Counter}
\end{figure}
\begin{remark}\label{rem:Claw_Flux}
The ``flux-capacitor'' homomesy of Theorem~\ref{thm:flux} also fails to generalize to the claw-graph
setting. Earlier versions of this work incorrectly asserted the following analogous
``flux-capacitor'' homomesy: 
``Let $B_i = \chi_{(i-1,\widehat{0})}+ \sum_{\ell=1}^n \chi_{(i,\ell)}$.  Then for the action of rowmotion
on $\mathcal{J}(\BB_n\times [k])$,  $B_{i}-B_{j}$ is $\frac{(j-i)(n+1)}{k+2}$-mesic  for all
$i,j\in [n]$.'' \\
\noindent
Unfortunately we only obtain this average in general for the subset of orbit boards which
are tiled entirely by one-tailed whorms. See Figure~\ref{fig:Claw_Flux_Counter} for a counterexample. 
\end{remark}

\subsection{Acknowledgments}\label{ss:ack}
We are grateful to Banff International Research Station for hosting 
online and hybrid workshops on dynamical algebraic combinatorics in October 2020 and
November 2021. Discussions with the following colleagues were particularly helpful:
Ben Adenbaum,
Chinmay Dharmendra,
Sam Hopkins, 
Michael Joseph, 
James Propp, 
Bruce Sagan, 
and Jessica Striker. We also thank two anonymous reviewers of an earlier (extended abstract)
version of this work. 

\printbibliography

\end{document}